\documentclass[11pt]{amsart}
\usepackage{mathrsfs}
\usepackage{amsmath}
\usepackage{amssymb}
\usepackage{graphicx}
\newtheorem{theorem}{Theorem}[section]
\newtheorem{corollary}[theorem]{Corollary}
\newtheorem{lemma}[theorem]{Lemma}
\newtheorem{proposition}[theorem]{Proposition}
\theoremstyle{definition}
\newtheorem{definition}[theorem]{Definition}

\numberwithin{equation}{section}
\title[On Scaled Hyperbolic Numbers Induced by Scaled Hypercomplex Rings]{On Scaled Hyperbolic Numbers Induced by Scaled Hypercomplex Rings}

\author[D. Alpay]{Daniel Alpay}

\address[D. Alpay]{Chapman Univ., Dept. of Math., 1 University Dr., Orange, CA, 92866,  U. S. A. }
\email{\tt alpay@chapman.edu}

\author[I. Cho]{Ilwoo Cho}
\address[I. Cho]{St. Ambrose Univ., Dept. of Math. and Stat., 421 Ambrose
Hall, 518 W. Locust St., Davenport, Iowa, 52803, U. S. A.}
\email{\tt choilwoo@sau.edu}



\begin{document}

\begin{abstract}
In this paper, we generalize the well-known hyperbolic numbers to
certain numeric structures scaled by the real numbers. Under our scaling
of $\mathbb{R}$, the usual hyperbolic numbers are understood to be
our $1$-scaled hyperbolic numbers. If a scale $t$ is not positive
in $\mathbb{R}$, then our $t$-scaled hyperbolic numbers have similar
numerical structures with those of the complex numbers of $\mathbb{C}$,
however, if a scale is positive in $\mathbb{R}$, then their numerical
properties are similar to those of the classical hyperbolic numbers.
We here understand scaled-hyperbolic numbers as elements of the scaled-hypercomplex
rings $\left\{ \mathbb{H}_{t}\right\} _{t\in\mathbb{R}}$, introduced
in {[}1{]}. This scaled-hyperbolic analysis is done by algebra, analysis,
operator theory, operator-algebra theory and free probability on scaled-hypercomplex
numbers. 
\end{abstract}

\maketitle

\section{Introduction}

In this paper, we generalize the well-known hyperbolic numbers (or,
split-complex numbers) by applying the algebra, analysis and representation
theory on the scaled hypercomplex structures recently introduced in
{[}1{]}. By doing so, one can generalize the hyperbolic analysis and
geometry in connection with the study of the scaled-hypercomplex numbers.

Throughout this paper, we let $\mathbb{C}^{2}$ be the usual 2-dimensional
vector space over the complex field $\mathbb{C}$. Every vector $\left(a,b\right)\in\mathbb{C}^{2}$
is understood as a hypercomplex number induced by the complex numbers
$a$ and $b$. Under a suitable scaling over the real numbers $\mathbb{R}$,
the set $\mathbb{C}^{2}$ of hypercomplex numbers forms a ring, 
\[
\mathbb{H}_{t}=\left(\mathbb{C}^{2},\;+,\;\cdot_{t}\right),
\]
where ($+$) is the usual vector addition on $\mathbb{C}^{2}$, and
($\cdot_{t}$) is the $t$-scaled vector-multiplication, 
\[
\left(a_{1},b_{1}\right)\cdot_{t}\left(a_{2},b_{2}\right)=\left(a_{1}a_{2}+tb_{1}\overline{b_{2}},\:a_{1}b_{2}+b_{1}\overline{a_{2}}\right),
\]
on $\mathbb{C}^{2}$, where $\overline{z}$ are the conjugates of
$z$ in $\mathbb{C}$. By the representation $\left(\mathbb{C}^{2},\pi_{t}\right)$
of the $t$-scaled hypercomplex ring $\mathbb{H}_{t}$, considered
in {[}1{]}, we regard a hypercomplex number $h=\left(a,b\right)\in\mathbb{H}_{t}$
as a (2$\times$2)-matrix, or a Hilbert-space operator acting on $\mathbb{C}^{2}$,
\[
\pi_{t}\left(h\right)\overset{\textrm{denote}}{=}[h]_{t}\overset{\textrm{def}}{=}\left(\begin{array}{cc}
a & tb\\
\overline{b} & \overline{a}
\end{array}\right)\;\mathrm{in\;}M_{2}\left(\mathbb{C}\right),
\]
where $M_{2}\left(\mathbb{C}\right)=B\left(\mathbb{C}^{2}\right)$
is the matricial algebra (or, the operator $C^{*}$-algebra acting
on the Hilbert space $\mathbb{C}^{2}$) over $\mathbb{C}$, for $t\in\mathbb{R}$.
Remark that the ring $\mathbb{H}_{-1}$ is nothing but the noncommutative
field $\mathbb{H}$ of all quaternions (e.g., {[}3{]} and {[}19{]}),
and the ring $\mathbb{H}_{1}$ is the ring of all split-quaternions
(e.g., {[}1{]}). The algebra, spectra theory, operator theory, and
free probability on $\left\{ \mathbb{H}_{t}\right\} _{t\in\mathbb{R}}$
are studied in {[}1{]}, under representation. We, however, have to
emphasize that the operator theory and free probability of {[}1{]},
and those of this paper are different; in {[}1{]}, we construct a
$C^{*}$-algebra $\mathfrak{H}_{2}^{t}$ generated by $\mathbb{H}_{t}$
in $M_{2}\left(\mathbb{C}\right)$ somewhat artificially, and then
study operator theory and free probability on $\mathfrak{H}_{2}^{t}$,
meanwhile, here, we find a natural complete semi-normed $*$-algebraic
structure of $\mathbb{H}_{t}$, and then investigate operator theory
and the corresponding analytic data directly on $\mathbb{H}_{t}$.

\subsection{Motivation}

The quaternions $\mathbb{H}=\mathbb{H}_{-1}$ has been studied in
various different fields in mathematics and applied science (e.g.,
{[}5{]}, {[}7{]}, {[}10{]}, {[}11{]}, {[}12{]}, {[}16{]}, {[}17{]},
{[}19{]}, {[}20{]}, {[}21{]}, {[}23{]}), as an extended algebraic
structure of the complex field $\mathbb{C}$, which also motivates
the construction and analysis on Clifford algebras (e.g., {[}4{]}).
From the quaternions $\mathbb{H}=\mathbb{H}_{-1}$, we extended them
to those on the scaled-hypercomplex rings $\left\{ \mathbb{H}_{t}\right\} _{t\in\mathbb{R}}$
in {[}1{]}, and generalized the main results of {[}3{]}.

Meanwhile, the hyperbolic numbers are considered in non-Euclidean
geometry, especially in hyperbolic geometry (e.g., {[}14{]} and {[}15{]}).
We cannot help emphasizing the importance of the hyperbolic analysis
and their applications to various applied areas (e.g., {[}13{]}).

The main purpose of this paper is to study generalized hyperbolic
algebra and analysis, and the corresponding operator theory under
representation, with help of those on scaled-hypercomplex rings. The
main tools and techniques are from the construction of the subring,
\[
\mathbb{D}_{t}=\left\{ \left(x,y\right)\in\mathbb{H}_{t}:x,y\in\mathbb{R}\right\} ,
\]
inside $\mathbb{H}_{t}$. It forms the collection of all ``scaled-hyperbolic''
numbers for all $t\in\mathbb{R}$.

\subsection{Overview}

In Section 2, we review the definitions and the main results of {[}1{]}. 

In Sections 3 and 4, we define and study the hypercomplex-conjugate
on $\mathbb{H}_{t}$ to understand $\mathbb{H}_{t}$ as a topological
vector space, for each $t\in\mathbb{R}$. It is shown that, by the
hypercomplex conjugation, one can define a bilinear form $\left\langle ,\right\rangle _{t}$
on $\mathbb{H}_{t}$ over $\mathbb{R}$, and it induces a semi-norm
$\left\Vert .\right\Vert _{t}$ on $\mathbb{H}_{t}$ over $\mathbb{R}$,
i.e., the scaled hypercomplex ring $\mathbb{H}_{t}$ is regarded as
a semi-normed vector space over $\mathbb{R}$, for $t\in\mathbb{R}$.
In particular, we show that if $t<0$, then $\mathbb{H}_{t}$ forms
a Hilbert space over $\mathbb{R}$, meanwhile, if $t\geq0$, then
it becomes a complete indefinite semi-inner product space over $\mathbb{R}$. 

In Section 5, by understanding $\mathbb{H}_{t}$ as a topological
vector space over $\mathbb{R}$, the multiplication operators acting
on $\mathbb{H}_{t}$ with their symbols in $\mathbb{H}_{t}$ are considered,
for each fixed $t\in\mathbb{R}$. It is shown that these operators
form a well-defined operator algebra over $\mathbb{R}$ (acting on
$\mathbb{H}_{t}$), i.e., $\mathbb{H}_{t}$ is understood to be an
operator algebra of multiplication operators acting on $\mathbb{H}_{t}$
over $\mathbb{R}$, for all $t\in\mathbb{R}$. Some interesting operator-theoretic
properties; self-adjointness, projection-property, normality, isometric
property and unitarity; and the corresponding analytic data are considered.

After considering the algebraic (as rings), analytic (as topological
vector spaces), and operator-theoretic (as multiplication-operator
algebras) of $\left\{ \mathbb{H}_{t}\right\} _{t\in\mathbb{R}}$ in
Sections 3, 4 and 5, we construct so-called the $t$-scaled hyperbolic
numbers $\mathbb{D}_{t}$ as a sub-structure of $\mathbb{H}_{t}$,
for $t\in\mathbb{R}$. It means that $\mathbb{D}_{t}$ is a subring,
a subspace, and a subalgebra of $\mathbb{H}_{t}$, for all $t\in\mathbb{R}$.
By definition, $\mathbb{D}_{-1}$ is nothing but the complex field
$\mathbb{C}$, and $\mathbb{D}_{1}$ is the classical hyperbolic numbers
$\mathbb{D}$. In Section 7, the polar decomposition of $\left\{ \mathbb{D}_{t}\right\} _{t\in\mathbb{R}}$
is characterized up to scales.

In Section 8, we act the matricial algebra $M_{2}\left(\mathbb{R}\right)$
(over $\mathbb{R}$) on $\mathbb{D}_{t}$ (as a vector space), for
$t\in\mathbb{R}$. Such an action provides a way how one may / can
consider geometry on $\mathbb{D}_{t}$ in future.

\section{The Hypercomplex Numbers Scaled by $\mathbb{R}$}

In this section, we define a ring $\mathbb{H}_{t}$ of hypercomplex
numbers for a fixed scale $t\in\mathbb{R}$. Throughout this section,
we let 
\[
\mathbb{C}^{2}=\left\{ \left(a,b\right):a,b\in\mathbb{C}\right\} 
\]
be the the usual 2-dimensional Hilbert space equipped with its canonical
orthonormal basis, $\left\{ \left(1,0\right),\;\left(0,1\right)\right\} .$

\subsection{Scaled Hypercomplex Rings}

Fix an arbitrarily scale $t$ in the real field $\mathbb{R}$. On
the Hilbert space $\mathbb{C}^{2}$, define the $t$-scaled vector-multiplication
($\cdot_{t}$) by

\medskip{}

\hfill{}$\left(a_{1},b_{1}\right)\cdot_{t}\left(a_{2},b_{2}\right)\overset{\textrm{def}}{=}\left(a_{1}a_{2}+tb_{1}\overline{b_{2}},\:a_{1}b_{2}+b_{1}\overline{a_{2}}\right),$\hfill{}(2.1.1)

\medskip{}

\noindent for $\left(a_{l},b_{l}\right)\in\mathbb{C}^{2}$, for all
$l=1,2$. 
\begin{proposition}
The algebraic structure $\left(\mathbb{C}^{2},+,\cdot_{t}\right)$
forms a unital ring with its unity, or the ($\cdot_{t}$)-identity,
$\left(1,0\right)$, where ($+$) is the usual vector addition on
$\mathbb{C}^{2}$, and ($\cdot_{t}$) is the vector multiplication
(2.1.1). 
\end{proposition}

\begin{proof}
The pair $\left(\mathbb{C}^{2},+\right)$ is an abelian group for
($+$) with its ($+$)-identity $\left(0,0\right)$. And the algebraic
pair $\left(\mathbb{C}^{2\times},\cdot_{t}\right)$ is a semigroup
with its ($\cdot_{t}$)-identity $\left(1,0\right)$ where $\mathbb{C}^{2\times}=\mathbb{C}^{2}\setminus\left\{ \left(0,0\right)\right\} $.
It is not difficult to check ($+$) and ($\cdot_{t}$) are distributed
on $\mathbb{C}^{2}$ (e.g., see {[}1{]} for details). So, the algebraic
triple $\left(\mathbb{C}^{2},+,\cdot_{t}\right)$ forms a unital ring
with its unity $\left(1,0\right)$. 
\end{proof}
Since $\mathbb{C}^{2}$ is a Hilbert space equipped with the usual-metric
topology, one can understand these unital rings $\left\{ \left(\mathbb{C}^{2},+,\cdot_{t}\right)\right\} _{t\in\mathbb{R}}$
as topological rings. 
\begin{definition}
For $t\in\mathbb{R}$, the ring $\mathbb{H}_{t}\overset{\textrm{denote}}{=}\left(\mathbb{C}^{2},+,\cdot_{t}\right)$
is called the $t$-scaled hypercomplex ring. 
\end{definition}

For a fixed $t\in\mathbb{R}$, let $\mathbb{H}_{t}$ be the $t$-scaled
hypercomplex ring. Define an injective map, 
\[
\pi_{t}:\mathbb{H}_{t}\rightarrow B\left(\mathbb{C}^{2}\right)=M_{2}\left(\mathbb{C}\right),
\]
by\hfill{}(2.1.2) 
\[
\pi_{t}\left(\left(a,b\right)\right)=\left(\begin{array}{cc}
a & tb\\
\overline{b} & \overline{a}
\end{array}\right),\;\forall\left(a,b\right)\in\mathbb{H}_{t},
\]
where $B\left(H\right)$ is the operator algebra, the $C^{*}$-algebra
acting on a Hilbert space $H$ equipped with the operator-norm topology,
consisting of all (bounded linear) operators on $H$, and $M_{k}\left(\mathbb{C}\right)$
is the matricial algebra of all $\left(k\times k\right)$-matrices
over $\mathbb{C}$, which is $*$-isomorphic to $B\left(\mathbb{C}^{k}\right)$,
for all $k\in\mathbb{N}$ (e.g., {[}8{]} and {[}9{]}). Such an injective
map $\pi_{t}$ satisfies that 
\[
\pi_{t}\left(h_{1}+h_{2}\right)=\pi_{t}\left(h_{1}\right)+\pi_{t}\left(h_{2}\right),
\]
and\hfill{}(2.1.3) 
\[
\pi_{t}\left(h_{1}\cdot_{t}h_{2}\right)=\pi_{t}\left(h_{1}\right)\pi_{t}\left(h_{2}\right),
\]
in $M_{2}\left(\mathbb{C}\right)$ (e.g., see {[}1{]} for details). 
\begin{proposition}
The pair $\left(\mathbb{C}^{2},\:\pi_{t}\right)$ forms an injective
Hilbert-space representation of our $t$-scaled hypercomplex ring
$\mathbb{H}_{t}$, where $\pi_{t}$ is an action (2.1.2). 
\end{proposition}

\begin{proof}
The injective function $\pi_{t}$ of (2.1.2) is indeed a well-defined
ring-action by (2.1.3). Since $\mathbb{C}^{2}$ and $M_{2}\left(\mathbb{C}\right)$
are finite-dimensional, the ring-action $\pi_{t}$ is continuous,
i.e., $\pi_{t}$ is a well-defined topological-ring-action of $\mathbb{H}_{t}$
acting $\mathbb{C}^{2}$. 
\end{proof}
By the above proposition, the realization,

\medskip{}

\hfill{}$\pi_{t}\left(\mathbb{H}_{t}\right)=\left\{ \left(\begin{array}{cc}
a & tb\\
\overline{b} & \overline{a}
\end{array}\right)\in M_{2}\left(\mathbb{C}\right):\left(a,b\right)\in\mathbb{H}_{t}\right\} ,$\hfill{}(2.1.4)

\medskip{}

\noindent of $\mathbb{H}_{t}$ is well-determined in $M_{2}\left(\mathbb{C}\right)$. 
\begin{definition}
The realization $\mathcal{H}_{2}^{t}\overset{\textrm{denote}}{=}\pi_{t}\left(\mathbb{H}_{t}\right)$
of (2.1.4) is called the $t$-scaled (hypercomplex-)realization of
$\mathbb{H}_{t}$ (in $M_{2}\left(\mathbb{C}\right)$) for $t\in\mathbb{R}$.
For convenience, we denote the realization $\pi_{t}\left(h\right)$
of $h\in\mathbb{H}_{t}$ by $\left[h\right]_{t}$ in $\mathcal{H}_{2}^{t}$. 
\end{definition}

By the above definition, we obtain the following algebraic characterization. 
\begin{proposition}
Let $t\in\mathbb{R}$ be arbitrary. Then

\medskip{}

\hfill{}$\mathbb{H}_{t}\overset{\textrm{T.R}}{=}\mathit{\mathcal{H}_{2}^{t}\;\;\;\mathrm{in\;\;\;}M_{2}\left(\mathbb{C}\right)},$\hfill{}(2.1.5)

\medskip{}

\noindent where ``$\overset{\textrm{T.R}}{=}$'' means ``being
topological-ring-isomorphic to.'' 
\end{proposition}

\begin{proof}
The relation (2.1.5) is proven by (2.1.3) and (2.1.4), by the injectivity
of $\pi_{t}$. 
\end{proof}
If $\mathbb{H}_{t}^{\times}\overset{\textrm{denote}}{=}\mathbb{H}_{t}\setminus\left\{ \left(0,0\right)\right\} ,$
where $\left(0,0\right)\in\mathbb{H}_{t}$ is the ($+$)-identity,
then, this set $\mathbb{H}_{t}^{\times}$ forms the maximal monoid,
\[
\mathbb{H}_{t}^{\times}\overset{\textrm{denote}}{=}\left(\mathbb{H}_{t}^{\times},\;\cdot_{t}\right),
\]
embedded in the ring $\mathbb{H}_{t}$, with its identity $\left(1,0\right)$
(e.g., see {[}1{]}). 
\begin{definition}
The maximal embedded monoid $\mathbb{H}_{t}^{\times}=\left(\mathbb{H}_{t}^{\times},\:\cdot_{t}\right)$
of $\mathbb{H}_{t}$ is called the $t$-scaled hypercomplex monoid. 
\end{definition}

By (2.1.5), the $t$-scaled hypercomplex monoid $\mathbb{H}_{t}^{\times}$
is monoid-isomorphic to $\mathcal{H}_{2}^{t\times}\overset{\textrm{denote}}{=}\left(\mathcal{H}_{2}^{t\times},\:\cdot\right)$
with its identity, $I_{2}=\left[\left(1,0\right)\right]_{t},$ the
$\left(2\times2\right)$-identity matrix of $M_{2}\left(\mathbb{C}\right)$,
where ($\cdot$) is the matricial multiplication, i.e.,

\medskip{}

\hfill{}$\mathbb{H}_{t}^{\times}=\left(\mathbb{H}_{t}^{\times},\:\cdot_{t}\right)\overset{\textrm{Monoid}}{=}\left(\mathcal{H}_{2}^{t\times},\:\cdot\right)=\mathcal{H}_{2}^{t\times},$\hfill{}(2.1.6)

\medskip{}

\noindent where ``$\overset{\textrm{Monoid}}{=}$'' means ``being
monoid-isomorphic.''

\subsection{Algebra on $\mathbb{H}_{t}$}

For an arbitrarily fixed $t\in\mathbb{R}$, let $\mathbb{H}_{t}$
be the corresponding $t$-scaled hypercomplex ring, isomorphic to
its $t$-scaled realization $\mathcal{H}_{2}^{t}$ by (2.1.5). Observe
that, for any $\left(a,b\right)\in\mathbb{H}_{t}$, one has

\medskip{}

\hfill{}$det\left(\left[\left(a,b\right)\right]_{t}\right)=det\left(\begin{array}{cc}
a & tb\\
\overline{b} & \overline{a}
\end{array}\right)=\left|a\right|^{2}-t\left|b\right|^{2}.$\hfill{}(2.2.1)

\medskip{}

\noindent where $det$ is the determinant, and $\left|.\right|$ is
the modulus on $\mathbb{C}$. 
\begin{lemma}
If $\left(a,b\right)\in\mathbb{H}_{t}$, then $\left|a\right|^{2}\neq t\left|b\right|^{2}$
in $\mathbb{C}$, if and only if $\left(a,b\right)$ is invertible
in $\mathbb{H}_{t}$ with its inverse,

\noindent 
\[
\left(a,b\right)^{-1}=\left(\frac{\overline{a}}{\left|a\right|^{2}-t\left|b\right|^{2}},\:\frac{-b}{\left|a\right|^{2}-t\left|b\right|^{2}}\right)\;\;\mathrm{in\;\;}\mathbb{H}_{t},
\]
satisfying\hfill{}(2.2.2) 
\[
\left[\left(a,b\right)^{-1}\right]_{t}=\left[\left(a,b\right)\right]_{t}^{-1}\;\;\mathrm{in\;\;\mathcal{H}_{2}^{t}.}
\]
\end{lemma}

\begin{proof}
Suppose $\left(a,b\right)\in\mathbb{H}_{t}$ satisfies the condition
$\left|a\right|^{2}\neq t\left|b\right|^{2}$ in $\mathbb{C}$, for
$t\in\mathbb{R}$. Then, by (2.2.1), 
\[
det\left(\left[\left(a,b\right)\right]_{t}\right)=\left|a\right|^{2}-t\left|b\right|^{2}\neq0.
\]
Thus, the matrix $\left[\left(a,b\right)\right]_{t}$ is invertible
``in $M_{2}\left(\mathbb{C}\right)$,'' with its inverse matrix,
\[
\left[\left(a,b\right)\right]_{t}^{-1}=\frac{1}{\left|a\right|^{2}-t\left|b\right|^{2}}\left(\begin{array}{cc}
\overline{a} & -tb\\
-\overline{b} & a
\end{array}\right)\;\mathrm{in\;}M_{2}\left(\mathbb{C}\right).
\]
Note that this inverse satisfies that 
\[
\left[\left(a,b\right)\right]_{t}^{-1}=\left[\left(\frac{\overline{a}}{\left|a\right|^{2}-t\left|b\right|^{2}},\:\frac{-b}{\left|a\right|^{2}-t\left|b\right|^{2}}\right)\right]_{t},
\]
''in the $t$-scaled realization $\mathcal{H}_{2}^{t}$.'' Therefore,
the relation (2.2.2) holds. 
\end{proof}
Recall that an algebraic triple $\left(X,+,\cdot\right)$ is a noncommutative
field, if it is a unital ring, where $\left(X^{\times},\cdot\right)$
is a non-abelian group (e.g., {[}1{]} and {[}3{]}). 
\begin{theorem}
We have the algebraic characterization,

\medskip{}

\hfill{}$t<0\;\mathrm{in\;}\mathbb{R}\Longleftrightarrow\mathbb{H}_{t}\textrm{ is a noncommutative field.}$\hfill{}(2.2.3) 
\end{theorem}

\begin{proof}
($\Rightarrow$) By the above theorem, if $t<0$ in $\mathbb{R}$,
then every hypercomplex number $\left(a,b\right)$ of the $t$-scaled
hypercomplex monoid $\mathbb{H}_{t}^{\times}$ automatically satisfies
the condition (2.2.2): $\left|a\right|^{2}\neq t\left|b\right|^{2}$,
because 
\[
\left|a\right|^{2}>t\left|b\right|^{2}\Longrightarrow\left|a\right|^{2}\neq t\left|b\right|^{2},
\]
since left-hand side is nonnegative, and the right-hand side is not
positive, and $\left(a,b\right)\neq\left(0,0\right)$. Thus, if $t<0$,
then every monoidal element $h\in\mathbb{H}_{t}^{\times}$ is invertible
in $\mathbb{H}_{t}$, equivalently, the monoid $\mathbb{H}_{t}^{\times}$
is a group. By (2.1.1), this group $\mathbb{H}_{t}^{\times}$ is noncommutative
up to ($\cdot_{t}$). Therefore, if $t<0$, then the ring $\mathbb{H}_{t}$
becomes a noncommutative field.

\noindent ($\Leftarrow$) Assume that $t\geq0$. First, let $t=0$.
If $\left(0,b\right)\in\mathbb{H}_{0}^{\times}$ (i.e., $b\neq0$),
then 
\[
det\left(\left[\left(0,b\right)\right]_{0}\right)=det\left(\left(\begin{array}{cc}
0 & 0\\
\overline{b} & 0
\end{array}\right)\right)=0,
\]
implying that $\left[\left(0,b\right)\right]_{0}\in\mathcal{H}_{2}^{t}$
is not invertible. So, the monoid $\mathbb{H}_{0}^{\times}$ is not
a group. Now, let $t>0$. If $\left(a,b\right)\in\mathbb{H}_{t}^{\times}$,
with $\left|b\right|^{2}=\frac{\left|a\right|^{2}}{t}$ in $\mathbb{C}$,
then 
\[
det\left(\left[\left(a,b\right)\right]_{t}\right)=\left|a\right|^{2}-t\left|b\right|^{2}=0,
\]
implying that $\left(a,b\right)$ is not invertible in $\mathbb{H}_{t}$.
So, if $t>0$, then $\mathbb{H}_{t}^{\times}$ cannot be a group.
So, if $t\geq0$, then $\mathbb{H}_{t}$ is not a noncommutative field. 
\end{proof}
The proof of the characterization (2.2.3) also proves the following
two corollaries. 
\begin{corollary}
An element $\left(a,b\right)\in\mathbb{H}_{0}^{\times}$ is invertible
in $\mathbb{H}_{0}$, with its inverse $\left(\frac{\overline{a}}{\left|a\right|^{2}},\;\frac{-b}{\left|a\right|^{2}}\right)\in\mathbb{H}_{0}^{\times},$
if and only if the only elements of the subset, 
\[
\left\{ \left(a,b\right)\in\mathbb{H}_{0}^{\times}:a\neq0\right\} \textrm{ of }\mathbb{H}_{0}^{\times}
\]
are invertible in $\mathbb{H}_{0}^{\times}$, if and only if $\left(0,b\right)\in\mathbb{H}_{0}^{\times}$
are not invertible in $\mathbb{H}_{0}^{\times}$, for all $b\in\mathbb{C}$. 
\end{corollary}

\begin{proof}
See the proof of the above theorem. 
\end{proof}
\begin{corollary}
If $t>0$ in $\mathbb{R}$, then an element $\left(a,b\right)\in\mathbb{H}_{t}^{\times}$
is invertible in $\mathbb{H}_{t}^{\times}$ with its inverse, 
\[
\left(\frac{\overline{a}}{\left|a\right|^{2}-t\left|b\right|^{2}},\;\frac{-b}{\left|a\right|^{2}-t\left|b\right|^{2}}\right)\in\mathbb{H}_{t}^{\times},
\]
if and only if $\left(a,b\right)$ is contained in the proper subset,
\[
\left\{ \left(a,b\right):\left|a\right|^{2}\neq t\left|b\right|^{2}\;\mathrm{in\;}\mathbb{R}_{0}^{+}\right\} \;\mathrm{of\;}\mathbb{H}_{t}^{\times}.
\]
\end{corollary}

\begin{proof}
See the proof of the above theorem. 
\end{proof}
By the above results, the negatively-scaled hypercomplex rings $\left\{ \mathbb{H}_{s}\right\} _{s<0}$
are noncommutative fields, meanwhile, the non-negatively scaled hypercomplex
rings $\left\{ \mathbb{H}_{t}\right\} _{t\geq0}$ cannot be noncommutative
fields by (2.2.3). So, for any scale $t\in\mathbb{R}$, the $t$-scaled
hypercomplex ring $\mathbb{H}_{t}$ is decomposed by 
\[
\mathbb{H}_{t}=\mathbb{H}_{t}^{inv}\sqcup\mathbb{H}_{t}^{sing}
\]
with\hfill{}(2.2.4) 
\[
\mathbb{H}_{t}^{inv}=\left\{ \left(a,b\right):\left|a\right|^{2}\neq t\left|b\right|^{2}\right\} ,
\]
and 
\[
\mathbb{H}_{t}^{sing}=\left\{ \left(a,b\right):\left|a\right|^{2}=t\left|b\right|^{2}\right\} ,
\]
where $\sqcup$ is the disjoint union. By (2.2.4), the $t$-scaled
hypercomplex monoid $\mathbb{H}_{t}^{\times}$ is decomposed to be
\[
\mathbb{H}_{t}^{\times}=\mathbb{H}_{t}^{inv}\sqcup\mathbb{H}_{t}^{\times sing},
\]
with\hfill{}(2.2.5) 
\[
\mathbb{H}_{t}^{\times sing}=\mathbb{H}_{t}^{sing}\setminus\left\{ \left(0,0\right)\right\} .
\]

\begin{proposition}
The subset $\mathbb{H}_{t}^{inv}$ forms a non-abelian group in the
monoid $\mathbb{H}_{t}^{\times}$. Meanwhile, the subset $\mathbb{H}_{t}^{\times sing}$
is a semigroup in $\mathbb{H}_{t}^{\times}$ without identity. 
\end{proposition}

\begin{proof}
Let $t\in\mathbb{R}$, and $\mathbb{H}_{t}^{\times}$, the $t$-scaled
hypercomplex monoid, decomposed by (2.2.5). If $h_{1},h_{2}\in\mathbb{H}_{t}^{inv}$,
then $h_{1}\cdot_{t}h_{2}\in\mathbb{H}_{t}^{inv}$, because 
\[
det\left(\left[h_{1}\cdot_{t}h_{2}\right]_{t}\right)=det\left(\left[h_{1}\right]_{t}\left[h_{2}\right]_{t}\right)=det\left(\left[h_{1}\right]_{t}\right)det\left(\left[h_{2}\right]_{t}\right)\neq0.
\]
So, the algebraic pair $\left(\mathbb{H}_{t}^{inv},\cdot_{t}\right)$
forms a group, embedded in the monoid $\mathbb{H}_{t}^{\times}$.
Meanwhile, if $h_{1},h_{2}\in\mathbb{H}_{t}^{\times sing}$, then
$h_{1}\cdot_{t}h_{2}\in\mathbb{H}_{t}^{\times sing}$, since 
\[
det\left(\left[h_{1}\cdot_{t}h_{2}\right]_{t}\right)=det\left(\left[h_{1}\right]_{t}\left[h_{2}\right]_{t}\right)=det\left(\left[h_{1}\right]_{t}\right)det\left(\left[h_{2}\right]_{t}\right)=0.
\]
This restricted operation ($\cdot_{t}$) is associative on $\mathbb{H}_{t}^{\times sing}$,
however, it does not have its identity $\left(1,0\right)$ in $\mathbb{H}_{t}^{\times sing}$,
since $\left(1,0\right)\in\mathbb{H}_{t}^{inv}$. Thus, the pair $\left(\mathbb{H}_{t}^{\times sing},\cdot_{t}\right)$
forms a semigroup without identity in $\mathbb{H}_{t}^{\times}$. 
\end{proof}
The above proposition characterizes the algebraic structures of the
blocks $\mathbb{H}_{t}^{inv}$ and $\mathbb{H}_{t}^{\times sing}$
of the monoid $\mathbb{H}_{t}^{\times}$, as a non-abelian group,
respectively, a semigroup. 
\begin{definition}
The block $\mathbb{H}_{t}^{inv}$ of (2.2.5) is called the group-part
of $\mathbb{H}_{t}^{\times}$ (or, of $\mathbb{H}_{t}$), and the
other algebraic block $\mathbb{H}_{t}^{\times sing}$ of (2.2.5) is
called the semigroup-part of $\mathbb{H}_{t}^{\times}$ (or, of $\mathbb{H}_{t}$). 
\end{definition}

One obtains the following refinements. 
\begin{corollary}
If $t<0$ in $\mathbb{R}$, then $\mathbb{H}_{t}^{\times sing}$ is
empty in $\mathbb{H}_{t}^{\times}$, equivalently, $\mathbb{H}_{t}^{\times}=\mathbb{H}_{t}^{inv}$.
Meanwhile, if $t\geq0$ in $\mathbb{R}$, then $\mathbb{H}_{t}^{\times sing}$
is a non-empty properly embedded semigroup of $\mathbb{H}_{t}^{\times}$,
without identity. 
\end{corollary}

\begin{proof}
The proof is immediately done by the above proposition. 
\end{proof}

\subsection{Spectral Properties of Scaled Hypercomplex Numbers}

By applying the invertibility discussed in Section 2.2, we consider
the spectral analysis on the $t$-scaled hypercomplex ring $\mathbb{H}_{t}$
for $t\in\mathbb{R}$. Let $\left(a,b\right)\in\mathbb{H}_{t}$. Then
its realization $\left[\left(a,b\right)\right]_{t}\in\mathcal{H}_{2}^{t}$
has its characteristic polynomial,

\medskip{}

\hfill{}$det\left(\left[\left(a,b\right)\right]_{t}-z\left[\left(1,0\right)\right]_{t}\right)=z^{2}-2Re\left(a\right)z+det\left(\left[a,b\right]_{t}\right)$,\hfill{}(2.3.1)

\medskip{}

\noindent in a variable $z$ on $\mathbb{C}$, where $Re\left(\bullet\right)$
is the real part in $\mathbb{C}$. Then it has its zeroes,

\medskip{}

\hfill{}$z=Re\left(a\right)\pm\sqrt{Re\left(a\right)^{2}-det\left(\left[\left(a,b\right)\right]_{t}\right)}$\hfill{}(2.3.2)

\noindent (e.g., see {[}1{]} for details).

Recall that the spectrum of an operator $T\in B\left(H\right)$, where
$H$ is a Hilbert space, is defined to be a nonempty compact subset,
\[
spec\left(T\right)=\left\{ \lambda\in\mathbb{C}:T-\lambda I\textrm{ is not invertible on }H\right\} ,
\]
in $\mathbb{C}$, where $I$ is the identity operator on $H$. If
$H$ is $n$-dimensional for $n\in\mathbb{N}$, and hence, if $B\left(H\right)$
is $*$-isomorphic to $M_{n}\left(\mathbb{C}\right)$, then 
\[
spec\left(T\right)=\left\{ w\in\mathbb{C}:det\left(T-wI_{n}\right)=0\right\} ,
\]
where $I_{n}$ is the identity matrix of $M_{n}\left(\mathbb{C}\right)$
(e.g., {[}8{]} and {[}9{]}). 
\begin{proposition}
If $\left(a,b\right)\in\mathbb{H}_{t}$, where 
\[
a=x+yi,\;b=u+vi\in\mathbb{C},
\]
with $x,y,u,v\in\mathbb{R}$ and $i=\sqrt{-1}$ in $\mathbb{C}$,
then

\medskip{}

\hfill{}$spec\left(\left[\left(a,b\right)\right]_{t}\right)=\left\{ x\pm i\sqrt{y^{2}-tu^{2}-tv^{2}}\right\} \;\mathrm{in\;}\mathbb{C}.$\hfill{}(2.3.3) 
\end{proposition}

\begin{proof}
By (2.3.2), we have that 
\[
spec\left(\left[\left(a,b\right)\right]_{t}\right)=\left\{ Re\left(a\right)\pm\sqrt{Re\left(a\right)^{2}-\left(\left|a\right|^{2}-t\left|b\right|^{2}\right)}\right\} ,
\]
in $\mathbb{C}$. Therefore, the set-equality (2.3.3) holds. See {[}1{]}
for details. 
\end{proof}
In the rest, we assume 
\[
a=x+yi\;\mathrm{and\;}b=u+vi\;\mathrm{in\;}\mathbb{C},
\]
with\hfill{}(2.3.4) 
\[
x,y,u,v\in\mathbb{R},\;\mathrm{and\;}i=\sqrt{-1}.
\]

\begin{corollary}
Let $\left(a,b\right)\in\mathbb{H}_{t}$ satisfy (2.3.4). If $Im\left(a\right)^{2}=t\left|b\right|^{2}$
in $\mathbb{R}$, where $Im\left(a\right)$ is the imaginary part
of $a$ in $\mathbb{C}$, then 
\[
spec\left(\left[\left(a,b\right)\right]_{t}\right)=\left\{ x\right\} =\left\{ Re\left(a\right)\right\} \;\mathrm{in\;}\mathbb{R}.
\]
Also, if $Im\left(a\right)^{2}<t\left|b\right|^{2}$ in $\mathbb{R}$,
then 
\[
spec\left(\left[\left(a,b\right)\right]_{t}\right)=\left\{ x\pm\sqrt{tu^{2}+tv^{2}-y^{2}}\right\} \;\mathrm{in\;}\mathbb{R}.
\]
Meanwhile, if $Im\left(a\right)^{2}>t\left|b\right|^{2}$ in $\mathbb{R}$,
then 
\[
spec\left(\left[\left(a,b\right)\right]_{t}\right)=\left\{ x\pm i\sqrt{y^{2}-tu^{2}-tv^{2}}\right\} \;\mathrm{in\;}\mathbb{C}\setminus\mathbb{R}.
\]
\end{corollary}

\begin{proof}
The above refined results are shown case-by-case, by (2.3.3). See
{[}1{]} for details. 
\end{proof}
The above corollary illustrates that: if $Im\left(a\right)^{2}\leq t\left|b\right|^{2}$,
then 
\[
spec\left(\left[\left(a,b\right)\right]_{t}\right)\subset\mathbb{R};
\]
meanwhile, if $Im\left(b\right)^{2}>t\left|b\right|^{2}$, then 
\[
spec\left(\left[\left(a,b\right)\right]_{t}\right)\subset\left(\mathbb{C}\setminus\mathbb{R}\right),\;\mathrm{in\;}\mathbb{C}.
\]
Define the subsets $\mathbb{H}_{t}^{+}$ and $\mathbb{H}_{t}^{-0}$
of $\mathbb{H}_{t}$ by 
\[
\mathbb{H}_{t}^{+}=\left\{ \left(a,b\right)\in\mathbb{H}_{t}:Im\left(a\right)^{2}>t\left|b\right|^{2}\right\} ,
\]
and\hfill{}(2.3.5) 
\[
\mathbb{H}_{t}^{-0}=\left\{ \left(a,b\right)\in\mathbb{H}_{t}:Im\left(a\right)^{2}\leq t\left|b\right|^{2}\right\} .
\]
Then $\mathbb{H}_{t}$ is decomposed by

\medskip{}

\hfill{}$\mathbb{H}_{t}=\mathbb{H}_{t}^{+}\sqcup\mathbb{H}_{t}^{-0}$,\hfill{}(2.3.6)

\medskip{}

\noindent by the above corollary and (2.3.5). 
\begin{corollary}
Let $\mathbb{H}_{t}$ be the $t$-scaled hypercomplex ring for $t\in\mathbb{R}$.
Then it is decomposed to be

\medskip{}

\hfill{}$\begin{array}{cc}
\mathbb{H}_{t}= & \left(\mathbb{H}_{t}^{inv}\cap\mathbb{H}_{t}^{+}\right)\sqcup\left(\mathbb{H}_{t}^{inv}\cap\mathbb{H}_{t}^{-0}\right)\\
\\
 & \sqcup\left(\mathbb{H}_{t}^{sing}\cap\mathbb{H}_{t}^{+}\right)\sqcup\left(\mathbb{H}_{t}^{sing}\cap\mathbb{H}_{t}^{-0}\right),
\end{array}$\hfill{}(2.3.7)

\medskip{}

\noindent where $\mathbb{H}_{t}^{sing}=\mathbb{H}_{t}^{\times sing}\cup\left\{ \left(0,0\right)\right\} $. 
\end{corollary}

\begin{proof}
The partition (2.3.7) of $\mathbb{H}_{t}$ is obtained by (2.2.4)
and (2.3.6). 
\end{proof}
Define a surjection, 
\[
\sigma_{t}:\mathbb{H}_{t}\rightarrow\mathbb{C},
\]
by\hfill{}(2.3.8) 
\[
\sigma_{t}\left(\left(a,b\right)\right)\overset{\textrm{def}}{=}\left\{ \begin{array}{ccc}
a=x+yi &  & \mathrm{if\;}b=0\;\mathrm{in\;}\mathbb{C}\\
\\
x+i\sqrt{y^{2}-tu^{2}-tv^{2}} &  & \mathrm{if\;}b\neq0\;\mathrm{in\;}\mathbb{C},
\end{array}\right.
\]
for all $\left(a,b\right)\in\mathbb{H}_{t}$ satisfying the condition
(2.3.4). 
\begin{definition}
The surjection $\sigma_{t}:\mathbb{H}_{t}\rightarrow\mathbb{C}$ of
(2.3.8) is called the $t$(-scaled)-spectralization on $\mathbb{H}_{t}$.
The images of $\sigma_{t}$ are said to be $t$(-scaled)-spectral
values. 
\end{definition}

\begin{definition}
Let $\xi\in\mathbb{H}_{t}$ have its $t$-spectral value $\sigma_{t}\left(\xi\right)\in\mathbb{C}$.
The realization of $\left(\sigma_{t}\left(\xi\right),0\right)$,

\medskip{}

\hfill{}$\left[\left(\sigma_{t}\left(\xi\right),0\right)\right]_{t}=\left(\begin{array}{ccc}
\sigma_{t}\left(\xi\right) &  & 0\\
\\
0 &  & \overline{\sigma_{t}\left(\xi\right)}
\end{array}\right)\in\mathcal{H}_{2}^{t},$\hfill{}(2.3.9)

\medskip{}

\noindent is called the $t$(-scaled)-spectral form of $\xi$, denoted
by $\Sigma_{t}\left(\xi\right)$. 
\end{definition}

Note that the conjugate-notation in (2.3.9) is symbolically understood
in the sense that: if 
\[
\sigma_{t}\left(\left(a,b\right)\right)=x+i\sqrt{y^{2}-tu^{2}-tv^{2}},
\]
satisfying 
\[
y^{2}-tu^{2}-tv^{2}=Im\left(a\right)^{2}-t\left|b\right|^{2}<0,
\]
under (2.3.4), equivalently, if 
\[
\sigma_{t}\left(\left(a,b\right)\right)=x-\sqrt{tu^{2}-tv^{2}-y^{2}}\in\mathbb{R},
\]
then the symbol, 
\[
\overline{\sigma_{t}\left(\left(a,b\right)\right)}\overset{\textrm{means}}{=}\overline{x+i\sqrt{R}}=x-i\sqrt{R}=x+\sqrt{tu^{2}-tv^{2}-y^{2}},
\]
in $\mathbb{R}$, where $R=y^{2}-tu^{2}-tv^{2}$ in $\mathbb{R}$.
While, if $R\geq0$, and hence, if 
\[
\sigma_{t}\left(\left(a,b\right)\right)=x+i\sqrt{R}\;\;\mathrm{in\;\;}\mathbb{C},
\]
then $\overline{\sigma_{t}\left(\left(a,b\right)\right)}=x-i\sqrt{R}$
is the usual conjugate of $\sigma_{t}\left(\left(a,b\right)\right)$
in $\mathbb{C}$. 
\begin{definition}
Two hypercomplex numbers $\xi,\eta\in\mathbb{H}_{t}$ are said to
be $t$(-scaled)-spectral-related, if $\sigma_{t}\left(\xi\right)=\sigma_{t}\left(\eta\right),\;\mathrm{in\;}\mathbb{C}.$ 
\end{definition}

The $t$-spectral relation is an equivalence relation on $\mathbb{H}_{t}$
(e.g., {[}1{]}). So, the equivalence class $\widetilde{\xi}$ of a
hypercomplex number $\xi\in\mathbb{H}_{t}$ is well-determined to
be 
\[
\widetilde{\xi}\overset{\textrm{def}}{=}\left\{ \eta\in\mathbb{H}_{t}:\eta\textrm{ is }t\textrm{-spectral related to }\xi\right\} \;\mathrm{in\;}\mathbb{H}_{t},
\]
as a subset of $\mathbb{H}_{t}$. So, the quotient set,

\medskip{}

\hfill{}$\widetilde{\mathbb{H}_{t}}\overset{\textrm{def}}{=}\left\{ \widetilde{\xi}:\xi\in\mathbb{H}_{t}\right\} ,$\hfill{}(2.3.10)

\medskip{}

\noindent is constructed. This quotient set $\widetilde{\mathbb{H}_{t}}$
of (2.3.10) is equipotent (or, bijective) to $\mathbb{C}$, by the
surjectivity of $\sigma_{t}$. So, each complex number $z\in\mathbb{C}$
represents all $t$-scaled hypercomplex numbers $\eta\in\mathbb{H}_{t}$,
having their $t$-spectral values $\sigma_{t}\left(\eta\right)=z$.

\subsection{The $t$-Scaled-Spectral Relation on $\mathbb{H}_{t}$}

Recall that two operators $T$ and $S$ are said to be similar in
the operator algebra $B\left(H\right)$ on a Hilbert space $H$, if
there is an invertible operator $U\in B\left(H\right)$, such that
\[
S=U^{-1}TU\;\;\;\mathrm{in\;\;\;}B(H).
\]

\noindent In {[}3{]}, it is shown that: the $\left(-1\right)$-spectral
form $\Sigma_{-1}\left(\eta\right)$ and the realization $\left[\eta\right]_{-1}$
are similar ``in $\mathcal{H}_{2}^{-1}$,'' in the sense that: there
exists $[h]_{-1}\in\mathcal{H}_{2}^{t}$, such that 
\[
\Sigma_{-1}\left(\eta\right)=\left[h\right]_{-1}^{-1}\left[\eta\right]_{-1}\left[h\right]_{-1}=\left[h^{-1}\cdot_{-1}\eta\cdot_{-1}h\right]_{-1},
\]
in $\mathcal{H}_{2}^{-1}$, for ``all'' $\eta\in\mathbb{H}_{-1}$,
where $\mathcal{H}_{2}^{-1}$ is the realization of the quaternions
$\mathbb{H}_{-1}=\mathbb{H}$. i.e., the realizations of all quaternions
of $\mathbb{H}_{-1}$ are similar to their $\left(-1\right)$-spectral
forms in $\mathcal{H}_{2}^{-1}$. In {[}1{]}, we studied such an invertibility
on $\left\{ \mathbb{H}_{t}\right\} _{t\in\mathbb{R}}$, in particular,
where $t\neq-1$ in $\mathbb{R}$. 
\begin{definition}
Two realizations $T,S\in\mathcal{H}_{2}^{t}$ of hypercomplex numbers
of $\mathbb{H}_{t}$ are similar ``in $\mathcal{H}_{2}^{t}$,''
if there exists an invertible $U\in\mathcal{H}_{2}^{t}$, such that
\[
S=U^{-1}TU\;\;\mathrm{in\;\;}\mathcal{H}_{2}^{t}.
\]
Two hypercomplex numbers $\xi\;\mathrm{and\;}\eta$ are said to be
similar in $\mathbb{H}_{t}$, if their realizations $\left[\xi\right]_{t}$
and $\left[\eta\right]_{t}$ are similar in $\mathcal{H}_{2}^{t}$. 
\end{definition}

If $\left(a,b\right)\in\mathbb{H}_{t}$ satisfies the condition (2.3.4),
then 
\[
\left[\left(a,b\right)\right]_{t}=\left(\begin{array}{cc}
a & tb\\
\overline{b} & \overline{a}
\end{array}\right)\in\mathcal{H}_{2}^{t},
\]
having its determinant, 
\[
det\left(\left[\left(a,b\right)\right]_{t}\right)=\left|a\right|^{2}-t\left|b\right|^{2}=\left(x^{2}+y^{2}\right)-t\left(u^{2}+v^{2}\right),
\]
meanwhile, 
\[
\Sigma_{t}\left(\left(a,b\right)\right)=\left(\begin{array}{ccc}
x+i\sqrt{y^{2}-tu^{2}-tv^{2}} &  & 0\\
\\
0 &  & x-i\sqrt{y^{2}-tu^{2}-tv^{2}}
\end{array}\right),
\]
in $\mathcal{H}_{2}^{t}$, having its determinant, 
\[
det\left(\Sigma_{t}\left(\left(a,b\right)\right)\right)=x^{2}+\left|y^{2}-tu^{2}-tv^{2}\right|.
\]
The above computations show that 
\[
det\left(\left[\left(a,b\right)\right]_{t}\right)\neq det\left(\Sigma_{t}\left(\left(a,b\right)\right)\right),\;\mathrm{in\;general},
\]
implying that $\left[\left(a,b\right)\right]_{t}$ and $\Sigma_{t}\left(\left(a,b\right)\right)$
are not similar in $M_{2}\left(\mathbb{C}\right)$, and hence, in
$\mathcal{H}_{2}^{t}$, in general, for some $t\in\mathbb{R}$, especially,
where $t\geq0$. 
\begin{lemma}
If $t<0$ in $\mathbb{R}$, then every hypercomplex number $h\in\mathbb{H}_{t}$
is similar to $\left(\sigma_{t}\left(h\right),0\right)\in\mathbb{H}_{t}$,
where $\sigma_{t}\left(h\right)$ is the $t$-spectral value of $h$,
i.e.,

\medskip{}

\hfill{}$t<0$$\Longrightarrow$ $\left[h\right]_{t}$ and $\Sigma_{t}\left(h\right)$
are similar in $\mathcal{H}_{2}^{t}$.\hfill{}(2.4.1) 
\end{lemma}

\begin{proof}
Suppose $t<0$. If $h=\left(a,0\right)\in\mathbb{H}_{t}$, then 
\[
\left[\left(a,0\right)\right]_{t}=\left(\begin{array}{cc}
a & 0\\
0 & \overline{a}
\end{array}\right)=\Sigma_{t}\left(\left(a,0\right)\right)\;\mathrm{in\;}\mathcal{H}_{2}^{t},
\]
since $\sigma_{t}\left(\left(a,0\right)\right)=a$ in $\mathbb{C}$.
So, $\left[\left(a,0\right)\right]_{t}$ and $\Sigma_{t}\left(\left(a,0\right)\right)$
are automatically similar in $\mathcal{H}_{2}^{t}$. Assume now that
$h=\left(a,b\right)\in\mathbb{H}_{t}$ with $b\neq0$, where $t<0$.
Then $\left[h\right]_{t}$ and $\Sigma_{t}\left(h\right)$ are similar
in $\mathcal{H}_{2}^{t}$, because there exists 
\[
q_{h}=\left(1,\;\;\overline{\frac{w-a}{tb}}\right)\in\mathbb{H}_{t},
\]
such that 
\[
\Sigma_{t}\left(h\right)=\left[q_{h}\right]_{t}^{-1}\left[h\right]_{t}\left[q_{h}\right]_{t}=\left[q_{h}^{-1}\cdot_{t}h\cdot_{t}q_{h}\right]_{t},
\]
in $\mathcal{H}_{2}^{t}$ (See {[}1{]} for details). 
\end{proof}
The relation (2.4.1) shows that if a scale $t$ is negative in $\mathbb{R}$,
then ``all'' hypercomplex numbers $h\in\mathbb{H}_{t}$ are similar
to $\left(\sigma_{t}\left(h\right),0\right)\in\mathbb{H}_{t}$ induced
by their $t$-spectral values $\sigma_{t}\left(h\right)$; and it
generalizes the quaternionic case where $t=-1$, proven independently
in {[}3{]}. 
\begin{proposition}
If $t<0$ in $\mathbb{R}$, then the $t$-spectral relation and the
similarity are equivalent on $\mathbb{H}_{t}$, i.e.,

\medskip{}

\hfill{}$t<0\Longrightarrow[t\textrm{-spectral relation}\overset{\textrm{equi}}{=}\textrm{similarity on }\mathbb{H}_{t},]$\hfill{}(2.4.2)

\medskip{}

\noindent where ``$\overset{\textrm{equi}}{=}$'' means ``being
equivalent to, as equivalence relations.'' 
\end{proposition}

\begin{proof}
The relation (2.4.2) is shown by (2.4.1). See {[}1{]} for details. 
\end{proof}
How about the cases where given scales $t$ are nonnegative in $\mathbb{R}$?
If $t\geq0$, then there do exist $h\in\mathbb{H}_{t}$, such that
$\left[h\right]_{t}$ and $\Sigma_{t}\left(h\right)$ are not similar
in $\mathcal{H}_{2}^{t}$ as we discussed at the beginning of this
section. However, one may / can verify that if $t\geq0$, and 
\[
h\in\mathbb{H}_{t}^{inv}\cap\mathbb{H}_{t}^{+},\;\textrm{an block of the partition (2.3.7),}
\]
in $\mathbb{H}_{t}$, then $\left[h\right]_{t}$ and $\Sigma_{t}\left(h\right)$
seem to be similar in $\mathcal{H}_{2}^{t}$.

\section{Scaled-Hypercomplex Conjugation}

In this section, we consider conjugation-processes on our scaled hypercomplex
rings $\left\{ \mathbb{H}_{t}\right\} _{t\in\mathbb{R}}$, motivated
by the adjoints introduced in {[}1{]}. Fix an arbitrary scale $t\in\mathbb{R}$
throughout this section, and let $\mathbb{H}_{t}$ be the corresponding
$t$-scaled hypercomplex ring. Define a function, 
\[
J:\mathbb{H}_{t}\rightarrow\mathbb{H}_{t},
\]
by\hfill{}(3.1) 
\[
J\left(\left(a,b\right)\right)\overset{\textrm{def}}{=}\left(\overline{a},\;-b\right),\;\;\forall\left(a,b\right)\in\mathbb{H}_{t}.
\]
By the very definition (3.1), indeed, this morphism $J$ is a well-defined
function. Moreover, if 
\[
h_{1}=\left(a_{1},b_{1}\right)\neq\left(a_{2},b_{2}\right)=h_{2}\;\;\mathrm{in\;\;}\mathbb{H}_{t},
\]
then 
\[
J\left(h_{1}\right)=\left(\overline{a_{1}},-b_{1}\right)\neq\left(\overline{a_{2}},-b_{2}\right)=J\left(h_{2}\right),
\]
in $\mathbb{H}_{2}$, and hence, this function $J$ is injective;
also, for any $\left(a,b\right)\in\mathbb{H}_{t}$, there exists $\left(\overline{a},-b\right)\in\mathbb{H}_{t}$,
such that 
\[
J\left(\left(\overline{a},-b\right)\right)=\left(\overline{\overline{a}},\:-\left(-b\right)\right)=\left(a,b\right),
\]
in $\mathbb{H}_{t}$, showing the surjectivity of $J$. Thus, the
function $J$ of (3.1) is a bijection.

Remember that $\mathbb{H}_{t}$ is topological-ring-isomorphic to
its realization $\mathcal{H}_{2}^{t}$ by (2.1.5). So, one may understand
the bijection $J$ of (3.1) as a function, say $J_{2}$, acting on
$\mathcal{H}_{2}^{t}$, i.e., 
\[
J_{2}:\mathcal{H}_{2}^{t}\rightarrow\mathcal{H}_{2}^{t}
\]
defined by\hfill{}(3.2) 
\[
J_{2}\left(\left[\left(a,b\right)\right]_{t}\right)\overset{\textrm{def}}{=}\left[J\left(\left(a,b\right)\right)\right]_{t}=\left[\left(\overline{a},-b\right)\right]_{t},
\]
for all $\left(a,b\right)\in\mathbb{H}_{t}$. Since the action $\pi_{t}:\mathbb{H}_{t}\rightarrow\mathcal{H}_{2}^{t}$
and the function $J$ are bijective, the above function $J_{2}$ of
(3.2) is not only well-defined, but also bijective on $\mathcal{H}_{2}^{t}$. 
\begin{theorem}
The bijection $J$ of (3.1) acting on $\mathbb{H}_{t}$ is an adjoint
on $\mathbb{H}_{t}$ in the sense that: for all $h_{1},h_{2}\in\mathbb{H}_{t}$,
\[
J^{2}\left(h_{1}\right)=J\left(J\left(h_{1}\right)\right)=h_{1}
\]
\[
J\left(h_{1}+h_{2}\right)=J\left(h_{1}\right)+J\left(h_{2}\right),
\]
\[
J\left(h_{1}\cdot_{t}h_{2}\right)=J\left(h_{2}\right)\cdot_{t}J\left(h_{1}\right),
\]
in addition to\hfill{}(3.3) 
\[
J\left(\left(z,0\right)\cdot_{t}h\right)=\left(\overline{z},0\right)\cdot_{t}J\left(h\right),
\]
for all $z\in\mathbb{C}$ and $h\in\mathbb{H}_{t}$. 
\end{theorem}

\begin{proof}
Since the $t$-scaled hypercomplex ring $\mathbb{H}_{t}$ is topological-ring-isomorphic
to its realization $\mathcal{H}_{2}^{t}$, by the definition (3.2)
of the bijection $J_{2}$ on $\mathcal{H}_{2}^{t}$, it is sufficient
to show that $J_{2}$ is an adjoint on $\mathcal{H}_{2}^{t}$ satisfying
(3.3) under representation. Observe that, for all $\left(a,b\right)\in\mathbb{H}_{t}$,
we have

\medskip{}

$\;\;\;\;$$J_{2}^{2}\left(\left[\left(a,b\right)\right]_{t}\right)=J_{2}\left(J_{2}\left(\left[\left(a,b\right)\right]_{t}\right)\right)=J_{2}\left(\left[J\left(\left(a,b\right)\right)\right]_{t}\right)$

\medskip{}

$\;\;\;\;\;\;\;\;\;\;\;\;$$=J_{2}\left(\left[\left(\overline{a},-b\right)\right]_{t}\right)=\left[J\left(\overline{a},-b\right)\right]_{t}=\left[\left(\overline{\overline{a}},\:-\left(-b\right)\right)\right]_{t}$

\medskip{}

$\;\;\;\;\;\;\;\;\;\;\;\;\;\;$$=\left[\left(a,b\right)\right]_{t}$;

\medskip{}

\noindent also, for any $\left(a_{l},b_{l}\right)\in\mathbb{H}_{t}$,
for $l=1,2$,

\medskip{}

$\;\;\;\;$$J_{2}\left(\left[\left(a_{1},b_{1}\right)\right]_{t}+\left[\left(a_{2},b_{2}\right)\right]_{t}\right)=J_{2}\left(\left[\left(a_{1}+a_{2},b_{1}+b_{2}\right)\right]_{t}\right)$

\medskip{}

$\;\;\;\;\;\;\;\;$$=\left[J\left(\left(a_{1}+a_{2},b_{1}+b_{2}\right)\right)\right]_{t}=\left[\left(\overline{a_{1}+a_{2}},\;-\left(b_{1}+b_{2}\right)\right)\right]_{t}$

\medskip{}

$\;\;\;\;\;\;\;\;$$=\left[\left(\overline{a_{1}},-b_{1}\right)\right]_{t}+\left[\left(\overline{a_{2}},-b_{2}\right)\right]_{t}=\left[J\left(\left(a_{1},b_{1}\right)\right)\right]_{t}+\left[J\left(a_{2},b_{2}\right)\right]_{t}$

\medskip{}

$\;\;\;\;\;\;\;\;$$=J_{2}\left(\left[\left(a_{1},b_{1}\right)\right]_{t}\right)+J_{2}\left(\left[\left(a_{2},b_{2}\right)\right]_{t}\right)$;

\medskip{}

\noindent and

\medskip{}

$\;\;\;$$J_{2}\left(\left[\left(a_{1},b_{1}\right)\right]_{t}\left[\left(a_{2},b_{2}\right)\right]_{t}\right)=J_{2}\left(\left(\begin{array}{ccc}
a_{1}a_{2}+tb_{1}\overline{b_{2}} &  & t\left(a_{1}b_{2}+b_{1}\overline{a_{2}}\right)\\
\\
\overline{a_{1}b_{2}+b_{1}\overline{a_{2}}} &  & \overline{a_{1}a_{2}+tb_{1}\overline{b_{2}}}
\end{array}\right)\right)$

\medskip{}

$\;\;\;\;\;\;\;\;\;\;\;\;\;\;$$=\left(\begin{array}{ccc}
\overline{a_{1}a_{2}+tb_{1}\overline{b_{2}}} &  & t\left(-a_{1}b_{2}-b_{1}\overline{a_{2}}\right)\\
\\
\overline{-a_{1}b_{2}-b_{1}\overline{a_{2}}} &  & a_{1}a_{2}+tb_{1}\overline{b_{2}}
\end{array}\right)$

\medskip{}

$\;\;\;\;\;\;\;\;\;\;\;\;\;\;$$=\left(\begin{array}{ccc}
\overline{a_{2}} &  & t\left(-b_{2}\right)\\
\\
\overline{-b_{2}} &  & a_{2}
\end{array}\right)\left(\begin{array}{ccc}
\overline{a_{1}} &  & t\left(-b_{1}\right)\\
\\
\overline{-b_{1}} &  & a_{1}
\end{array}\right)$

\medskip{}

$\;\;\;\;\;\;\;\;\;\;\;\;\;\;$$=\left(J_{2}\left(\left[\left(a_{2},b_{2}\right)\right]_{t}\right)\right)\left(J_{2}\left(\left[\left(a_{1},b_{1}\right)\right]_{t}\right)\right)$.

\medskip{}

\noindent Moreover, if $z\in\mathbb{C}$ inducing $\left(z,0\right)\in\mathbb{H}_{t}$
and $\left(a,b\right)\in\mathbb{H}_{t}$, then

\medskip{}

$\;\;\;\;$$J_{2}\left(\left[\left(z,0\right)\right]_{t}\left[\left(a,b\right)\right]_{t}\right)=J_{2}\left(\left(\begin{array}{ccc}
za &  & tzb\\
\\
\overline{zb} &  & \overline{za}
\end{array}\right)\right)$

\medskip{}

$\;\;\;\;\;\;\;\;\;\;\;\;$$=\left(\begin{array}{ccc}
\overline{za} &  & t\left(-zb\right)\\
\\
\overline{-zb} &  & za
\end{array}\right)=\left[\left(\overline{z},0\right)\right]_{t}\left[\left(\overline{a},-b\right)\right]_{t}$

\medskip{}

$\;\;\;\;\;\;\;\;\;\;\;\;$$=\left(\left[\left(\overline{z},0\right)\right]_{t}\right)\left(J_{2}\left(\left[\left(a,b\right)\right]_{t}\right)\right)$.

\medskip{}

\noindent Therefore, the bijection $J_{2}$ of (3.2) is an adjoint
on $\mathcal{H}_{2}^{t}$. Therefore, the bijection $J$ of (3.1)
is an adjoint on $\mathbb{H}_{t}$, satisfying the conditions of (3.3),
under the inverse-action of $\pi_{t}^{-1}:\mathcal{H}_{2}^{t}\rightarrow\mathbb{H}_{t}$. 
\end{proof}
The above theorem shows that the bijection $J$ is a well-defined
adjoint on $\mathbb{H}_{t}$ by (3.3). To distinguish the difference
between the usual adjoint (which is the conjugate-transpose) on $M_{2}\left(\mathbb{C}\right)$,
and the adjoint $J$ on $\mathbb{H}_{t}$ (or, the adjoint $J_{2}$
on $\mathcal{H}_{2}^{t}$), we call $J$ (and $J_{2}$) the hypercomplex-conjugate.
Here, note that such a bijection $J$ (and $J_{2}$) is free from
the choice of scales $t\in\mathbb{R}$. So, instead of calling $J$
(and $J_{2}$) the $t$-scaled hypercomplex-conjugate, we simply call
it the hypercomplex-conjugate. 
\begin{definition}
The bijection $J$ of (3.1), 
\[
J:\left(a,b\right)\in\mathbb{H}_{t}\longmapsto\left(\overline{a},-b\right)\in\mathbb{H}_{t},
\]
is called the hypercomplex-conjugate on $\mathbb{H}_{t}$, for ``all''
scales $t\in\mathbb{R}$. The bijection $J_{2}$ of (3.2), 
\[
J_{2}:\left[\left(a,b\right)\right]_{t}\in\mathcal{H}_{2}^{t}\longmapsto\left[\left(\overline{a},-b\right)\right]_{t}\in\mathcal{H}_{2}^{t},
\]
is also called the hypercomplex-conjugate on $\mathcal{H}_{2}^{t}$,
for all $t\in\mathbb{R}$. For convenience, we denote the images $J\left(\left(a,b\right)\right)$
by $\left(a,b\right)^{\dagger}$ in $\mathbb{H}_{t}$, and similarly,
the images $J_{2}\left(\left[\left(a,b\right)\right]_{t}\right)$
are denoted by $\left[\left(a,b\right)\right]_{t}^{\dagger}$ in $\mathcal{H}_{2}^{t}$,
i.e., 
\[
\left(a,b\right)^{\dagger}=\left(\overline{a},-b\right)\;\;\mathrm{in\;\;}\mathbb{H}_{t},
\]
and 
\[
\left[\left(a,b\right)\right]_{t}^{\dagger}=\left[\left(\overline{a},-b\right)\right]_{t}=\left[\left(a,b\right)^{\dagger}\right]_{t}\;\mathrm{in\;}\mathcal{H}_{2}^{t},
\]
for all $\left(a,b\right)\in\mathbb{H}_{t}$, for all $t\in\mathbb{R}$. 
\end{definition}

If $\left(a,b\right)\in\mathbb{H}_{t}$, then one can get that 
\[
\left[\left(a,b\right)\right]_{t}^{\dagger}\left[\left(a,b\right)\right]_{t}=\left(\begin{array}{ccc}
\left|a\right|^{2}-t\left|b\right|^{2} &  & 0\\
\\
0 &  & \left|a\right|^{2}-t\left|b\right|^{2}
\end{array}\right)=\left[\left(a,b\right)\right]_{t}\left[\left(a,b\right)\right]_{t}^{\dagger},
\]
implying that

\medskip{}

\hfill{}$\left[h\right]_{t}^{\dagger}\left[h\right]_{t}=\left[\left(\left|a\right|^{2}-t\left|b\right|^{2},\;0\right)\right]_{t}=\left[h\right]_{t}\left[h\right]_{t}^{\dagger},$\hfill{}(3.4)

\medskip{}

\noindent for all $h=\left(a,b\right)\in\mathbb{H}_{t}$, for ``all''
$t\in\mathbb{R}$. 
\begin{theorem}
If $\left(a,b\right)\in\mathbb{H}^{t}$, then

\medskip{}

\hfill{}$\left(a,b\right)^{\dagger}\cdot_{t}\left(a,b\right)=\left(\left|a\right|^{2}-t\left|b\right|^{2},\;0\right)=\left(a,b\right)\cdot_{t}\left(a,b\right)^{\dagger},$\hfill{}(3.5)

\medskip{}

\noindent in $\mathbb{H}^{t}$, for all $t\in\mathbb{R}$. It implies
that 
\[
\sigma_{t}\left(\left(a,b\right)^{\dagger}\cdot_{t}\left(a,b\right)\right)=\left|a\right|^{2}-t\left|b\right|^{2}=\sigma_{t}\left(\left(a,b\right)\cdot_{t}\left(a,b\right)^{\dagger}\right),
\]

\noindent i.e.,\hfill{}(3.6) 
\[
\sigma_{t}\left(\left(a,b\right)^{\dagger}\cdot_{t}\left(a,b\right)\right)=det\left(\left[\left(a,b\right)\right]_{t}\right)=\sigma_{t}\left(\left(a,b\right)\cdot_{t}\left(a,b\right)^{\dagger}\right),
\]
for all $\left(a,b\right)\in\mathbb{H}_{t}$, for all $t\in\mathbb{R}$. 
\end{theorem}

\begin{proof}
The relation (3.5) is proven by (3.4). By (3.5), the first $t$-spectral-value
relation of (3.6) is immediately obtained. Also, since 
\[
det\left(\left[\left(a,b\right)\right]_{t}\right)=\left|a\right|^{2}-t\left|b\right|^{2},
\]
the second formula of (3.6) holds true, too. 
\end{proof}
The formulas of (3.6) show the connections between the determinant
on $\mathcal{H}_{2}^{t}$ and our hypercomplex-conjugation (3.1),
or (3.2).

\section{Scaled Hypercomplex Rings $\mathbb{H}_{t}$ as Topological Vector
Spaces}

Throughout this section, we fix a scale $t\in\mathbb{R}$, and the
corresponding $t$-scaled hypercomplex ring $\mathbb{H}_{t}$. In
Section 3, we showed that the ring $\mathbb{H}_{t}$ has its well-defined
adjoint, the hypercomplex-conjugate ($\dagger$) on it, defined by
\[
\left(a,b\right)^{\dagger}=\left(\overline{a},\:-b\right),\;\;\forall\left(a,b\right)\in\mathbb{H}_{t},
\]
inducing an adjoint, also denoted by ($\dagger$), on the $t$-scaled
realization $\mathcal{H}_{2}^{t}$, 
\[
\left[\left(a,b\right)\right]_{t}^{\dagger}=\left[\left(a,b\right)^{\dagger}\right]_{t}=\left[\left(\overline{a},-b\right)\right]_{t},
\]
satisfying that\hfill{}(4.1) 
\[
\sigma_{t}\left(\left(a,b\right)^{\dagger}\cdot_{t}\left(a,b\right)\right)=det\left(\left[\left(a,b\right)\right]_{t}\right)=\sigma_{t}\left(\left(a,b\right)\cdot_{t}\left(a,b\right)^{\dagger}\right).
\]
by (3.5) and (3.6), for all $\left(a,b\right)\in\mathbb{H}_{t}$ (in
fact, for ``all'' scales).

Recall that the $t$-scaled realization $\mathcal{H}_{2}^{t}$ is
a sub-structure of $M_{2}\left(\mathbb{C}\right)$. So, the normalized
trace, 
\[
\tau=\frac{1}{2}tr\;\;\;\;\;\mathrm{on\;\;\;\;\;}M_{2}\left(\mathbb{C}\right),
\]
can be restricted to $\tau\mid_{\mathcal{H}_{2}^{t}}$, also denoted
simply by $\tau$, on $\mathcal{H}_{2}^{t}$, where $tr$ is the usual
trace on $M_{2}\left(\mathbb{C}\right)$, 
\[
tr\left(\left(\begin{array}{cc}
a_{11} & a_{12}\\
a_{21} & a_{22}
\end{array}\right)\right)=a_{11}+a_{22},\;\forall\left[a_{ij}\right]_{2\times2}\in M_{2}\left(\mathbb{C}\right).
\]
Observe that, for any $\left[\left(a,b\right)\right]_{t}\in\mathcal{H}_{2}^{t}$,
one has 
\[
\tau\left(\left[\left(a,b\right)\right]_{t}\right)=\frac{1}{2}tr\left(\left(\begin{array}{cc}
a & tb\\
\overline{b} & \overline{a}
\end{array}\right)\right)=\frac{1}{2}\left(a+\overline{a}\right),
\]
i.e.,\hfill{}(4.2) 
\[
\tau\left(\left[\left(a,b\right)\right]_{t}\right)=Re\left(a\right),\;\;\forall\left(a,b\right)\in\mathbb{H}_{t}.
\]
Therefore, without loss of generality, one can define a linear functional,
also denoted by $\tau$, on $\mathbb{H}_{t}$, by

\medskip{}

\hfill{}$\tau\left(\left(a,b\right)\right)\overset{\textrm{def}}{=}Re\left(a\right),\;\;\;\forall\left(a,b\right)\in\mathbb{H}_{t},$\hfill{}(4.3)

\medskip{}

\noindent by (4.2). Then this well-defined linear functional $\tau$
of (4.3) has the following relation with our $t$-spectralization
$\sigma_{t}$. 
\begin{proposition}
Let $\left(a,b\right)\in\mathbb{H}_{t}$ satisfy the condition (2.3.4),
and let $\tau$ be the linear functional (4.3) on $\mathbb{H}_{t}$.
Then

\medskip{}

\hfill{}$\tau\left(\left(a,b\right)\right)=x=Re\left(\sigma_{t}\left(\left(a,b\right)\right)\right).$\hfill{}(4.4)

\medskip{}

\noindent In particular, if the $t$-spectral value $w\overset{\textrm{denote}}{=}\sigma_{t}\left(\left(a,b\right)\right)$
is polar-decomposed to be $w=re^{i\theta}$ in $\mathbb{C}$, with
$r=\left|w\right|$ and $\theta=Arg\left(w\right)$, then

\medskip{}

\hfill{}$\tau\left(\left(a,b\right)\right)=rcos\theta.$\hfill{}(4.5) 
\end{proposition}

\begin{proof}
Suppose a hypercomplex number $\left(a,b\right)\in\mathbb{H}_{t}$
satisfies the condition (2.3.4), i.e., 
\[
a=x+yi,\;\mathrm{and\;}b=u+vi,\;\mathrm{in\;\mathbb{C},}
\]
with $x,y,u,v\in\mathbb{R}$ and $i=\sqrt{-1}$. Then one obtains
that 
\[
\sigma_{t}\left(\left(a,b\right)\right)=x+i\sqrt{y^{2}-tu^{2}-tv^{2}},
\]
by (2.3.3) and (2.3.9). So, by the definition (4.3) of $\tau$, 
\[
\tau\left(\left(a,b\right)\right)=Re\left(a\right)=x=Re\left(\sigma_{t}\left(\left(a,b\right)\right)\right).
\]
Therefore, the relation (4.4) holds. So, the equivalent formula (4.5)
of (4.4) holds under polar decomposition of $\sigma_{t}\left(\left(a,b\right)\right)$. 
\end{proof}
Define now a form, 
\[
\left\langle ,\right\rangle _{t}:\mathbb{H}_{t}\times\mathbb{H}_{t}\rightarrow\mathbb{R}\subset\mathbb{C},
\]
by\hfill{}(4.6) 
\[
\left\langle h_{1},h_{2}\right\rangle _{t}\overset{\textrm{def}}{=}\tau\left(h_{1}\cdot_{t}h_{2}^{\dagger}\right),\;\;\;\forall h_{1},h_{2}\in\mathbb{H}_{t},
\]
where the linear functional $\tau$ in (4.6) is in the sense of (4.3),
whose range is contained in $\mathbb{R}$ in $\mathbb{C}$. By the
bijectivity of our hypercomplex-conjugation ($\dagger$), this form
(4.6) is well-defined as a function. Moreover, it satisfies that

\medskip{}

$\left\langle \left(a_{1},b_{1}\right)+\left(a_{2},b_{2}\right),\:\left(a_{3},b_{3}\right)\right\rangle _{t}$

\medskip{}

$=\left\langle \left(a_{1}+a_{2},\:b_{1}+b_{2}\right),\:\left(a_{3},b_{3}\right)\right\rangle _{t}$

\medskip{}

$=\tau\left(\left(a_{1}+a_{2},b_{1}+b_{2}\right)\cdot_{t}\left(a_{3},b_{3}\right)^{\dagger}\right)$

\medskip{}

$=\tau\left(\left(\begin{array}{ccc}
a_{1}+a_{2} &  & t\left(b_{1}+b_{2}\right)\\
\\
\overline{b_{1}+b_{2}} &  & \overline{a_{1}+a_{2}}
\end{array}\right)\left(\begin{array}{cc}
\overline{a_{3}} & t\left(-b_{3}\right)\\
\overline{-b_{3}} & a_{3}
\end{array}\right)\right)$

\medskip{}

\noindent where $\tau$ is in the sense of (4.2)

\medskip{}

$=\tau\left(\left(\begin{array}{ccc}
a_{1}\overline{a_{3}}+a_{2}\overline{a_{3}}-t\left(b_{1}\overline{b_{3}}+b_{2}\overline{b_{3}}\right) &  & t\left(-a_{1}b_{3}-a_{2}b_{3}+a_{3}b_{1}+a_{3}b_{2}\right)\\
\\
\overline{-a_{1}b_{3}-a_{2}b_{3}+a_{3}b_{1}+a_{3}b_{2}} &  & \overline{a_{1}\overline{a_{3}}+a_{2}\overline{a_{3}}-t\left(b_{1}\overline{b_{3}}+b_{2}\overline{b_{3}}\right)}
\end{array}\right)\right)$

\medskip{}

$=Re\left(a_{1}\overline{a_{3}}+a_{2}\overline{a_{3}}-t\left(b_{1}\overline{b_{3}}+b_{2}\overline{b_{3}}\right)\right)$

\medskip{}

$=Re\left(a_{1}\overline{a_{3}}-tb_{1}\overline{b_{3}}\right)+Re\left(a_{2}\overline{a_{3}}-tb_{2}\overline{b_{3}}\right)$

\medskip{}

$=\tau\left(\left(a_{1},b_{1}\right)\cdot_{t}\left(a_{3},b_{3}\right)^{\dagger}\right)+\tau\left(\left(a_{2},b_{2}\right)\cdot_{t}\left(a_{3},b_{3}\right)^{\dagger}\right)$

\medskip{}

$=\left\langle \left(a_{1},b_{1}\right),\left(a_{3},b_{3}\right)\right\rangle _{t}+\left\langle \left(a_{2},b_{2}\right),\left(a_{3},b_{3}\right)\right\rangle _{t}$,

\medskip{}

\noindent for all $\left(a_{l},b_{l}\right)\in\mathbb{H}_{t}$, for
$l=1,2,3$, i.e.,

\medskip{}

\hfill{}$\left\langle h_{1}+h_{2},h_{3}\right\rangle _{t}=\left\langle h_{1},h_{3}\right\rangle _{t}+\left\langle h_{2},h_{3}\right\rangle _{t},$\hfill{}(4.7)

\medskip{}

\noindent for all $h_{1},h_{2},h_{3}\in\mathbb{H}_{t}$. Similarly,
one obtains that

\medskip{}

\hfill{}$\left\langle h_{1},h_{2}+h_{3}\right\rangle _{t}=\left\langle h_{1},h_{2}\right\rangle _{t}+\left\langle h_{1},h_{3}\right\rangle _{t}$,\hfill{}(4.8)

\medskip{}

\noindent for all $h_{1},h_{2},h_{3}\in\mathbb{H}_{t}$. Also, if
$h_{l}=\left(a_{l},b_{l}\right)\in\mathbb{H}_{t}$, for $l=1,2$,
and $r\in\mathbb{R}$, then

\medskip{}

$\;\;\;\;$$\left\langle rh_{1},h_{2}\right\rangle _{t}=\tau\left(\left(\left(r,0\right)\cdot_{t}h_{1}\right)\cdot_{t}h_{2}^{\dagger}\right)$

\medskip{}

$\;\;\;\;\;\;\;\;\;\;$$=\tau\left(\left(\left(\begin{array}{cc}
r & 0\\
0 & \overline{r}
\end{array}\right)\left(\begin{array}{cc}
a_{1} & tb_{1}\\
\overline{b_{1}} & \overline{a_{1}}
\end{array}\right)\right)\left(\begin{array}{cc}
a_{2} & tb_{2}\\
\overline{b_{2}} & \overline{a_{2}}
\end{array}\right)^{\dagger}\right)$

\medskip{}

$\;\;\;\;\;\;\;\;\;\;$$=\tau\left(\left(\begin{array}{cc}
ra_{1} & trb_{1}\\
\overline{rb_{1}} & \overline{ra_{1}}
\end{array}\right)\left(\begin{array}{cc}
\overline{a_{2}} & t\left(-b_{2}\right)\\
\overline{-b_{2}} & a_{2}
\end{array}\right)\right)$

\medskip{}

\noindent since $r=\overline{r}$ in $\mathbb{R}$

\medskip{}

$\;\;\;\;\;\;\;\;\;\;$$=\tau\left(\left(\begin{array}{ccc}
ra_{1}\overline{a_{2}}-trb_{1}\overline{b_{2}} &  & t\left(-ra_{1}b_{2}+ra_{2}b_{1}\right)\\
\\
\overline{-ra_{1}b_{2}+ra_{2}b_{1}} &  & \overline{ra_{1}\overline{a_{2}}-trb_{1}\overline{b_{2}}}
\end{array}\right)\right)$

\medskip{}

$\;\;\;\;\;\;\;\;\;\;$$=Re\left(ra_{1}\overline{a_{2}}-trb_{1}\overline{b_{2}}\right)=rRe\left(a_{1}\overline{a_{2}}-tb_{1}\overline{b_{2}}\right)$

\medskip{}

$\;\;\;\;\;\;\;\;\;\;$$=rtr\left(\left(a_{1},b_{1}\right)\cdot_{t}\left(a_{2},b_{2}\right)^{\dagger}\right)=r\left\langle h_{1},h_{2}\right\rangle _{t}$,

\noindent i.e.,\hfill{}(4.9) 
\[
\left\langle rh_{1},h_{2}\right\rangle _{t}=r\left\langle h_{1},h_{2}\right\rangle _{t},\;\;\forall r\in\mathbb{R}\;\mathrm{and\;}h_{1},h_{2}\in\mathbb{H}_{t}.
\]
Similarly, one obtains that

\medskip{}

\hfill{}$\left\langle h_{1},rh_{2}\right\rangle _{t}=r\left\langle h_{1},h_{2}\right\rangle _{t},\;\;\forall r\in\mathbb{R}\;\mathrm{and\;}h_{1},h_{2}\in\mathbb{H}_{t}.$\hfill{}(4.10)

\medskip{}

\begin{lemma}
The form $\left\langle ,\right\rangle _{t}$ of (4.6) is a well-defined
bilinear form on $\mathbb{H}_{t}$ ``over $\mathbb{R}$.'' 
\end{lemma}

\begin{proof}
Indeed, the form $\left\langle ,\right\rangle _{t}$ of (4.6) is a
bilinear form over $\mathbb{R}$, by (4.7), (4.8), (4.9) and (4.10). 
\end{proof}
By the above lemma, the $t$-scaled hypercomplex ring $\mathbb{H}_{t}$
is equipped with a well-defined bilinear form $\left\langle ,\right\rangle _{t}$
of (4.6) over $\mathbb{R}$. 
\begin{lemma}
If $h_{1},h_{2}\in\mathbb{H}_{t}$, then

\medskip{}

\hfill{}$\left\langle h_{1},h_{2}\right\rangle _{t}=\left\langle h_{2},h_{1}\right\rangle _{t}$
in $\mathbb{R}$.\hfill{}(4.11) 
\end{lemma}

\begin{proof}
Let $h_{l}=\left(a_{l},b_{l}\right)\in\mathbb{H}_{t}$, for $l=1,2$.
Then

\medskip{}

$\;\;\;$$\left\langle h_{1},h_{2}\right\rangle _{t}=\tau\left(\left[h_{1}\cdot_{t}h_{2}^{\dagger}\right]_{t}\right)=\tau\left(\left(\begin{array}{cc}
a_{1} & tb_{1}\\
\overline{b_{1}} & \overline{a_{1}}
\end{array}\right)\left(\begin{array}{cc}
\overline{a_{2}} & t\left(-b_{2}\right)\\
-\overline{b_{2}} & a_{2}
\end{array}\right)\right)$

\medskip{}

$\;\;\;\;\;\;\;\;\;\;\;\;\;\;\;\;$$=\tau\left(\left(\begin{array}{ccc}
a_{1}\overline{a_{2}}-tb_{1}\overline{b_{2}} &  & t\left(a_{2}b_{1}-a_{1}b_{2}\right)\\
\\
\overline{a_{2}b_{1}-a_{1}b_{2}} &  & \overline{a_{1}\overline{a_{2}}-tb_{1}\overline{b_{2}}}
\end{array}\right)\right)$

\medskip{}

$\;\;\;\;\;\;\;\;\;\;\;\;\;\;\;\;$$=Re\left(a_{1}\overline{a_{2}}-tb_{1}\overline{b_{2}}\right)$,

\noindent and\hfill{}(4.12)

$\;\;\;$$\left\langle h_{2},h_{1}\right\rangle _{t}=\tau\left(\left[h_{2}\cdot_{t}h_{1}^{\dagger}\right]\right)=\tau\left(\left(\begin{array}{cc}
a_{2} & tb_{2}\\
\overline{b_{2}} & \overline{a_{2}}
\end{array}\right)\left(\begin{array}{cc}
\overline{a_{1}} & t\left(-b_{1}\right)\\
-\overline{b_{1}} & a_{1}
\end{array}\right)\right)$

\medskip{}

$\;\;\;\;\;\;\;\;\;\;\;\;\;\;\;\;$$=\tau\left(\left(\begin{array}{ccc}
\overline{a_{1}}a_{2}-t\overline{b_{1}}b_{2} &  & t\left(a_{1}b_{2}-a_{2}b_{1}\right)\\
\\
\overline{a_{1}b_{2}-a_{2}b_{1}} &  & \overline{\overline{a_{1}}a_{2}-t\overline{b_{1}}b_{2}}
\end{array}\right)\right)$

\medskip{}

$\;\;\;\;\;\;\;\;\;\;\;\;\;\;\;\;$$=Re\left(\overline{a_{1}}a_{2}-t\overline{b_{1}}b_{2}\right)=Re\left(\overline{\overline{a_{1}}a_{2}-t\overline{b_{1}}b_{2}}\right)=\overline{\left\langle h_{1},h_{2}\right\rangle _{t}},$

\medskip{}

\noindent implying that 
\[
\left\langle h_{2},h_{1}\right\rangle _{t}=\overline{\left\langle h_{1},h_{2}\right\rangle _{t}}=\left\langle h_{1},h_{2}\right\rangle _{t},
\]
in $\mathbb{R}$. In particular, the second equality of the above
formula is satisfied, since $\left\langle h_{1},h_{2}\right\rangle _{t}\in\mathbb{R}$
by (4.6), for all $h_{1},h_{2}\in\mathbb{H}_{t}$. Therefore, the
symmetry (4.11) of $\left\langle ,\right\rangle _{t}$ holds by (4.12).

Independently, the symmetry (4.11) is obtained, because the linear
functional $tr$ is a trace on $\mathcal{H}_{2}^{t}$ satisfying 
\[
tr\left(AB\right)=tr\left(BA\right),\;\forall A,B\in\mathcal{H}_{2}^{t},
\]
implying that 
\[
\tau\left(AB\right)=\tau\left(BA\right),\;\forall A,B\in\mathcal{H}_{2}^{t},
\]
and hence, 
\[
\tau\left(h_{1}\cdot_{t}h_{2}\right)=\tau\left(h_{2}\cdot_{t}h_{1}\right),\;\forall h_{1},h_{2}\in\mathbb{H}_{t},
\]
in $\mathbb{R}$. Therefore, one obtains the symmetry (4.11). 
\end{proof}
The above lemma shows our bilinear form $\left\langle ,\right\rangle _{t}$
of (4.6) is symmetric by (4.11). 
\begin{lemma}
If $h_{1},h_{2}\in\mathbb{H}_{t}$, then

\medskip{}

\hfill{}$\left|\left\langle h_{1},h_{2}\right\rangle _{2}\right|^{2}\leq\left|\left\langle h_{1},h_{1}\right\rangle _{t}\right|^{2}\left|\left\langle h_{2},h_{2}\right\rangle _{t}\right|^{2}$,\hfill{}(4.13)

\medskip{}

\noindent where $\left|.\right|$ is the absolute value on $\mathbb{R}$. 
\end{lemma}

\begin{proof}
By (4.12), if $h_{l}=\left(a_{l},b_{l}\right)\in\mathbb{H}_{t}$ for
$l=1,2$, then one has 
\[
\left|\left\langle h_{1},h_{2}\right\rangle _{t}\right|=\left|Re\left(a_{1}\overline{a_{2}}-tb_{1}\overline{b_{2}}\right)\right|,
\]
and, similarly, 
\[
\left|\left\langle h_{l},h_{l}\right\rangle _{t}\right|=\left|\left|a_{l}\right|^{2}-t\left|b_{l}\right|^{2}\right|,
\]
for $l=1,2$. Therefore, the inequality (4.13) holds. 
\end{proof}
Observe now that if $h=\left(a,b\right)\in\mathbb{H}_{t}$, then 
\[
\left\langle h,h\right\rangle _{t}=\tau\left(\left(a,b\right)\cdot_{t}\left(a,b\right)^{\dagger}\right)=Re\left(\left|a\right|^{2}-t\left|b\right|^{2}\right),
\]
by (4.1) and (4.4), implying that

\medskip{}

\hfill{}$\left\langle h,h\right\rangle _{t}=\left|a\right|^{2}-t\left|b\right|^{2}=det\left(\left[h\right]_{t}\right).$\hfill{}(4.14) 
\begin{lemma}
If $h\in\mathbb{H}_{t}$, then $\left\langle h,h\right\rangle _{t}=det\left(\left[h\right]_{t}\right)$. 
\end{lemma}

\begin{proof}
It is shown by (4.14). 
\end{proof}
The above lemma shows that the bilinear form $\left\langle ,\right\rangle _{t}$
of (4.6) is not positively defined in general. i.e., 
\[
\left\langle h,h\right\rangle _{t}=det\left(\left[h\right]_{t}\right)\in\mathbb{R}.
\]
Observe that, by (4.14), 
\[
\left\langle h,h\right\rangle _{t}=0\Longleftrightarrow det\left(\left[h\right]_{t}\right)=0\;\mathrm{in\;}\mathbb{R},
\]
if and only if\hfill{}(4.15) 
\[
h=\left(a,b\right)\in\mathbb{H}_{t},\;\mathrm{with\;}\left|a\right|^{2}=t\left|b\right|^{2}.
\]

\begin{lemma}
Let $h=\left(a,b\right)\in\mathbb{H}_{t}$, and $\left\langle ,\right\rangle _{t}$,
the bilinear form (4.6). Then

\medskip{}

\hfill{}$\left\langle h,h\right\rangle _{t}=0,\Longleftrightarrow\left|a\right|^{2}=t\left|b\right|^{2}.$\hfill{}(4.16) 
\end{lemma}

\begin{proof}
The relation (4.16) is shown by (4.14) and (4.15). 
\end{proof}
Define now a subset $Q_{t}$ on $\mathbb{H}_{t}$ by

\medskip{}

\hfill{}$Q_{t}=\left\{ h\in\mathbb{H}_{t}:\left\langle h,h\right\rangle _{t}=0\right\} $.\hfill{}(4.17)

\medskip{}

\begin{proposition}
Let $Q_{t}$ be the subset (4.17) of our $t$-scaled hypercomplex
ring $\mathbb{H}_{t}$. Then 
\[
Q_{t}=\left\{ h\in\mathbb{H}_{t}:det\left(\left[h\right]_{t}\right)=0\right\} ,
\]
equivalently,\hfill{}(4.18) 
\[
Q_{t}=\left\{ h\in\mathbb{H}_{t}:h\textrm{ is not invertible in }\mathbb{H}_{t}\right\} ,
\]
equivalently,\hfill{}(4.19) 
\[
Q_{t}=\mathbb{H}_{t}^{sing}=\mathbb{H}_{t}^{\times sing}\cup\left\{ \left(0,0\right)\right\} ,
\]
where $\mathbb{H}_{t}^{\times sing}$ is the semigroup-part of $\mathbb{H}_{t}$. 
\end{proposition}

\begin{proof}
By (4.14), (4.16) and (4.17), the first set-equality of (4.18) holds.
It implies the second set-equality of (4.18) by (2.2.2). Therefore,
the set-equality (4.19) holds by (2.2.4) and (2.2.5). 
\end{proof}
The above proposition shows that the subset $Q_{t}$ of $\mathbb{H}_{t}$,
in the sense of (4.17) induced by the bilinear form $\left\langle ,\right\rangle _{t}$
of (4.6), is nothing but the block $\mathbb{H}_{t}^{sing}$ of $\mathbb{H}_{t}$,
in the sense of (2.2.4). Thus, this proposition illustrates connections
between the invertibility on $\mathbb{H}_{t}$, and the bilinear form
$\left\langle ,\right\rangle _{t}$ on $\mathbb{H}_{t}$.

\medskip{}

\noindent $\mathbf{Notation.}$ From below, we denote the subset $Q_{t}$
of (4.17) by $\mathbb{H}_{t}^{sing}$, by (4.18) and (4.19). i.e.,
\[
\mathbb{H}_{t}^{sing}=\left\{ \left(a,b\right)\in\mathbb{H}_{t}:\left|a\right|^{2}=t\left|b\right|^{2}\right\} =\left\{ h\in\mathbb{H}_{t}:\left\langle h,h\right\rangle _{t}=0\right\} ,
\]
in $\mathbb{H}_{t}$.\hfill{}{$\square$}

\medskip{}

Consider now that if $h_{1},h_{2}\in\mathbb{H}_{t}$, and if either
$h_{1}$ or $h_{2}$ is contained in $\mathbb{H}_{t}^{sing}$, then
$h_{1}\cdot_{t}h_{2}^{\dagger}\in\mathbb{H}_{t}^{sing}$, by (2.2.6).
Indeed, if either $h_{1}$ or $h_{2}$ is in $\mathbb{H}_{t}^{sing}$,
then 
\[
det\left(\left[h_{1}\cdot_{t}h_{2}^{\dagger}\right]_{t}\right)=det\left(\left[h_{1}\right]_{t}\left[h_{2}\right]_{t}^{\dagger}\right),
\]
implying that\hfill{}(4.20) 
\[
det\left(\left[h_{1}\cdot_{t}h_{2}^{\dagger}\right]_{t}\right)=\left(det\left(\left[h_{1}\right]_{t}\right)\right)\left(det\left(\left[h_{2}\right]_{t}^{\dagger}\right)\right)=0,
\]
showing that $h_{1}\cdot_{t}h_{2}^{\dagger}\in\mathbb{H}_{t}^{sing}$,
too, in $\mathbb{H}_{t}$. 
\begin{definition}
For a vector space $X$ over $\mathbb{R}$, a form $\left\langle ,\right\rangle :X\times X\rightarrow\mathbb{R}$
is called a (definite) semi-inner product on $X$ over $\mathbb{R}$,
if (i) it is a bilinear form on $X$ over $\mathbb{R}$, (ii) 
\[
\left\langle x_{1},x_{2}\right\rangle =\left\langle x_{2},x_{1}\right\rangle ,\;\forall x_{1},x_{2}\in X,
\]
and (iii) $\left\langle x,x\right\rangle \geq0$, for all $x\in X$.
If such a semi-inner product $\left\langle ,\right\rangle $ satisfies
an additional condition (iv) 
\[
\left\langle x,x\right\rangle =0,\textrm{ if and only if }x=0_{X},
\]
where $0_{X}$ is the zero vector of $X$, then it is called an inner
product on $X$ over $\mathbb{R}$. If $\left\langle ,\right\rangle $
is a semi-inner product (or, an inner product) on the $\mathbb{R}$-vector
space $X$, then the pair $\left(X,\left\langle ,\right\rangle \right)$
is said to be a semi-inner product space (respectively, an inner product
space) over $\mathbb{R}$ (in short, a $\mathbb{R}$-SIPS, respectively,
a $\mathbb{R}$-IPS). 
\end{definition}

By definition, all $\mathbb{R}$-IPSs are automatically $\mathbb{R}$-SIPSs.
But not all $\mathbb{R}$-SIPSs are $\mathbb{R}$-IPSs. 
\begin{definition}
For a vector space $X$ over $\mathbb{R}$, a form $\left\langle ,\right\rangle :X\times X\rightarrow\mathbb{R}$
is called an indefinite semi-inner product on $X$ over $\mathbb{R}$,
if (i) it is a bilinear form on $X$ over $\mathbb{R}$, (ii) 
\[
\left\langle x_{1},x_{2}\right\rangle =\left\langle x_{2},x_{1}\right\rangle ,\;\;\forall x_{1},x_{2}\in X,
\]
and (iii) $\left\langle x,x\right\rangle \in\mathbb{R}$, for all
$x\in X$. If such an indefinite semi-inner product $\left\langle ,\right\rangle $
satisfies an additional condition (iv) 
\[
\left\langle x,x\right\rangle =0,\textrm{ if and only if }x=0_{X},
\]
then it is said to be an indefinite inner product on $X$ over $\mathbb{R}$.
If $\left\langle ,\right\rangle $ is an indefinite semi-inner product
(or, an indefinite inner product) on the $\mathbb{R}$-vector space
$X$, then the pair $\left(X,\left\langle ,\right\rangle \right)$
is called an indefinite-semi-inner product space (respectively, an
indefinite-inner product space) over $\mathbb{R}$ (in short, a $\mathbb{R}$-ISIPS,
respectively, $\mathbb{R}$-IIPS). 
\end{definition}

By definition, every $\mathbb{R}$-IIPS is an $\mathbb{R}$-ISIPS,
however, not all $\mathbb{R}$-ISIPSs are $\mathbb{R}$-IIPSs. With
help of the above two definitions, we obtain the following result.
Remember that, set-theoretically, $\mathbb{H}_{t}=\mathbb{C}^{2}$,
implying that the set $\mathbb{H}_{t}$ is clearly a vector space
over $\mathbb{R}$, which is (isomorphic to) a subspace of 4-dimensional
$\mathbb{R}$-vector space $\mathbb{R}^{4}$. 
\begin{theorem}
If $t<0$ in $\mathbb{R}$, then the bilinear form $\left\langle ,\right\rangle _{t}$
of (4.6) is a continuous semi-inner product on $\mathbb{H}_{t}$ over
$\mathbb{R}$, i.e.,

\medskip{}

\hfill{}$t<0\Longrightarrow\left\langle ,\right\rangle _{t}\textrm{ is an inner product on }\mathbb{H}_{t}.$\hfill{}(4.21)

\medskip{}

\noindent Meanwhile, if $t\geq0$, then $\left\langle ,\right\rangle _{t}$
is a continuous indefinite semi-inner product on $\mathbb{H}_{t}$
over $\mathbb{R}$, i.e.,

\medskip{}

\hfill{}$t\geq0\Longrightarrow\left\langle ,\right\rangle _{t}\textrm{ is an indefinite semi-inner product on }\mathbb{H}_{t}.$\hfill{}(4.22) 
\end{theorem}

\begin{proof}
Recall first that, by the partition (2.2.4), if $t<0$ in $\mathbb{R}$,
then the semigroup-part $\mathbb{H}_{t}^{\times sing}=\mathbb{H}_{t}^{sing}\setminus\left\{ \left(0,0\right)\right\} $
is empty in $\mathbb{H}_{t}$, and hence, 
\[
\mathbb{H}_{t}^{sing}=\left\{ \left(0,0\right)\right\} =Q_{t},
\]
in the sense of (4.17), and 
\[
\mathbb{H}_{t}=\mathbb{H}_{t}^{inv}\cup\left\{ \left(0,0\right)\right\} .
\]
It implies that 
\[
\left\langle h,h\right\rangle _{t}=0\Longleftrightarrow h=\left(0,0\right)\in\mathbb{H}_{t},
\]
whenever $t<0$. Moreover, for any $h=\left(a,b\right)\in\mathbb{H}_{t}$,
\[
det\left(\left[h\right]_{t}\right)=\left|a\right|^{2}-t\left|b\right|^{2}=\tau\left(\left[h\cdot_{t}h^{\dagger}\right]\right)=\left\langle h,h\right\rangle _{t},
\]
satisfying\hfill{}(4.23) 
\[
\left\langle h,h\right\rangle _{t}\geq0,\;\;\mathrm{whenever\;\;}t<0.
\]
Therefore, if $t<0$, then the form $\left\langle ,\right\rangle _{t}$
is an inner product on $\mathbb{H}_{t}$ over $\mathbb{R}$, by (4.7),
(4.8), (4.9), (4.10), (4.11), (4.16), (4.19), and (4.23). Especially,
the continuity of $\left\langle ,\right\rangle _{t}$ is guaranteed
by (4.13).

Assume now that $t\geq0$ in $\mathbb{R}$. Then the semigroup-part
$\mathbb{H}_{t}^{\times sing}$ is not empty in $\mathbb{H}_{t}$,
and hence, 
\[
\mathbb{H}_{t}^{sing}\supset\left\{ \left(0,0\right)\right\} \;\;\mathrm{in\;\;}\mathbb{H}_{t}.
\]
Furthermore, by (4.16) and (4.20), the above non-negativity (4.23)
does not hold true on $\mathbb{H}_{t}$, whenever $t\geq0$. However,
by (4.7), (4.8), (4.9), (4.10), (4.11), (4.16) and (4.19), it becomes
a well-defined indefinite semi-inner product on $\mathbb{H}_{t}$
over $\mathbb{R}$, for $t\geq0$. The continuity of $\left\langle ,\right\rangle _{t}$
is obtained by (4.13). 
\end{proof}
The following corollary is an immediate consequence of the above theorem. 
\begin{corollary}
If $t<0$, then the pair $\left(\mathbb{H}_{t},\left\langle ,\right\rangle _{t}\right)$
is a $\mathbb{R}$-IPS, meanwhile, if $t\geq0$, then $\left(\mathbb{H}_{t},\left\langle ,\right\rangle _{t}\right)$
is a $\mathbb{R}$-ISIPS. 
\end{corollary}

\begin{proof}
It is proven by (4.21) for $t<0$, and by (4.22) for $t\geq0$. 
\end{proof}
The above theorem and corollary show that $\left\{ \mathbb{H}_{t}\right\} _{t<0}$
are $\mathbb{R}$-SIPSs, while, $\left\{ \mathbb{H}_{t}\right\} _{t\geq0}$
form $\mathbb{R}$-ISIPSs, and such definite, or indefinite inner-product
processes are continuous by (4.13).

Recall that a pair $\left(X,\left\Vert .\right\Vert \right)$ of a
vector space $X$ over $\mathbb{R}$, and a map $\left\Vert .\right\Vert :X\rightarrow\mathbb{R}$
is called a semi-normed space, if $\left\Vert .\right\Vert $ is a
semi-norm, satisfying the three conditions; (i) $\left\Vert x\right\Vert \geq0$,
for all $x\in X$, (ii) $\left\Vert rx\right\Vert =\left|r\right|\left\Vert x\right\Vert $,
for all $r\in\mathbb{R}$ and $x\in X$, and (iii) 
\[
\left\Vert x_{1}+x_{2}\right\Vert \leq\left\Vert x_{1}\right\Vert +\left\Vert x_{2}\right\Vert ,\;\forall x_{1},x_{2}\in X.
\]
If the semi-norm $\left\Vert .\right\Vert $ satisfies an additional
condition (iv) 
\[
\left\Vert x\right\Vert =0\Longleftrightarrow x=0_{X}\;\;\mathrm{in\;\;}X,
\]
then it is said to be a norm on $X$. In such a case, the semi-normed
space $\left(X,\left\Vert .\right\Vert \right)$ is called a normed
space over $\mathbb{R}$. 
\begin{definition}
If a pair $\left(X,\left\Vert .\right\Vert _{X}\right)$ of a vector
space $X$ over $\mathbb{R}$, and a semi-norm (respectively, a norm)
$\left\Vert .\right\Vert $ on $X$, is complete under its semi-norm
topology (respectively, norm topology) induced by $\left\Vert .\right\Vert $,
then it is said to be a complete semi-normed space (respectively,
a Banach space) over $\mathbb{R}$, in short, a complete $\mathbb{R}$-SNS
(respectively, $\mathbb{R}$-Banach space). 
\end{definition}

Suppose $t<0$, and hence, the pair $\left(\mathbb{H}_{t},\;\left\langle ,\right\rangle _{t}\right)$
is a $\mathbb{R}$-IPS. Then, as usual, the semi-inner product $\left\langle ,\right\rangle _{t}$
induces the corresponding semi-norm, 
\[
\left\Vert .\right\Vert _{t}:\mathbb{H}_{t}\rightarrow\mathbb{R},
\]
by\hfill{}(4.24) 
\[
\left\Vert h\right\Vert _{t}\overset{\textrm{def}}{=}\sqrt{\left\langle h,h\right\rangle _{t}}=\sqrt{det\left(\left[h\right]_{t}\right)}=\sqrt{\tau\left(\left[h\cdot_{t}h^{\dagger}\right]\right)},
\]
for all $h\in\mathbb{H}_{t}$. Indeed, if $t<0$, then the map $\left\Vert .\right\Vert _{t}$
of (4.24) satisfies that 
\[
\left\Vert h\right\Vert _{t}\geq0,\;\mathrm{since\;}\left\langle h,h\right\rangle _{t}=\left|a\right|^{2}-t\left|b\right|^{2}\geq0,
\]
and 
\[
\left\Vert rh\right\Vert _{t}=\sqrt{\left\langle rh,rh\right\rangle _{t}}=\sqrt{r^{2}}\sqrt{\left\langle h,h\right\rangle _{t}}=\left|r\right|\left\Vert h\right\Vert _{t},
\]
for all $r\in\mathbb{R}$ and $h\in\mathbb{H}_{t}$, and 
\[
\left\Vert h_{1}+h_{2}\right\Vert _{t}^{2}=\left\langle h_{1}+h_{2},h_{1}+h_{2}\right\rangle _{t}\leq\left(\left\Vert h_{1}\right\Vert _{t}+\left\Vert h_{2}\right\Vert _{t}\right)^{2},
\]
implying that 
\[
\left\Vert h_{1}+h_{2}\right\Vert _{t}\leq\left\Vert h_{1}\right\Vert _{t}+\left\Vert h_{2}\right\Vert _{t},\;\mathrm{in\;}\mathbb{R},
\]
for all $h_{1},h_{2}\in\mathbb{H}_{t}$.By the above definition, one
obtains the following result. 
\begin{theorem}
Let $t<0$ in $\mathbb{R}$, and $\left(\mathbb{H}_{t},\left\langle ,\right\rangle _{t}\right)$,
the corresponding $\mathbb{R}$-IPS. Then the map $\left\Vert .\right\Vert _{t}$
of (4.24) is a norm on $\mathbb{H}_{t}$ over $\mathbb{R}$. Furthermore,
$\mathbb{H}_{t}$ is complete under $\left\Vert .\right\Vert _{t}$.
i.e.,

\medskip{}

\hfill{}$t<0$$\;$$\Longrightarrow$$\;$$\left(\mathbb{H}_{t},\left\Vert .\right\Vert _{t}\right)$
is a $\mathbb{R}$-Banach space.\hfill{}(4.25)

\medskip{}

\noindent Equivalently, if $t<0$, then the $\mathbb{R}$-IPS $\left(\mathbb{H}_{t},\left\langle ,\right\rangle _{t}\right)$
is a Hilbert space over $\mathbb{R}$ (or, a $\mathbb{R}$-Hilbert
space). 
\end{theorem}

\begin{proof}
As we discussed in the very above paragraph, if a given scale $t$
is negative in $\mathbb{R}$, then the map $\left\Vert .\right\Vert _{t}$
induced by the semi-inner product $\left\langle ,\right\rangle _{t}$
forms a well-defined semi-norm. So, the pair $\left(\mathbb{H}_{t},\left\Vert .\right\Vert _{t}\right)$
is a $\mathbb{R}$-SNS. Moreover, if $t<0$, then 
\[
\left\Vert h\right\Vert _{t}=0\Longleftrightarrow\left\langle h,h\right\rangle _{t}=0\Longleftrightarrow h=\left(0,0\right),
\]
by (4.19), because $\mathbb{H}_{t}^{sing}=\left\{ \left(0,0\right)\right\} $
whenever $t<0$. So, the semi-norm $\left\Vert .\right\Vert _{t}$
is a norm on $\mathbb{H}_{t}$ over $\mathbb{R}$. Since $\mathbb{H}_{t}$
is a subspace of a finite-dimensional $\mathbb{R}$-vector space $\mathbb{R}^{4}$,
this norm $\left\Vert .\right\Vert _{t}$ is complete over $\mathbb{R}$,
by (4.13) and (4.24), i.e., $\left(\mathbb{H}_{t},\left\Vert .\right\Vert _{t}\right)$
is a complete $\mathbb{R}$-NS, equivalently, it is a $\mathbb{R}$-Banach
space. Therefore, the statement (4.25) holds true.

Note that every $\mathbb{R}$-IPS $\left(X,\left\langle ,\right\rangle \right)$
induces the corresponding normed space $\left(X,\left\Vert .\right\Vert _{X}\right)$
over $\mathbb{R}$, where 
\[
\left\Vert x\right\Vert _{X}=\sqrt{\left\langle x,x\right\rangle },\;\;\;\forall x\in X.
\]
If this $\mathbb{R}$-normed space $\left(X,\left\Vert .\right\Vert _{X}\right)$
forms a $\mathbb{R}$-Banach space, then the $\mathbb{R}$-IPS, $\left(X,\left\langle ,\right\rangle \right)$,
is said to be a $\mathbb{R}$-Hilbert space. Thus, our $\mathbb{R}$-IPS,
$\left(\mathbb{H}_{t},\left\langle ,\right\rangle _{t}\right)$, forms
a $\mathbb{R}$-Hilbert space by (4.25). 
\end{proof}
The above theorem shows that our $\mathbb{R}$-IPSs $\left\{ \left(\mathbb{H}_{t},\left\langle ,\right\rangle _{t}\right)\right\} _{t<0}$
induce the corresponding $\mathbb{R}$-Banach spaces $\left\{ \left(\mathbb{H}_{t},\left\Vert .\right\Vert _{t}\right)\right\} _{t<0}$,
where $\left\{ \left\Vert .\right\Vert _{t}\right\} _{t<0}$ are in
the sense of (4.24).

Now, assume that $t\geq0$, and hence, the pair $\left(\mathbb{H}_{t},\left\langle ,\right\rangle _{t}\right)$
is a $\mathbb{R}$-ISIPS. Now, we define a map, also denoted by $\left\Vert .\right\Vert _{t}$,
\[
\left\Vert .\right\Vert _{t}:\mathbb{H}_{t}\rightarrow\mathbb{R},
\]
by\hfill{}(4.26) 
\[
\left\Vert h\right\Vert _{t}\overset{\textrm{def}}{=}\sqrt{\left|\left\langle h,h\right\rangle _{t}\right|}=\sqrt{\left|\tau\left(h\cdot_{t}h^{\dagger}\right)\right|}=\sqrt{\left|det\left(\left[h\right]_{t}\right)\right|},
\]
for all $h\in\mathbb{H}_{t}$, where $\left|.\right|$ in (4.26) is
the absolute value on $\mathbb{R}$. Then, by the very definition
(4.26), one obtains the similar result like (4.25). 
\begin{theorem}
Let $t\geq0$, and $\left(\mathbb{H}_{t},\left\langle ,\right\rangle _{t}\right)$,
the corresponding $\mathbb{R}$-ISIPS. Then the map $\left\Vert .\right\Vert _{t}$
of (4.26) is a well-defined complete semi-norm over $\mathbb{R}$.
i.e.,

\medskip{}

\hfill{}$t\geq0\Longrightarrow$$\left(\mathbb{H}_{t},\left\Vert .\right\Vert _{t}\right)$
is a complete $\mathbb{R}$-SNS.\hfill{}(4.27) 
\end{theorem}

\begin{proof}
Assume that a given scale $t$ is non-negative in $\mathbb{R}$. For
any $h\in\mathbb{H}_{t}$, 
\[
\left\Vert h\right\Vert _{t}=\sqrt{\left|\left\langle h,h\right\rangle _{t}\right|}\geq0,
\]
since $\left|\left\langle h,h\right\rangle _{t}\right|\geq0$, for
all $h\in\mathbb{H}_{t}$. Also, for any $r\in\mathbb{R}$ and $h\in\mathbb{H}_{t}$,
\[
\left\Vert rh\right\Vert _{t}=\sqrt{\left|\left\langle rh,rh\right\rangle _{t}\right|}=\sqrt{\left|r^{2}\right|}\sqrt{\left|\left\langle h,h\right\rangle _{t}\right|}=\left|r\right|\left\Vert h\right\Vert _{t}.
\]
And, we have that 
\[
\left\Vert h_{1}+h_{2}\right\Vert _{t}=\sqrt{\left|\left\langle h_{1}+h_{2},h_{1}+h_{2}\right\rangle _{t}\right|},
\]
and 
\[
\left|\left\langle h_{1}+h_{2},h_{1}+h_{2}\right\rangle _{t}\right|^{2}\leq\left(\left|\left\langle h_{1},h_{1}\right\rangle _{t}\right|+\left|\left\langle h_{2},h_{2}\right\rangle _{t}\right|\right)^{2},
\]
implying that 
\[
\left\Vert h_{1}+h_{2}\right\Vert _{t}\leq\left\Vert h_{1}\right\Vert _{t}+\left\Vert h_{2}\right\Vert _{t}\;\;\mathrm{in\;\;}\mathbb{R},
\]
for all $h_{1},h_{2}\in\mathbb{H}_{t}$. Moreover, 
\[
\left\Vert h\right\Vert _{t}=0\Longleftrightarrow\left\langle h,h\right\rangle _{t}=0\Longleftrightarrow h\in\mathbb{H}_{t}^{sing},
\]
by (4.19), and the semigroup-part $\mathbb{H}_{t}^{\times sing}=\mathbb{H}_{t}^{sing}\setminus\left\{ \left(0,0\right)\right\} $
is not empty in the $t$-scaled monoid $\mathbb{H}_{t}^{\times}$.
So, the map $\left\Vert .\right\Vert _{t}$ of (4.26) is a well-defined
semi-norm on $\mathbb{H}_{t}$, which is not a norm, over $\mathbb{R}$,
and hence, the pair $\left(\mathbb{H}_{t},\left\Vert .\right\Vert _{t}\right)$
forms a $\mathbb{R}$-SNS. The completeness is similarly shown as
in the proof of (4.25), by (4.13) and (4.26). Therefore, the statement
(4.27) holds true. 
\end{proof}
By (4.25) and (4.27), we obtain the following corollary. 
\begin{corollary}
For any $t\in\mathbb{R}$, let $\left\Vert .\right\Vert _{t}:\mathbb{H}_{t}\rightarrow\mathbb{R}$
be the map, defined by

\medskip{}

\hfill{}$\left\Vert h\right\Vert _{t}\overset{\textrm{def}}{=}\sqrt{\left|\left\langle h,h\right\rangle _{t}\right|},\;\;\;\forall h\in\mathbb{H}_{t}.$\hfill{}(4.28)

\medskip{}

\noindent Then the pair $\left(\mathbb{H}_{t},\left\Vert .\right\Vert _{t}\right)$
forms a complete $\mathbb{R}$-SNS. Especially, if $t<0$, then it
is a $\mathbb{R}$-Banach space, meanwhile, if $t\geq0$, then it
forms a complete $\mathbb{R}$-SNS. 
\end{corollary}

\begin{proof}
By the very definition (4.28), this map $\left\Vert .\right\Vert _{t}$
is identical to the semi-norm (4.24), if $t<0$, meanwhile, it is
identified with the semi-norm (4.26), if $t\geq0$. Therefore, by
(4.25) and (4.27), the corresponding pair $\left(\mathbb{H}_{t},\left\Vert .\right\Vert _{t}\right)$
is a complete $\mathbb{R}$-SNS, for all $t\in\mathbb{R}$. 
\end{proof}
The above corollary shows that a bilinear form $\left\langle ,\right\rangle _{t}$
of (4.6) induces the semi-norm $\left\Vert .\right\Vert _{t}$ of
(4.28), making the pair $\left(\mathbb{H}_{t},\left\Vert .\right\Vert _{t}\right)$
be a complete $\mathbb{R}$-SNS, for all $t\in\mathbb{R}$.\medskip{}

\noindent $\mathbf{Notation.}$ From below, we denote a $\mathbb{R}$-Hilbert
space $\left(\mathbb{H}_{t},\left\langle ,\right\rangle _{t}\right)$
by $\mathbf{H}_{t}$, whenever $t<0$. And, by $\mathbf{K}_{t}$,
we denote a complete $\mathbb{R}$-ISIPS $\left(\mathbb{H}_{t},\left\langle ,\right\rangle _{t}\right)$,
whenever $t\geq0$, where the completeness is up to the complete $\mathbb{R}$-SNS
$\left(\mathbb{H}_{t},\left\Vert .\right\Vert _{t}\right)$. Universally,
we write $\mathbf{X}_{t}$ as either $\mathbf{H}_{t}$, or $\mathbf{K}_{t}$,
for $t\in\mathbb{R}$. For convenience, we call $\mathbf{X}_{t}\in\left\{ \mathbf{H}_{t},\mathbf{K}_{t}\right\} $,
the $t$(-scaled)-hypercomplex (semi-normed) $\mathbb{R}$(-vector-)space
(over $\mathbb{R}$).\hfill{}{$\square$}

\medskip{}

By the above notation, $\left\{ \mathbf{X}_{t}=\mathbf{H}_{t}\right\} _{t<0}$
are $\mathbb{R}$-Hilbert spaces, while $\left\{ \mathbf{X}_{t}=\mathbf{K}_{t}\right\} _{t\geq0}$
are complete $\mathbb{R}$-ISIPSs.

\section{Certain Operators on $\mathbf{X}_{t}=\left(\mathbb{H}_{t},\left\langle ,\right\rangle _{t}\right)$}

Let $t\in\mathbb{R}$ be an arbitrary scale, and $\mathbb{H}_{t}$,
the $t$-scaled hypercomplex ring, understood as the hypercomplex
$\mathbb{R}$-space $\mathbf{X}_{t}=\left(\mathbb{H}_{t},\left\langle ,\right\rangle _{t}\right)$
over $\mathbb{R}$, where 
\[
\left\langle h_{1},h_{2}\right\rangle _{t}\overset{\textrm{def}}{=}\tau\left(\left[h_{1}\cdot_{t}h_{2}^{\dagger}\right]\right),\;\;\forall h_{1},h_{2}\in\mathbb{H}_{t},
\]
is a bilinear form on $\mathbb{H}_{t}$ over $\mathbb{R}$, which
is either a semi-Hilbert space $\mathbf{H}_{t}$ (if $t<0$), or a
semi-Krein space $\mathbf{K}_{t}$ (if $t\geq0$) over $\mathbb{R}$,
equipped with the semi-norm, 
\[
\left\Vert h\right\Vert _{t}\overset{\textrm{def}}{=}\sqrt{\left|\left\langle h,h\right\rangle _{t}\right|}=\sqrt{\left|det\left(\left[h\right]_{t}\right)\right|},
\]
in the sense of (4.28), for all $h\in\mathbb{H}_{t}$, satisfying
that: $h\in\mathbb{H}_{t}^{sing}$, if and only if $\left\Vert h\right\Vert _{t}=0$.

By the construction of $\mathbf{X}_{t}\in\left\{ \mathbf{H}_{t},\mathbf{K}_{t}\right\} $,
one can understand each hypercomplex number $h\in\mathbb{H}_{t}$
as an operator acting on this $t$-hypercomplex space $\mathbf{X}_{t}$.
Indeed, one can define a linear transformation, 
\[
M_{\eta}:\mathbf{X}_{t}\rightarrow\mathbf{X}_{t},\;\;\;\forall\eta\in\mathbb{H}_{t},
\]
by\hfill{}(5.1) 
\[
M_{\eta}\left(h\right)\overset{\textrm{def}}{=}\eta\cdot_{t}h,\;\;\;\forall h\in\mathbf{X}_{t},
\]
i.e., this morphism $M_{\eta}$ is a multiplication operator on $\mathbf{X}_{t}$
with its symbol $\eta\in\mathbb{H}_{t}$. By (4.13) and (4.28), this
operator $M_{\eta}$ of (5.1) is bounded (equivalently, continuous
under linearity) on $\mathbf{X}_{t}$. 
\begin{proposition}
Every hypercomplex number $\eta\in\mathbb{H}_{t}$ is a well-defined
bounded operator (over $\mathbb{R}$) on the $t$-hypercomplex space
$\mathbf{X}_{t}\in\left\{ \mathbf{H}_{t},\mathbf{K}_{t}\right\} $,
for all $t\in\mathbb{R}$. 
\end{proposition}

\begin{proof}
The proof is done by construction. See the definition (5.1) of the
bounded multiplication operators $\left\{ M_{\eta}\right\} _{\eta\in\mathbb{H}_{t}}$. 
\end{proof}
By the above proposition, one can realize that each hypercomplex number
$\eta\in\mathbb{H}_{t}$ is acting on the $t$-hypercomplex $\mathbb{R}$-space
$\mathbf{X}_{t}\in\left\{ \mathbf{H}_{t},\mathbf{K}_{t}\right\} $,
as a bounded operator (5.1). More precisely, if $\mathbf{X}_{t}=\mathbf{H}_{t}$,
equivalently, if $t<0$, then $\eta$ is acting on the semi-Hilbert
space $\mathbf{H}_{t}$, meanwhile, if $\mathbf{X}_{t}=\mathbf{K}_{t}$,
equivalently, if $t\geq0$, then it is acting on the semi-Krein space
$\mathbf{K}_{t}$. In the rest of this section, we study operator-theoretic
properties of hypercomplex numbers.

Let $B\left(\mathbf{X}_{t}\right)$ be an operator space consisting
of all bounded operators on $\mathbf{X}_{t}$ ``over $\mathbb{R}$,''
i.e., 
\[
B\left(\mathbf{X}_{t}\right)\overset{\textrm{def}}{=}\left\{ T:\mathbf{X}_{t}\rightarrow\mathbf{X}_{t}\mid T\textrm{ is bounded linear on }\mathbf{X}_{t}\right\} ,
\]
with its operator-norm $\left\Vert .\right\Vert $,\hfill{}(5.2)
\[
\left\Vert T\right\Vert =\mathrm{sup}\left\{ \left\Vert T\left(h\right)\right\Vert _{t}:\left\Vert h\right\Vert _{t}=1\right\} ,\;\forall T\in B\left(\mathbf{X}_{t}\right).
\]

\noindent Then, by the above proposition, the family

\medskip{}

\hfill{}$\mathcal{M}_{t}\overset{\textrm{def}}{=}\left\{ M_{\eta}:\eta\in\mathbb{H}_{t}\right\} ,$\hfill{}(5.3)

\medskip{}

\noindent consisting of the multiplication operators of (5.1) are
properly contained in the operator space $B\left(\mathbf{X}_{t}\right)$
of (5.2). It is not hard to check that 
\[
M_{h_{1}}+M_{h_{2}}=M_{h_{1}+h_{2}},\;\mathrm{and\;}M_{h_{1}}M_{h_{2}}=M_{h_{1}\cdot_{t}h_{2}},
\]
in $\mathcal{M}_{t}$, for all $h_{1},h_{2}\in\mathbb{H}_{t}$, and
\[
rM_{h}=M_{rh}=M_{\left(r,0\right)\cdot_{t}h},\;\;\mathrm{in\;\;}\mathcal{M}_{t},
\]
for all $r\in\mathbb{R}$, and $h\in\mathbb{H}_{t}$ (See the proof
of the following theorem). Thus, on this family $\mathcal{M}_{t}$
of (5.3), one can define a unary operation ($*$) by

\medskip{}

\hfill{}$*:M_{\eta}\in\mathcal{M}_{t}\rightarrow M_{\eta}^{*}\overset{\textrm{def}}{=}M_{\eta^{\dagger}}\in\mathcal{M}_{t},$\hfill{}(5.4)

\medskip{}

\noindent for all $\eta\in\mathbb{H}_{t}$. Then this operation (5.4)
is a well-defined adjoint on $\mathcal{M}_{t}$, because (i) $\mathcal{M}_{t}$
is equipotent (or, bijective) to $\mathbb{H}_{t}$, set-theoretically,
and (ii) the hypercomplex-conjugate ($\dagger$) is a well-defined
adjoint on $\mathbb{H}_{t}$ by (3.3). 
\begin{theorem}
The family $\mathcal{M}_{t}$ of (5.3) is an embedded complete semi-normed
$*$-algebra over $\mathbb{R}$ in the operator space $B\left(\mathbf{X}_{t}\right)$
of (5.2). 
\end{theorem}

\begin{proof}
Clearly, by definitions, one has the set-inclusion, $\mathcal{M}_{t}\subset B\left(\mathbf{X}_{t}\right)$.
If $M_{h_{1}},\:M_{h_{2}}\in\mathcal{M}_{t}$, then, for any $h\in\mathbf{X}_{t}$,

\medskip{}

$\;\;\;\;$$\left(M_{h_{1}}+M_{h_{2}}\right)\left(h\right)=M_{h_{1}}\left(h\right)+M_{h_{2}}\left(h\right)$

\medskip{}

$\;\;\;\;\;\;\;\;\;\;\;\;\;\;\;\;\;\;$$=\left(h_{1}\cdot_{t}h\right)+\left(h_{2}\cdot_{t}h\right)=\left(h_{1}+h_{2}\right)\cdot_{t}h$

\medskip{}

\noindent since $\mathbb{H}_{t}=\mathbf{X}_{t}$ is a ring satisfying
the distributedness of ($+$) and ($\cdot_{t}$)

\medskip{}

$\;\;\;\;\;\;\;\;\;\;\;\;\;\;\;\;\;\;$$=M_{h_{1}+h_{2}}\left(h\right)$,

\medskip{}

\noindent implying that\hfill{}(5.5) 
\[
M_{h_{1}}+M_{h_{2}}=M_{h_{1}+h_{2}}\in\mathcal{M}_{t},\;\forall h_{1},h_{2}\in\mathbb{H}_{t}.
\]
Also, we have 
\[
M_{h_{1}}M_{h_{2}}\left(h\right)=M_{h_{1}}\left(h_{2}\cdot_{t}h\right)=h_{1}\cdot_{t}\left(h_{2}\cdot_{t}h\right)=\left(h_{1}\cdot_{t}h_{2}\right)\cdot_{t}h,
\]
implying that\hfill{}(5.6) 
\[
M_{h_{1}}M_{h_{2}}=M_{h_{1}\cdot_{t}h_{2}}\in\mathcal{M}_{t},\;\;\forall h_{1},h_{2}\in\mathbb{H}_{t},
\]
where $M_{h_{1}}M_{h_{2}}$ means the operator-multiplication (or,
the operator-composition) on $B\left(\mathbf{X}_{t}\right)$. Therefore,
for any $h_{1},h_{2},h_{3}\in\mathcal{M}_{t}$, we have 
\[
M_{h_{1}}\left(M_{h_{2}}+M_{h_{3}}\right)=M_{h_{1}\cdot_{t}\left(h_{2}+h_{3}\right)}=M_{h_{1}}M_{h_{2}}+M_{h_{1}}M_{h_{3}},
\]
and similarly\hfill{}(5.7) 
\[
\left(M_{h_{1}}+M_{h_{2}}\right)M_{h_{3}}=M_{h_{1}}M_{h_{3}}+M_{h_{2}}M_{h_{3}},
\]
in $\mathcal{M}_{t}$, by (5.5) and (5.6). Also, by (5.6), for every
$r\in\mathbb{R}$ and $h\in\mathbf{X}_{t}$, we have 
\[
rM_{\eta}\left(h\right)=r\left(\eta\cdot_{t}h\right)=\left(r\eta\right)\cdot_{t}h=\left(\left(r,0\right)\cdot_{t}\eta\right)\cdot_{t}h,
\]
implying that\hfill{}(5.8) 
\[
rM_{\eta}=M_{r\eta}=M_{(r,0)\cdot_{t}\eta}\in\mathcal{M}_{t}.
\]
These show that the family $\mathcal{M}_{t}$ is a vector space over
$\mathbb{R}$ by (5.5) and (5.8), and this $\mathbb{R}$-vector space
$\mathcal{M}_{t}$ forms an algebra over $\mathbb{R}$, by (5.6) and
(5.7). Moreover, as we have discussed in the very above paragraph,
this $\mathbb{R}$-algebra $\mathcal{M}_{t}$ is equipped with a well-defined
adjoint ($*$) of (5.4), defined by 
\[
M_{\eta}^{*}=M_{\eta^{\dagger}}\;\;\mathrm{on\;\;}\mathcal{M}_{t},\;\;\forall\eta\in\mathbb{H}_{t}.
\]
Therefore, $\mathcal{M}_{t}$ is a well-determined $*$-algebra over
$\mathbb{R}$, embedded in $B\left(\mathbf{X}_{t}\right)$.

Now, observe that, for any $M_{\eta}\in\mathcal{M}_{t}$ with $\eta\in\mathbb{H}_{t}$,

\medskip{}

$\;\;\;\;$$\left\Vert M_{\eta}\right\Vert =\mathrm{sup}\left\{ \left\Vert M_{\eta}\left(h\right)\right\Vert _{t}:\left\Vert h\right\Vert _{t}=1\right\} $

\medskip{}

$\;\;\;\;\;\;\;\;$$=\mathrm{sup}\left\{ \left\Vert \eta\cdot_{t}h\right\Vert _{t}:\left\Vert h\right\Vert _{t}=1\right\} \leq\mathrm{sup}\left\{ \left\Vert \eta\right\Vert _{t}\left\Vert h\right\Vert _{t}:\left\Vert h\right\Vert _{t}=1\right\} $

\medskip{}

\noindent implying that\hfill{}(5.8) 
\[
\left\Vert M_{\eta}\right\Vert =\left\Vert \eta\right\Vert _{t}.
\]
So, this $*$-algebra $\mathcal{M}_{t}$ is complete over $\mathbb{R}$
by (5.8), since $\mathbf{X}_{t}\in\left\{ \mathbf{H}_{t},\mathbf{K}_{t}\right\} $
is complete under $\left\Vert .\right\Vert _{t}$. 
\end{proof}
The above theorem characterizes the family $\mathcal{M}_{t}$ as a
complete semi-normed $*$-algebra over $\mathbb{R}$ embedded in the
operator space $B\left(\mathbf{X}_{t}\right)$, where $\mathbf{X}_{t}\in\left\{ \mathbf{H}_{t},\mathbf{K}_{t}\right\} $,
for $t\in\mathbb{R}$. It means that (i) our $t$-scaled hypercomplex
ring $\mathbb{H}_{t}$ can be understood as a complete semi-normed
$*$-algebra over $\mathbb{R}$, acting on the $t$-hypercomplex space
$\mathbf{X}_{t}$, and (ii) all hypercomplex numbers $\mathit{h}\in\mathbb{H}_{t}$
are regarded as adjointable operators $M_{h}\in\mathcal{M}_{t}$ in
the operator space $B\left(\mathbf{X}_{t}\right)$. 
\begin{definition}
We call the complete semi-normed $*$-algebra $\mathcal{M}_{t}$ over
$\mathbb{R}$, the hypercomplex ($*$-)algebra over $\mathbb{R}$
in the operator space $B\left(\mathbf{X}_{t}\right)$, where $\mathbf{X}_{t}\in\left\{ \mathbf{H}_{t},\mathbf{K}_{t}\right\} $
for $t\in\mathbb{R}$. 
\end{definition}

By regarding each hypercomplex number $h\in\mathbb{H}_{t}$ as an
element $M_{h}$ of the hypercomplex algebra $\mathcal{M}_{t}$ over
$\mathbb{R}$, acting on the $t$-hypercomplex $\mathbb{R}$-space
$\mathbf{X}_{t}$, we study some operator-theoretic properties on
$\mathcal{M}_{t}$ under the following generalized concepts from the
usual $C^{*}$-algebra theory. 
\begin{definition}
Let $\mathscr{A}$ be a complete semi-normed $*$-algebra over $\mathbb{R}$,
and $A\in\mathscr{A}$.

\noindent (1) $A$ is self-adjoint in $\mathscr{A}$, if $A^{*}=A$
in $\mathscr{A}$.

\noindent (2) $A$ is a projection in $\mathscr{A}$, if $A^{*}=A=A^{2}$
in $\mathscr{A}$.

\noindent (3) $A$ is normal in $\mathscr{A}$, if $A^{*}A=AA^{*}$
in $\mathscr{A}$.

\noindent (4) $A$ is an isometry in $\mathscr{A}$, if $A^{*}A=1_{\mathscr{A}}$,
the identity element of $\mathscr{A}$.

\noindent (5) $A$ is a unitary in $\mathscr{A}$, if $A^{*}A=1_{\mathscr{A}}=AA^{*}$
in $\mathscr{A}$. 
\end{definition}

The following result characterizes the self-adjointness on the hypercomplex
$\mathbb{R}$-algebra $\mathcal{M}_{t}$. 
\begin{theorem}
An element $M_{\left(a,b\right)}$ is self-adjoint in $\mathcal{M}_{t}$,
if and only if 
\[
\left(a,b\right)=\left(Re\left(a\right),0\right)\;\;\mathrm{in\;\;}\mathbb{H}_{t},
\]
$\Longleftrightarrow$\hfill{}(5.9) 
\[
a\in\mathbb{R},\;\mathrm{and\;}b=0,\;\;\mathrm{in\;\;}\mathbb{C}.
\]
\end{theorem}

\begin{proof}
Let $h=\left(a,b\right)\in\mathbb{H}_{t}$, and $M_{h}\in\mathcal{M}_{t}$.
Then $M_{h}$ is self-adjoint in $\mathcal{M}_{t}$, if and only if
\[
M_{h}^{*}=M_{h^{\dagger}}=M_{h}\Longleftrightarrow h^{\dagger}=h\;\;\mathrm{in\;\;}\mathbb{H}_{t},
\]
if and only if 
\[
h^{\dagger}=\left(\overline{a},\:-b\right)=\left(a,b\right)=h\;\;\mathrm{in\;\;}\mathbb{H}_{t},
\]
if and only if 
\[
\overline{a}=a\;\;\mathrm{and\;\;}-b=b,\;\;\mathrm{in}\;\;\mathbb{C},
\]
if and only if 
\[
a\in\mathbb{R},\;\;\mathrm{and\;\;}b=0,\;\;\mathrm{in\;\;}\mathbb{C},
\]
if and only if 
\[
h=\left(a,b\right)=\left(Re\left(a\right),0\right)\;\;\mathrm{in\;\;}\mathbb{H}_{t}.
\]
Therefore, the self-adjointness (5.9) holds. 
\end{proof}
The above self-adjointness (5.9) on the hypercomplex $\mathbb{R}$-algebra
$\mathcal{M}_{t}$ actually shows that a hypercomplex number $\left(a,b\right)$
is self-conjugate, if and only if 
\[
\left(a,b\right)=\left(Re\left(a\right),0\right)\;\;\mathrm{in\;\;}\mathbb{H}_{t}.
\]

\begin{theorem}
An operator $M_{\left(a,b\right)}$ is a projection in the hypercomplex
$\mathbb{R}$-algebra $\mathcal{M}_{t}$, if and only if

\medskip{}

\hfill{}$\mathrm{either\;}\left(a,b\right)=\left(1,0\right),\;\mathrm{or\;}\left(a,b\right)=\left(0,0\right),\;\mathrm{in\;}\mathbb{H}_{t}.$\hfill{}(5.10) 
\end{theorem}

\begin{proof}
By definition, an operator $M_{\left(a,b\right)}\in\mathcal{M}_{t}$
is a projection, if and only if it is both self-adjoint and idempotent
in the sense that: $M_{\left(a,b\right)}^{2}=M_{\left(a,b\right)}$
in $\mathcal{M}_{t}$, by the definition of projections in complete
semi-normed $*$-algebras. So, it is a projection, if and only if
\[
\left(a,b\right)=\left(Re\left(a\right),0\right)\overset{\textrm{say}}{=}\left(x,0\right),\;\mathrm{in}\;\mathbb{H}_{t},
\]
with $x\in\mathbb{R}$, by (5.9), and 
\[
\left(x,0\right)\cdot_{t}\left(x,0\right)=\left(x,0\right),\;\mathrm{in\;}\mathbb{H}_{t},
\]
if and only if 
\[
\left(x^{2},0\right)=\left(x,0\right)\;\;\mathrm{in\;\;}\mathbb{H}_{t},
\]
if and only if 
\[
x^{2}=x\;\;\mathrm{in\;\;}\mathbb{R},\Longleftrightarrow x=0,\;\mathrm{or\;}1.
\]
Therefore, the projection-characterization (5.10) holds. 
\end{proof}
The above theorem shows that the only projections in the hypercomplex
$\mathbb{R}$-algebra $\mathcal{M}_{t}$ are the identity operator
$1_{t}=M_{\left(1,0\right)}$ and the zero operator $0_{t}=M_{\left(0,0\right)}$,
by (5.10). 
\begin{theorem}
Every operator $M_{h}$ is normal in $\mathcal{M}_{t}$ for all $h\in\mathbb{H}_{t}$. 
\end{theorem}

\begin{proof}
Observe that 
\[
M_{\left(a,b\right)}^{*}M_{\left(a,b\right)}=M_{\left(a,b\right)^{\dagger}}M_{\left(a,b\right)}=M_{\left(a,b\right)^{\dagger}\cdot_{t}\left(a,b\right)},
\]
and 
\[
M_{\left(a,b\right)}M_{\left(a,b\right)}^{*}=M_{\left(a,b\right)}M_{\left(a,b\right)^{\dagger}}=M_{\left(a,b\right)\cdot_{t}\left(a,b\right)^{\dagger}},
\]
in $\mathcal{M}_{t}$, and

\medskip{}

$\;\;$$\left(a,b\right)^{\dagger}\cdot_{t}\left(a,b\right)=\left(\overline{a},-b\right)\cdot_{t}\left(a,b\right)=\pi_{t}^{-1}\left(\left(\begin{array}{cc}
\overline{a} & t\left(-b\right)\\
-\overline{b} & a
\end{array}\right)\left(\begin{array}{cc}
a & tb\\
\overline{b} & \overline{a}
\end{array}\right)\right)$

\medskip{}

$\;\;\;\;\;\;\;\;\;\;\;\;$$=\pi_{t}^{-1}\left(\left(\begin{array}{ccc}
\left|a\right|^{2}-t\left|b\right|^{2} &  & 0\\
\\
0 &  & \left|a\right|^{2}-t\left|b\right|^{2}
\end{array}\right)\right)$

\medskip{}

$\;\;\;\;\;\;\;\;\;\;\;\;$$=\left(\left|a\right|^{2}-t\left|b\right|^{2},\;0\right)=\pi_{t}^{-1}\left(\left(\begin{array}{cc}
a & tb\\
\overline{b} & \overline{a}
\end{array}\right)\left(\begin{array}{cc}
\overline{a} & t\left(-b\right)\\
-\overline{b} & a
\end{array}\right)\right)$

\medskip{}

$\;\;\;\;\;\;\;\;\;\;\;\;$$=\left(a,b\right)\cdot_{t}\left(\overline{a},-b\right)=\left(a,b\right)\cdot_{t}\left(a,b\right)^{\dagger},$

\medskip{}

\noindent in $\mathbb{H}_{t}$. Therefore,

\medskip{}

\hfill{}$M_{\left(a,b\right)}^{*}M_{\left(a,b\right)}=M_{\left(a,b\right)}M_{\left(a,b\right)}^{*}\;\;\mathrm{in\;\;}\mathcal{M}_{t},$\hfill{}(5.11)

\medskip{}

\noindent guaranteeing the normality of $M_{\left(a,b\right)}$ in
$\mathcal{M}_{t}$. 
\end{proof}
The above theorem shows that all operators of the hypercomplex $\mathbb{R}$-algebra
$\mathcal{M}_{t}$ are normal in $\mathcal{M}_{t}$, by (5.11). 
\begin{corollary}
An operator $M_{h}$ is an isometry in the hypercomplex $\mathbb{R}$-algebra
$\mathcal{M}_{t}$, if and only if it is unitary in $\mathcal{M}_{t}$,
i.e.,

\medskip{}

\hfill{}Isometry Property $=$ Unitarity on $\mathcal{M}_{t}$.\hfill{}(5.12) 
\end{corollary}

\begin{proof}
By the above theorem, all operators of $\mathcal{M}_{t}$ are automatically
normal. In general, every unitary $U$ is an isometry in a complete
semi-normed $*$-algebra $\mathscr{A}$, by definition, i.e., $U$
is unitary in $\mathscr{A}$, if and only if $U^{*}U=1_{\mathscr{A}}=UU^{*}$,
and hence, it satisfies the isometry property, $U^{*}U=1_{\mathscr{A}}$,
in $\mathscr{A}$. Of course, the converse does not hold, in general.

However, in the hypercomplex $\mathbb{R}$-algebraic case, since all
operators of $\mathcal{M}_{t}$ are normal, the isometry property,
\[
M_{h}^{*}M_{h}=1_{t},
\]
automatically includes 
\[
M_{h}^{*}M_{h}=1_{t}=M_{h}M_{h}^{*}\;\;\mathrm{in\;\;}\mathcal{M}_{t},
\]
i.e., an isometry becomes a unitary in $\mathcal{M}_{t}$. So, the
equivalence (5.12) holds true. 
\end{proof}
The above corollary shows that the isometry property and the unitarity
are equivalent on our hypercomplex $\mathbb{R}$-algebras $\left\{ \mathcal{M}_{t}\right\} _{t\in\mathbb{R}}$. 
\begin{theorem}
An operator $M_{h}$ is a unitary in $\mathcal{M}_{t}$, if and only
if

\medskip{}

\hfill{}$\left\Vert h\right\Vert _{t}^{2}=\left\langle h,h\right\rangle _{t}=\tau\left(h\cdot_{t}h^{\dagger}\right)=det\left(\left[h\right]_{t}\right)=1.$\hfill{}(5.13) 
\end{theorem}

\begin{proof}
An operator $M_{\left(a,b\right)}$ is a unitary in $\mathcal{M}_{t}$,
if and only if 
\[
M_{\left(a,b\right)}^{*}M_{\left(a,b\right)}=M_{\left(\left|a\right|^{2}-t\left|b\right|^{2},0\right)}=M_{\left(1,0\right)}=1_{t},
\]
in $\mathcal{M}_{t}$, by (5.11) and (5.12), if and only if 
\[
\left|a\right|^{2}-t\left|b\right|^{2}=1=det\left(\left[\left(a,b\right)\right]_{t}\right).
\]
Therefore, the unitarity (5.13) holds by (5.12). 
\end{proof}
The above theorem characterizes both the unitarity and the isometry
property on hypercomplex $\mathbb{R}$-algebras by (5.13).

\medskip{}

\noindent $\mathbf{Observation.}$ It is interesting that the operator-theoretic
properties, the self-adjointness (5.9), the projection property (5.10),
the normality (5.11), and the unitarity (and the isometry property)
(5.13) are determined free from the choices of scales $t\in\mathbb{R}$.
i.e., these operator-theoretic properties are characterized universally
on $\left\{ \mathcal{M}_{t}\right\} _{t\in\mathbb{R}}$, even though
the $t$-hypercomplex space $\mathbf{X}_{t}$ can be either 
\[
\mathbf{X}_{t}=\mathbf{H}_{t},\;\textrm{the semi-Hilbert space, if }t<0,
\]
or 
\[
\mathbf{X}_{t}=\mathbf{K}_{t},\;\textrm{the semi-Krein space, if }t\geq0,
\]
over $\mathbb{R}$. e.g., see (4.28).\hfill{}{$\square$}

\medskip{}

In the rest of this section, we consider a naturally determined free
probability on the hypercomplex $\mathbb{R}$-algebra $\mathcal{M}_{t}$.
Note again that this complete semi-normed $*$-algebra $\mathcal{M}_{t}$
is equipotent to the $t$-scaled hypercomplex ring $\mathbb{H}_{t}$,
which means $\mathbb{H}_{t}$ is $*$-isomrophic to $\mathcal{M}_{t}$.
So, our free probability on $\mathcal{M}_{t}$ is actually studying
free-probabilistic information on $\mathbb{H}_{t}$. For more about
basic free probability theory, see e.g., {[}18{]} and {[}22{]}. We
here focus on free-distributional data on $\mathcal{M}_{t}$ by computing
joint free moments of elements of $\mathcal{M}_{t}$ under a canonically
defined linear functional.

Define a linear functional on the hypercomplex $\mathbb{R}$-algebra,
also denoted by $\tau$, on $\mathcal{M}_{t}$ by a linear morphism
satisfying

\medskip{}

\hfill{}$\tau\left(M_{\left(a,b\right)}\right)=\tau\left(\left(a,b\right)\right)=Re\left(a\right),\;\;\forall M_{\left(a,b\right)}\in\mathcal{M}_{t}$,\hfill{}(5.14)

\medskip{}

\noindent where the linear functional $\tau$ on the right-hand side
of (5.14) is in the sense of (4.3). Then one can have a well-defined
noncommutative free probability space $\left(\mathcal{M}_{t},\tau\right)$
(e.g., {[}22{]}), and the free distribution of an element $T\in\mathcal{M}_{t}$
is characterized by the joint free moments of $\left\{ T,T^{*}\right\} $,
\[
\tau\left(\overset{n}{\underset{l=1}{\prod}}T^{r_{l}}\right)=\tau\left(T^{r_{1}}T^{r_{2}}...T^{r_{n}}\right),
\]
for all $\left(r_{1},...,r_{n}\right)\in\left\{ 1,*\right\} ^{n}$,
for all $n\in\mathbb{N}$. So, if $T$ is self-adjoint in $\mathcal{M}_{t}$,
then the free distribution of $T$ is characterized by the free-moment
sequence, 
\[
\left(\tau\left(T^{n}\right)\right)_{n=1}^{\infty}=\left(\tau\left(T\right),\tau\left(T^{2}\right),\tau\left(T^{3}\right),...\right).
\]

Observe that, for any arbitrary $M_{h}\in\mathcal{M}_{t}$, with $h=\left(a,b\right)\in\mathbb{H}^{t}$,

\medskip{}

\hfill{}$\tau\left(\overset{n}{\underset{l=1}{\prod}}M_{h}^{r_{l}}\right)=\tau\left(\overset{n}{\underset{l=1}{\prod}}M_{h^{e_{l}}}\right)=\tau\left(M_{\overset{n}{\underset{l=1}{\cdot_{t}}}h^{e_{l}}}\right)=\tau\left(\overset{n}{\underset{l=1}{\cdot_{t}}}h^{e_{l}}\right),$\hfill{}(5.15)

\medskip{}

\noindent where 
\[
e_{l}=\left\{ \begin{array}{ccc}
1 &  & \mathrm{if\;}r_{l}=1\\
\dagger &  & \mathrm{if\;}r_{l}=*,
\end{array}\right.
\]
for all $l=1,...,n$, for all $\left(r_{1},...,r_{n}\right)\in\left\{ 1,*\right\} ^{n}$,
for all $n\in\mathbb{N}$. As one can see, computing (5.15) seems
complicated as $n$ is increasing in $\mathbb{N}$. So, we restrict
our interests to a special case of (5.15). 
\begin{theorem}
Assume that a hypercomplex number $h$ is similar to $\left(\sigma_{t}\left(h\right),0\right)$
in $\mathbb{H}_{t}$, where $\sigma_{t}\left(h\right)\in\mathbb{C}$
is the $t$-spectral value of $h$, and suppose $\sigma_{t}\left(h\right)$
is polar-decomposed to be 
\[
\sigma_{t}\left(h\right)=re^{i\theta},\;\mathrm{with\;}r=\left|\sigma_{t}\left(h\right)\right|,\;\theta=Arg\left(\sigma_{t}\left(h\right)\right),
\]
in $\mathbb{C}$. Then 
\[
\tau\left(\overset{n}{\underset{l=1}{\prod}}M_{h}^{r_{l}}\right)=r^{n}cos\left(\overset{n}{\underset{l=1}{\sum}}\varepsilon_{l}\theta\right),
\]
with\hfill{}(5.16) 
\[
\varepsilon_{l}=\left\{ \begin{array}{ccc}
1 &  & \mathrm{if\;}r_{l}=1\\
-1 &  & \mathrm{if\;}r_{l}=*,
\end{array}\right.
\]

\noindent for $l=1,...,n$, for all $\left(r_{1},...,r_{n}\right)\in\left\{ 1,*\right\} ^{n}$,
for all $n\in\mathbb{N}$. 
\end{theorem}

\begin{proof}
By the assumption that $h$ and $\left(\sigma_{t}\left(h\right),0\right)$
are similar in $\mathbb{H}_{t}$, the elements, 
\[
\overset{n}{\underset{l=1}{\cdot_{t}}}h^{e_{l}}\;\mathrm{and\;}\overset{n}{\underset{l=1}{\cdot_{t}}}\left(\sigma_{t}\left(h\right),0\right)^{e_{l}},
\]
are similar in $\mathbb{H}_{t}$, implying that 
\[
\tau\left(\overset{n}{\underset{l=1}{\cdot_{t}}}h^{e_{l}}\right)=\tau\left(\overset{n}{\underset{l=1}{\cdot_{t}}}\left(\sigma_{t}\left(h\right),0\right)^{e_{l}}\right),
\]
for all $\left(e_{1},...,e_{n}\right)\in\left\{ 1,\dagger\right\} ^{n}$,
for all $\mathit{n\in\mathbb{N}}$, where $\tau$, here, is in the
sense of (4.3).

For any $\left(e_{1},...,e_{n}\right)\in\left\{ 1,\dagger\right\} ^{n}$,
for $n\in\mathbb{N}$, one has that

\medskip{}

$\;\;\;\;$$\tau\left(\overset{n}{\underset{l=1}{\cdot_{t}}}\left(\sigma_{t}\left(h\right),0\right)^{e_{l}}\right)=\tau\left(\left(\overset{n}{\underset{l=1}{\prod}}\sigma_{t}\left(h\right)^{e_{l}},\;0\right)\right)$

\medskip{}

\noindent where 
\[
\sigma_{t}\left(h\right)^{e_{l}}=\left\{ \begin{array}{ccc}
\sigma_{t}\left(h\right) &  & \mathrm{if\;}e_{l}=1\\
\\
\overline{\sigma_{t}\left(h\right)} &  & \mathrm{if\;}e_{l}=\dagger,
\end{array}\right.
\]
for all $l=1,...,n$, and hence, it goes to

\medskip{}

$\;\;\;\;\;\;\;\;$$=\tau\left(\left(\overset{n}{\underset{l=1}{\prod}}\left(re^{i\theta}\right)^{e_{l}},\;0\right)\right)=Re\left(\overset{n}{\underset{l=1}{\prod}}\left(re^{i\theta}\right)^{e_{l}}\right)$

\medskip{}

$\;\;\;\;\;\;\;\;$$=r^{n}Re\left(\overset{n}{\underset{l=1}{\prod}}\left(e^{i\theta}\right)^{e_{l}}\right)=r^{n}Re\left(cos\left(\overset{n}{\underset{l=1}{\sum}}\varepsilon_{l}\theta\right)+isin\left(\overset{n}{\underset{l=1}{\sum}}\varepsilon_{l}\theta\right)\right),$

\medskip{}

\noindent where 
\[
\varepsilon_{l}=\left\{ \begin{array}{ccc}
1 &  & \mathit{\mathrm{if\;}e_{l}=1}\\
-1 &  & \mathrm{if\;}e_{l}=\dagger,
\end{array}\right.
\]
for all $l=1,...,n$, and hence,

\medskip{}

$\;\;\;\;\;\;\;\;$$=r^{n}cos\left(\overset{n}{\underset{l=1}{\sum}}\varepsilon_{l}\theta\right)$.

\medskip{}

\noindent Therefore, the free-distributional data (5.16) holds. 
\end{proof}
The above theorem fully characterizes the free-distributional information
on $\left(\mathcal{M}_{t},\tau\right)$ by (5.16). Note that the similarity-assumption
of $h$ and $\left(\sigma_{t}\left(h\right),0\right)$ in $\mathbb{H}_{t}$
is crucial in the above proof of (5.16). 
\begin{corollary}
If $t<0$, then every $M_{h}\in\mathcal{M}_{t}$ has the joint free
moments (5.16). 
\end{corollary}

\begin{proof}
Recall that, if $t<0$, then $h$ and $\left(\sigma_{t}\left(h\right),0\right)$
are similar in $\mathbb{H}_{t}$, for ``all'' $h\in\mathbb{H}_{t}$.
So, for all $h\in\mathbb{H}_{t}$, the corresponding multiplication
operators $M_{h}\in\mathcal{M}_{t}$ have their joint free moments
(5.16). 
\end{proof}
Different from the above corollary, if $t\geq0$, then we still need
the similarity-assumption to obtain the free-distributional data (5.16).
Without the similarity-assumption like the above theorem, the free-probabilistic
information is determined by (5.15) for $t\geq0$, in general.

\section{Scaled Hyperbolic Numbers}

In this section, we study the so-called scaled hyperbolic numbers
and related mathematical structures based on the main results of Sections
2 through 5.

\subsection{Motivation: The Hyperbolic Numbers}

In this section, we briefly review the hyperbolic numbers (or, the
split-complex numbers). For $x,y\in\mathbb{R}$, a hyperbolic number
$w$ is defined by 
\[
w=x+yj,\;\;\mathrm{with\;\;}j^{2}=1,
\]
where $j$ is an imaginary number satisfying $j^{2}=1$. By $\mathcal{D}$,
we denote the set of all such hyperbolic numbers, i.e., 
\[
\mathcal{D}=\left\{ x+yj:j^{2}=1,\textrm{ and }x,y\in\mathbb{R}\right\} .
\]
It is well-known that this set $\mathcal{D}$ is isomorphic to the
$2$-dimensional vector space $\mathbb{R}^{2}$ over $\mathbb{R}$,
and it is realized to be 
\[
\mathscr{D}=\left\{ \left(\begin{array}{cc}
x & y\\
y & x
\end{array}\right)\in M_{2}\left(\mathbb{R}\right):x,y\in\mathbb{R}\right\} ,
\]
i.e., each $x+yj\in\mathcal{D}$, assigned to be $\left(x,y\right)\in\mathbb{R}^{2}$,
is realized to be $\left(\begin{array}{cc}
x & y\\
y & x
\end{array}\right)\in\mathscr{D}$ in $M_{2}\left(\mathbb{R}\right)$.

The hyperbolic addition ($+$) is defined to be 
\[
\left(x_{1}+y_{1}j\right)+\left(x_{2}+y_{2}j\right)=\left(x_{1}+x_{2}\right)+\left(y_{1}+y_{2}\right)j,
\]
as in the vector addition on $\mathbb{R}^{2}$, and the hyperbolic
multiplication ($\cdot$) is defined by 
\[
\left(x_{1}+y_{1}j\right)\left(x_{2}+y_{2}j\right)=\left(x_{1}x_{2}+y_{1}y_{2}\right)+\left(x_{1}y_{2}+x_{2}y_{1}\right)j,
\]
by the matrix multiplication on $\mathscr{D}$, for all $x_{l}+y_{l}j\in\mathcal{D}$,
for $l=1,2$. Also, the $\mathbb{R}$-scalar product is defined by
\[
\left(r,x+yj\right)\in\mathbb{R}\times\mathcal{D}\longmapsto\left(r+0j\right)\left(x+yj\right)\in\mathcal{D}.
\]
So, the hyperbolic numbers $\mathcal{D}$ forms an algebra over $\mathbb{R}$.

Just like the complex numbers $\mathbb{C}$, one can define the hyperbolic-conjugate
($\dagger$) on $\mathcal{D}$ by 
\[
\left(x+yj\right)^{\dagger}=x-yj=x+\left(-y\right)j,
\]
for all $x+yj\in\mathcal{D}$, which is a well-defined adjoint on
$\mathcal{D}$, making the $\mathbb{R}$-algebra $\mathcal{D}$ to
the $*$-algebra over $\mathbb{R}$, equipped with its adjoint ($\dagger$).
By this adjoint, one can have that 
\[
\left(x+yj\right)^{\dagger}\left(x+yj\right)=x^{2}-y^{2}=\left(x+yj\right)\left(x+yj\right)^{\dagger},
\]
in $\mathcal{D}$, for all $x+yj\in\mathcal{D}$. Actually, the above
equality provides a difference between the complex numbers $\mathbb{C}$
and the hyperbolic numbers $\mathcal{D}$, i.e., 
\[
w^{\dagger}w\in\mathbb{R},\;\;\forall w\in\mathcal{D},
\]
which means $w^{\dagger}w$ can be negative in $\mathbb{R}$. So,
one can define the indefinite semi-inner product, 
\[
\left\langle w_{1},w_{2}\right\rangle \overset{\textrm{def}}{=}w_{1}w_{2}^{\dagger},\;\;\forall w_{1},w_{2}\in\mathcal{D},
\]
on $\mathcal{D}$, inducing the semi-norm, 
\[
\left\Vert w\right\Vert =\sqrt{\left|\left\langle w,w\right\rangle \right|},\;\;\forall w\in\mathcal{D}.
\]
Since $\mathcal{D}$ is finite-dimensional over $\mathbb{R}$, the
pair $\left(\mathcal{D},\left\langle ,\right\rangle \right)$ forms
a semi-Krein space over $\mathbb{R}$ (in our sense).

Similar to our constructions of Sections 4 and 5, the hyperbolic numbers
$\mathcal{D}$ forms a complete semi-normed $*$-algebra over $\mathbb{R}$,
acting on the semi-Krein space $\left(\mathcal{D},\left\langle ,\right\rangle \right)$. 
\begin{theorem}
Let $\mathcal{D}$ be the hyperbolic numbers, understood to be a complete
semi-normed $*$-algebra over $\mathbb{R}$ acting on the semi-Krein
space $\left(\mathcal{D},\left\langle ,\right\rangle \right)$. Then
it is a $*$-subalgebra of the hypercomplex algebra $\mathcal{M}_{1}$. 
\end{theorem}

\begin{proof}
First of all, note that the $1$-scaled hypercomplex ring $\mathbb{H}_{1}$
(which is the set of all split-quaternions of {[}1{]}) forms the semi-Krein
space $\mathbf{K}_{1}=\left(\mathbb{H}_{1},\left\langle ,\right\rangle _{1}\right)$,
and the corresponding hypercomplex $\mathbb{R}$-algebra $\mathcal{M}_{1}$,
which is equipotent to $\mathbb{H}_{1}$, acts on $\mathbf{K}_{1}$
as multiplication operators. Define a subset $\mathbb{D}_{1}$ of
$\mathbb{H}_{1}$ by 
\[
\mathbb{D}_{1}\overset{\textrm{def}}{=}\left\{ \left(x,y\right)\in\mathbb{H}_{1}:x,y\in\mathbb{R}\right\} ,
\]
realized to be 
\[
\pi_{1}\left(\mathbb{D}_{1}\right)=\left\{ \left(\begin{array}{cc}
x & y\\
y & x
\end{array}\right)\in\mathcal{H}_{2}^{1}:x,y\in\mathbb{R}\right\} \;\mathrm{in\;}\mathcal{H}_{2}^{1}.
\]
Then, by the very construction, it is not difficult to check that
(i) the semi-Krein space $\left(\mathcal{D},\left\langle ,\right\rangle \right)$
is isomorphic to the closed subspace $\left(\mathbb{D}_{1},\left\langle ,\right\rangle _{1}\right)$
of the semi-Krein space $\mathbf{K}_{1}$ by the isomorphism, 
\[
x+yj\in\left(\mathcal{D},\left\langle ,\right\rangle \right)\longmapsto\left(x,y\right)\in\left(\mathbb{D}_{1},\left\langle ,\right\rangle _{1}\right),
\]
and (ii) the subset, 
\[
\mathcal{D}_{1}=\left\{ M_{h}\in\mathcal{M}_{1}:h\in\mathbb{D}_{1}\;\mathrm{in\;}\mathbb{H}_{1}\right\} ,
\]
forms a $*$-subalgebra of the hypercomplex algebra $\mathcal{M}_{1}$
over $\mathbb{R}$, acting on $\left(\mathbb{D}_{1},\left\langle ,\right\rangle _{1}\right)$.
So, one can define a function, 
\[
\Phi:\mathcal{D}\rightarrow\mathcal{D}_{1},
\]
by\hfill{}(6.1.1) 
\[
\Phi\left(x+yj\right)=M_{\left(x,y\right)},\;\;\forall x+yj\in\mathcal{D}.
\]
Then this bijection $\Phi$ of (6.1.1) is a $*$-isomorphism over
$\mathbb{R}$. Therefore, the hyperbolic numbers $\mathcal{D}$ is
$*$-isomorphic to the $*$-subalgebra $\mathcal{D}_{1}$ of our hypercomplex
algebra $\mathcal{M}_{1}$. 
\end{proof}
The above theorem shows that the hyperbolic numbers $\mathcal{D}$
is understood to be a sub-structure of the $1$-scaled hypercomplex
ring $\mathbb{H}_{1}$; the semi-Krein space $\left(\mathcal{D},\left\langle ,\right\rangle \right)$
is a subspace of the semi-Krein space $\mathbf{K}_{1}$; and $\mathcal{D}$
is $*$-isomorphic to the $*$-subalgebra of the hypercomplex algebra
$\mathcal{M}_{1}$, as a complete semi-normed $*$-algebra over $\mathbb{R}$.

Motivated by the above theorem, we generalize the hyperbolic numbers
to the scaled-hyperbolic-numbers.

\subsection{Scaled Hyperbolic Numbers}

Let $t\in\mathbb{R}$, and $\mathbb{H}_{t}$, the $t$-scaled hypercomplex
ring, inducing the $t$-hypercomplex $\mathbb{R}$-space, 
\[
\mathbf{X}_{t}=\left(\mathbb{H}_{t},\left\langle ,\right\rangle _{t}\right)\in\left\{ \mathbf{H}_{t},\mathbf{K}_{t}\right\} ,
\]
and the hypercomplex $\mathbb{R}$-algebra $\mathcal{M}_{t}$. Define
a subset $\mathbb{D}_{t}$ of $\mathbb{H}_{t}$ by

\medskip{}

\hfill{}$\mathbb{D}_{t}\overset{\textrm{def}}{=}\left\{ \left(x,y\right)\in\mathbb{H}_{t}:x,y\in\mathbb{R}\right\} ,$\hfill{}(6.2.1)

\medskip{}

\noindent realized to be

\medskip{}

\hfill{}$\mathcal{D}_{2}^{t}\overset{\textrm{def}}{=}\pi_{t}\left(\mathbb{D}_{t}\right)=\left\{ \left[\left(x,y\right)\right]_{t}=\left(\begin{array}{cc}
x & ty\\
y & x
\end{array}\right):\left(x,y\right)\in\mathbb{D}_{t}\right\} ,$\hfill{}(6.2.2)

\medskip{}

\noindent in $\mathcal{H}_{2}^{t}=\pi_{t}\left(\mathbb{H}_{t}\right)$.
Then it is not hard to show that ($+$) and ($\cdot_{t}$) are closed
on $\mathbb{D}_{t}$, and hence, $\mathbb{D}_{t}$ is a subring of
$\mathbb{H}_{t}$. 
\begin{definition}
The subring $\mathbb{D}_{t}$ of (6.2.1) is called the $t$-scaled
hyperbolic (sub)ring of the $t$-scaled hypercomplex ring $\mathbb{H}_{t}$.
And the realization $\mathcal{D}_{2}^{t}$ of (6.2.2) is said to be
the $t$-(scaled-)hyperbolic realization of $\mathbb{D}_{t}$, as
a subring of the $t$-scaled realization $\mathcal{H}_{2}^{t}$. 
\end{definition}

Recall that our $t$-scaled hypercomplex ring $\mathbb{H}_{t}$ forms
the hypercomplex $\mathbb{R}$-algebra, 
\[
\mathcal{M}_{t}=\left\{ M_{h}\in B\left(\mathbf{X}_{t}\right):h\in\mathbb{H}_{t}\right\} ,
\]
over $\mathbb{R}$, equipped with its adjoint, 
\[
M_{h}^{*}=M_{h^{\dagger}}\;\;\mathrm{in\;\;}\mathcal{M}_{t},\;\;\forall h\in\mathbb{H}_{t},
\]
acting on the $t$-hypercomplex space $\mathbf{X}_{t}=\left(\mathbb{H}_{t},\left\langle ,\right\rangle _{t}\right)$,
which is either the semi-Hilbert space $\mathbf{H}_{t}$ (if $t<0$),
or the semi-Krein space $\mathbf{K}_{t}$ (if $t\geq0$). Therefore,
as a sub-structure of $\mathbb{H}_{t}$, the $t$-scaled hyperbolic
ring $\mathbb{D}_{t}$ forms the closed $*$-subalgebra of $\mathcal{M}_{t}$,
as a complete semi-normed $*$-algebra over $\mathbb{R}$ acting on
the closed subspace $\left(\mathbb{D}_{t},\left\langle ,\right\rangle _{t}\right)$
of $\mathbf{X}_{t}$. 
\begin{theorem}
Let $\mathbb{D}_{t}$ be the $t$-scaled hyperbolic ring of $\mathbb{H}_{t}$.
Then there exists a closed subspace of the $t$-hypercomplex space
$\mathbf{X}_{t}$,

\medskip{}

\hfill{}$\mathbf{Y}_{t}\overset{\textrm{def}}{=}\left(\mathbb{D}_{t},\:\left\langle ,\right\rangle _{t}\right)\;\;\mathrm{in\;\;}\mathbf{X}_{t}=\left(\mathbb{H}_{t},\left\langle ,\right\rangle _{t}\right),$\hfill{}(6.2.3)

\medskip{}

\noindent and a closed $*$-subalgebra,

\medskip{}

\hfill{}$\mathscr{D}_{t}\overset{\textrm{def}}{=}\left\{ M_{\left(x,y\right)}\in\mathcal{M}_{t}:\left(x,y\right)\in\mathbb{D}_{t}\right\} ,$\hfill{}(6.2.4)

\medskip{}

\noindent of the hypercomplex $\mathbb{R}$-algebra $\mathcal{M}_{t}$
over $\mathbb{R}$, acting on $\mathbf{Y}_{t}$. 
\end{theorem}

\begin{proof}
The proof is done by the very constructions (6.2.1) and (6.2.2). Indeed,
the semi-normed space $\mathbf{Y}_{t}$ of (6.2.3) is a closed subspace
of the $t$-hypercomplex space $\mathbf{X}_{t}$, by (6.2.1) and (6.2.2).
And hence, the $*$-algebra $\mathscr{D}_{t}$ of (6.2.4) forms a
complete semi-normed $*$-algebra over $\mathbb{R}$, acting on $\mathbf{Y}_{t}$,
as a closed $*$-subalgebra of $\mathcal{M}_{t}$. 
\end{proof}
The above theorem characterizes the functional-analytic, and operator-algebraic
properties of the $t$-hyperbolic ring $\mathbb{D}_{t}$, as a sub-structure
of the $t$-scaled hypercomplex ring $\mathbb{H}_{t}$. Note that
the $*$-subalgebra $\mathscr{D}_{t}$ of (6.2.4) has its adjoint,
\[
M_{\left(x,y\right)}^{*}=M_{\left(x,-y\right)}\;\;\mathrm{in\;\;}\mathscr{D}_{t},\;\;\forall\left(x,y\right)\in\mathbb{D}_{t},
\]
since 
\[
\left(x,y\right)^{\dagger}=\left(\overline{x},-y\right)=\left(x,-y\right)\;\mathrm{in\;}\mathbb{D}_{t},
\]
because $x,y\in\mathbb{R}$.

One immediately obtains the following result, as examples of our scaled
hyperbolic rings $\left\{ \mathbb{D}_{t}\right\} _{t\in\mathbb{R}}$. 
\begin{corollary}
The $\left(-1\right)$-scaled hyperbolic ring $\mathbb{D}_{-1}$ is
the complex field $\mathbb{C}$, and the $1$-scaled hyperbolic ring
$\mathbb{D}_{1}$ is the ring $\mathcal{D}$ of all hyperbolic numbers
of Section 6.1. 
\end{corollary}

\begin{proof}
By construction, the ring $\mathbb{D}_{-1}$ is isomorphic to the
complex field $\mathbb{C}$ in the noncommutative field $\mathbb{H}_{-1}$
of all quaternions. Also, the ring $\mathbb{D}_{1}$ is isomorphic
to the ring $\mathcal{D}$ of all hyperbolic numbers, by (6.2.1),
(6.2.2), (6.2.3) and (6.2.4). 
\end{proof}
The above corollary shows that, indeed, the hyperbolic numbers $\mathbb{D}_{1}$,
the $1$-scaled hyperbolic ring, induces the semi-Krein subspace,
\[
\mathbf{Y}_{1}=\left(\mathbb{D}_{t},\left\langle ,\right\rangle _{1}\right),
\]
by (6.2.3), equipped with the indefinite semi-inner product $\left\langle ,\right\rangle _{1}$,
inherited from that on the semi-Krein space $\mathbf{K}_{1}=\left(\mathbb{H}_{1},\left\langle ,\right\rangle _{1}\right)$;
and it induces the complete semi-normed $*$-algebra over $\mathbb{R}$,
\[
\mathscr{D}_{1}=\left\{ M_{\left(x,y\right)}\in\mathcal{M}_{1}:\left(x,y\right)\in\mathbb{D}_{1}\right\} ,
\]
by (6.2.4), acting on $\mathbf{Y}_{1}$, as multiplication operators,
satisfying 
\[
\left\Vert M_{\left(x,y\right)}\right\Vert =\mathrm{sup}\left\{ \left\Vert M_{\left(x,y\right)}v\right\Vert _{1}:\left\Vert v\right\Vert _{1}=1\right\} =\left\Vert \left(x,y\right)\right\Vert _{1},
\]
and hence, 
\[
\left\Vert M_{\left(x,y\right)}\right\Vert =\left\Vert \left(x,y\right)\right\Vert _{1}=\sqrt{\left|\left\langle \left(x,y\right),\left(x,y\right)\right\rangle _{1}\right|}=\sqrt{\left|x^{2}-y^{2}\right|},
\]
where 
\[
\left\langle \left(x,y\right),\left(x,y\right)\right\rangle _{1}=\left(x,y\right)\cdot_{1}\left(x,y\right)^{\dagger}=x^{2}-y^{2}\in\mathbb{R},
\]
for all $\left(x,y\right)\in\mathbb{D}_{1}$. So, the study of our
scaled hyperbolic rings $\left\{ \mathbb{D}_{t}\right\} _{t\in\mathbb{R}}$
generalizes the classical hyperbolic theory. 
\begin{definition}
We call the closed subspace $\mathbf{Y}_{t}$ of $\mathbf{X}_{t}$,
in the sense of (6.2.3), the $t$(-scaled-)hyperbolic (vector) space
(over $\mathbb{R}$). And the closed $*$-subalgebra $\mathscr{D}_{t}$
of the hypercomplex $\mathbb{R}$-algebra $\mathcal{M}_{t}$, in the
sense of (6.2.4), is said to be the ($t$-scaled) hyperbolic (complete
semi-normed $*$-)algebra (over $\mathbb{R}$). 
\end{definition}

By the above theorem, one obtains the following result. 
\begin{corollary}
An operator $M_{\left(x,y\right)}$ is self-adjoint in the $t$-hyperbolic
algebra $\mathscr{D}_{t}$, if and only if

\medskip{}

\hfill{}$\left(x,y\right)=\left(x,0\right)\;\;\;\mathrm{in\;\;\;}\mathbb{D}_{t}$.\hfill{}(6.2.5) 
\end{corollary}

\begin{proof}
An operator $M_{\left(a,b\right)}$ is self-adjoint in $\mathcal{M}_{t}$,
if and only if 
\[
\left(a,b\right)=\left(Re\left(a\right),0\right)\;\;\mathrm{in\;\;}\mathbb{H}_{t},
\]
by (5.9). So, by (6.2.4), an operator $M_{\left(x,y\right)}$ is self-adjoint
in $\mathscr{D}_{t}$, if and only if the condition (6.2.5) holds. 
\end{proof}
The above result characterizes the self-adjointness on the $t$-hyperbolic
algebra $\mathscr{D}_{t}$ by (6.2.5). 
\begin{corollary}
The only zero operator $0_{t}=M_{\left(0,0\right)}$ and the identity
operator $1_{t}=M_{\left(1,0\right)}$ are projections in $\mathscr{D}_{t}$. 
\end{corollary}

\begin{proof}
An operator $M_{\left(a,b\right)}$ is a projection in $\mathcal{M}_{t}$,
if and only if either 
\[
\left(a,b\right)=\left(0,0\right),\;\mathrm{or\;}\left(a,b\right)=\left(1,0\right),
\]
by (5.10). Note that $\left(0,0\right),\left(1,0\right)\in\mathbb{D}_{t}$,
too. So, as a closed $*$-subalgebra, the $t$-hyperbolic algebra
$\mathscr{D}_{t}$ has its projections $\left\{ 0_{t},1_{t}\right\} $. 
\end{proof}
The above corollary shows that the hyperbolic algebras $\left\{ \mathscr{D}_{t}\right\} _{t\in\mathbb{R}}$
and the hypercomplex algebras $\left\{ \mathcal{M}_{t}\right\} _{t\in\mathbb{R}}$
share the same projections $\left\{ 0_{t},1_{t}\right\} _{t\in\mathbb{R}}$. 
\begin{corollary}
Every operator of $\mathscr{D}_{t}$ is normal. 
\end{corollary}

\begin{proof}
By (5.11), all elements of the hypercomplex $\mathbb{R}$-algebra
$\mathcal{M}_{t}$ are automatically normal. Since $\mathscr{D}_{t}\subset\mathcal{M}_{t}$,
all operators of $\mathscr{D}_{t}$ are normal. 
\end{proof}
Indeed, one has 
\[
M_{h}^{*}M_{h}=M_{h^{\dagger}\cdot_{t}h}=M_{\left(x^{2}-ty^{2},0\right)}=M_{h\cdot_{t}h^{\dagger}}=M_{h}M_{h}^{*},
\]
in $\mathscr{D}_{t}$, for all $h=\left(x,y\right)\in\mathbb{D}_{t}$,
with $x,y\in\mathbb{R}$. 
\begin{corollary}
An operator $M_{\left(x,y\right)}$ is unitary in $\mathscr{D}_{t}$,
if and only if it is an isometry in $\mathscr{D}_{t}$, if and only
if

\medskip{}

\hfill{}$x^{2}-ty^{2}=1,\;\;\;\mathrm{in\;\;\;}\mathbb{R}.$\hfill{}(6.2.6) 
\end{corollary}

\begin{proof}
By the above corollary, every element of $\mathscr{D}_{t}$ is normal.
So, if $M_{\left(x,y\right)}$ is an isometry in $\mathscr{D}_{t}$,
then it is unitary by the normality. Therefore, $M_{\left(x,y\right)}$
is an isometry in $\mathscr{D}_{t}$, if and only if it is unitary
in $\mathscr{D}_{t}$ by (5.12), if and only if the condition (6.2.6)
holds in $\mathbb{R}$, by (5.13). 
\end{proof}
As a $*$-subalgebra, the $t$-hyperbolic algebra $\mathscr{D}_{t}$
has a canonical free-probabilistic structure, $\left(\mathscr{D}_{t},\tau\right)$,
where $\tau=\tau\mid_{\mathscr{D}_{t}}$ is the restriction of the
linear functional $\tau$ of (5.14). 
\begin{corollary}
Let $\left(x,y\right)\in\mathbb{D}_{t}$ be similar to $\left(\sigma_{t}\left(\left(x,y\right)\right),0\right)$
in $\mathbb{H}_{t}$, where 
\[
\sigma_{t}\left(\left(x,y\right)\right)=x+i\sqrt{-ty^{2}},\;\;\mathrm{in\;\;}\mathbb{C},
\]
where\hfill{}(6.2.7) 
\[
x+i\sqrt{-ty^{2}}=\left\{ \begin{array}{ccc}
x+i\sqrt{-ty^{2}}\in\mathbb{C} &  & \mathrm{if\;}t<0\\
x\in\mathbb{R} &  & \mathrm{if\;}t=0\\
x-\sqrt{ty^{2}}\in\mathbb{R} &  & \mathrm{if\;}t>0.
\end{array}\right.
\]
Then\hfill{}(6.2.8) 
\[
\tau\left(\overset{n}{\underset{l=1}{\prod}}M_{\left(x,y\right)}^{r_{l}}\right)=\left\{ \begin{array}{ccc}
x^{n} &  & \mathrm{if\;}t\leq0\\
\\
\left(x-\sqrt{ty^{2}}\right)^{n} &  & \mathrm{if\;}t>0,
\end{array}\right.
\]
for all $\left(r_{1},...,r_{n}\right)\in\left\{ 1,*\right\} ^{n}$,
for all $n\in\mathbb{N}$. 
\end{corollary}

\begin{proof}
By (5.16), if $\left(a,b\right)$ and $\left(\sigma_{t}\left(\left(a,b\right)\right),0\right)$
are similar in $\mathbb{H}_{t}$, and $\sigma_{t}\left(\left(a,b\right)\right)$
is polar-decomposed to be $re^{i\theta}$, with $r=\left|\sigma_{t}\left(h\right)\right|$,
and $\theta=Arg\left(\sigma_{t}\left(h\right)\right)$, in $\mathbb{C}$,
then 
\[
\tau\left(\overset{n}{\underset{l=1}{\prod}}M_{\left(a,b\right)}^{r_{l}}\right)=r^{n}cos\left(\overset{n}{\underset{l=1}{\sum}}\varepsilon_{l}\theta\right),
\]
for all $\left(r_{1},...,r_{n}\right)\in\left\{ 1,*\right\} ^{n}$,
for all $n\in\mathbb{N}$, where $\left\{ \varepsilon_{1},...,\varepsilon_{n}\right\} $
satisfies the condition of (5.16).

If a hyperbolic number $\left(x,y\right)\in\mathbb{D}_{t}$ is similar
to 
\[
\left(\sigma_{t}\left(\left(x,y\right)\right),\:0\right)=\left(x+i\sqrt{-ty^{2}},\:0\right),
\]
by (2.3.3), whose first entry $x+t\sqrt{-ty^{2}}$ is refined by (6.2.7)
case-by-case for a scale $t\in\mathbb{R}$. So, one has 
\[
\tau\left(\overset{n}{\underset{l=1}{\prod}}M_{\left(x,y\right)}^{r_{l}}\right)=\left\{ \begin{array}{ccc}
x^{n} &  & \mathrm{if\;}t<0\\
\\
x^{n} &  & \mathrm{if\;}t=0\\
\\
\left(x-\sqrt{ty^{2}}\right)^{n} &  & \mathrm{if\;}t>0,
\end{array}\right.
\]
by (6.2.7). Therefore, the free-distributional data (6.2.8) holds. 
\end{proof}
In the proof of (6.2.8), the similarity-assumption is crucial. By
the above corollary, if $t<0$, then all elements $M_{\left(x,y\right)}\in\mathscr{D}_{t}$
for all $t$-scaled hyperbolic numbers $\left(x,y\right)\in\mathbb{D}_{t}$
have their joint free moments determined by the first formula in (6.2.8),
because if $t<0$, then all $t$-scaled hyperbolic numbers $\left(x,y\right)$
are automatically similar to $\left(\sigma_{t}\left(\left(x,y\right)\right),\:0\right)$
in the $t$-scaled hyperbolic ring $\mathbb{D}_{t}$.

\section{Unit Elements of $\mathbb{D}_{t}$}

In this section, we consider unit elements of our scaled hyperbolic
rings $\left\{ \mathbb{D}_{t}\right\} _{t\in\mathbb{R}}$. By definition,
the $\left(-1\right)$-hyperbolic ring $\mathbb{D}_{-1}$ is isomorphic
to the complex field $\mathbb{C}$, by the isomorphism, 
\[
\left(x,y\right)\in\mathbb{D}_{-1}\longmapsto x+yi\in\mathbb{C},\;\mathrm{with\;}i=\sqrt{-1},
\]
since the $\left(-1\right)$-hypercomplex ring $\mathbb{H}_{-1}$
is isomorphic to the noncommutative field $\mathbb{H}$ of all quaternions.
Similarly, as we discussed in Section 6.1, the $1$-hyperbolic ring
$\mathbb{D}_{1}$ is isomorphic to the hyperbolic numbers $\mathcal{D}$,
as a sub-structure of the $1$-hypercomplex ring $\mathbb{H}_{1}$,
which is isomorphic to the ring of all split-quaternions.

On the complex field $\mathbb{C}$, one can construct the unit circle,
\[
\mathbb{T}=\left\{ z\in\mathbb{C}:\left|z\right|=1\right\} ,
\]
where $\left|.\right|$ is the modulus on $\mathbb{C}$. Then $z\in\mathbb{T}$,
if and only if 
\[
z=e^{i\theta}=cos\theta+isin\theta,
\]
where $\theta=Arg\left(z\right)$ is the argument of $z$. Equivalently,
$z\in\mathbb{T}$, if and only if 
\[
\left\langle z,z\right\rangle _{-1}=1\Longleftrightarrow\left\Vert z\right\Vert _{-1}=1=\left|z\right|,
\]
by regarding $\mathbb{C}$ as its isomorphic structure $\mathbb{D}_{-1}$,
in $\mathbb{H}_{-1}$. These unit elements of $\mathbb{T}$ play important
roles in the study of $\mathbb{C}$ by the polar decomposition on
$\mathbb{C}$.

Similarly, on the hyperbolic numbers $\mathcal{D}$, one can define
the unit hyperbola, 
\[
\mathcal{T}=\left\{ w\in\mathcal{D}:\left\Vert w\right\Vert _{1}=1\right\} ,
\]
where $\left\Vert .\right\Vert _{1}$ is the semi-norm on the $1$-scaled
hyperbolic ring $\mathbb{D}_{1}$, regarding it as its isomorphic
structure $\mathcal{D}$, i.e., 
\[
\left\Vert x+yj\right\Vert _{1}=\sqrt{\left|\left\langle x+yj,x+yj\right\rangle _{1}\right|}=\sqrt{\left|x^{2}-y^{2}\right|},
\]
for all $x+yj\in\mathcal{D}$, understood as $\left(x,y\right)\in\mathbb{D}_{1}$.
Then this unit hyperbola $\mathcal{T}$ of $\mathcal{D}$ can be re-defined
by 
\[
\mathcal{T}=\left\{ w\in\mathcal{D}:\left\langle w,w\right\rangle _{1}=\pm1\right\} ,
\]
since 
\[
\left\langle x+yj,x+yj\right\rangle _{1}=x^{2}-y^{2}\in\mathbb{R}.
\]
i.e., different from the complex-number case, the unit hyperbolic
numbers $w\in\mathcal{T}$ satisfies either 
\[
\left\langle w,w\right\rangle _{1}=1,\;\mathrm{or\;}\left\langle w,w\right\rangle _{1}=-1.
\]
It is well-known that $w\in\mathcal{T}$, if and only if 
\[
w=e^{j\theta}=cosh\theta+jsinh\theta,
\]
where $\theta=Arg\left(w\right)$ as in the complex case, and 
\[
coshx=\frac{e^{x}+e^{-x}}{2},\;\mathrm{and\;}sinhx=\frac{e^{x}-e^{-x}}{2},
\]
for all $x\in\mathbb{R}$. These unit elements also play key roles
in hyperbolic analysis (e.g., see {[}14{]} and {[}15{]}).

Just like the above complex, and hyperbolic cases, equivalently, like
the $\mathbb{D}_{-1}$ and $\mathbb{D}_{1}$ cases, one may consider
unit elements of $\left\{ \mathbb{D}_{t}\right\} _{t\in\mathbb{R}}$.
Let's fix an arbitrary scale $t\in\mathbb{R}$ and the corresponding
hyperbolic ring $\mathbb{D}_{t}$, realized to be $\mathcal{D}_{2}^{t}$.
Since we already know the cases where $t\in\left\{ \pm1\right\} $,
we restrict our interests to the cases where $t\in\mathbb{R}\setminus\left\{ \pm1\right\} $.

Observe that $\left(x,y\right)\in\mathbb{D}_{t}$ satisfies 
\[
\left\Vert \left(x,y\right)\right\Vert _{t}=\sqrt{\left|\left\langle \left(x,y\right),\left(x,y\right)\right\rangle _{t}\right|}=1,
\]
by (4.2.8), if and only if 
\[
\left\langle \left(x,y\right),\left(x,y\right)\right\rangle _{t}=1,\;\mathrm{whenever\;}t<0,
\]
and\hfill{}(7.1) 
\[
\left\langle \left(x,y\right),\left(x,y\right)\right\rangle _{t}\in\left\{ \pm1\right\} ,\;\mathrm{whenever\;}t\geq0.
\]

\begin{theorem}
If $t<0$ in $\mathbb{R}$, and $\left(x,y\right)\in\mathbb{D}_{t}$,
then 
\[
\left\Vert \left(x,y\right)\right\Vert _{t}=1\Longleftrightarrow y^{2}=\frac{x^{2}-1}{t}\Longleftrightarrow y^{2}=\frac{1-x^{2}}{\left|t\right|},
\]
and hence,\hfill{}(7.2) 
\[
\left\Vert \left(x,y\right)\right\Vert _{t}=1\Longleftrightarrow y=\pm\sqrt{\frac{1-x^{2}}{\left|t\right|}},\;\;\mathrm{in\;\;}\mathbb{R},
\]
where $\left|.\right|$ is the absolute value on $\mathbb{R}$. 
\end{theorem}

\begin{proof}
By (7.1), if $t<0$, then $\left\Vert \left(x,y\right)\right\Vert _{t}=1$,
if and only if 
\[
\left\langle \left(x,y\right),\left(x,y\right)\right\rangle _{t}=1,
\]
if and only if 
\[
\tau\left(\left(x,y\right)\cdot_{t}\left(x,y\right)^{\dagger}\right)=det\left(\left[\left(x,y\right)\right]_{t}\right)=x^{2}-ty^{2}=1,
\]
and hence, the unit characterization (7.2) holds on $\mathbb{D}_{t}$. 
\end{proof}
It is not hard to check that the characterization (7.2) actually contains
the complex case where $t=-1$. 
\begin{theorem}
If $t\geq0$, and $\left(x,y\right)\in\mathbb{D}_{t}$, then

\medskip{}

\hfill{}$\left\Vert \left(x,y\right)\right\Vert _{t}=1\Longleftrightarrow x^{2}-ty^{2}=\pm1.$\hfill{}(7.3)

\medskip{}

\noindent More precisely,

\medskip{}

\hfill{}$t=0\Longrightarrow\left(x,y\right)=\left(1,y\right),\;\mathrm{or\;}\left(-1,y\right),\;\forall y\in\mathbb{R},$\hfill{}(7.4)

\medskip{}

\noindent meanwhile,

\medskip{}

\hfill{}$t>0\Longrightarrow\left(x,y\right)\in\left\{ \left(x,\pm\sqrt{\frac{x^{2}+1}{t}}\right),\left(x,\pm\sqrt{\frac{x^{2}-1}{t}}\right)\right\} $.\hfill{}(7.5) 
\end{theorem}

\begin{proof}
Assume that $t\geq0$. Then $\left\Vert \left(x,y\right)\right\Vert _{1}=1$,
if and only if either 
\[
\left\langle \left(x,y\right),\left(x,y\right)\right\rangle _{t}=x^{2}-ty^{2}=1,
\]
or 
\[
\left\langle \left(x,y\right),\left(x,y\right)\right\rangle _{t}=x^{2}-ty^{2}=-1,
\]
by (7.1), if and only if 
\[
x^{2}-ty^{2}=\pm1,\;\;\textrm{notationally.}
\]
Therefore, the general characterization (7.3) holds.

Suppose $t=0$. Then the above characterization (7.3) becomes that
\[
x^{2}-0\cdot y^{2}=\pm1\Longleftrightarrow x^{2}=1,
\]
if and only if $x=\pm1$ in $\mathbb{R}$. Therefore, the sub-characterization
(7.4) of (7.3) holds, where $t=0$.

Now, let $t>0$. Then the characterization (7.3) says that 
\[
x^{2}-ty^{2}=\pm1\Longleftrightarrow y^{2}=\frac{x^{2}\pm1}{t},
\]
if and only if 
\[
y=\pm\sqrt{\frac{x^{2}\pm1}{t}},\;\;\;\mathrm{notationally}\;\mathrm{in\;\;\;}\mathbb{R}.
\]
Therefore, the sub-characterization (7.5) of (7.3) holds, where $t>0$. 
\end{proof}
The above results fully characterizes the units of $\left\{ \mathbb{D}_{t}\right\} _{t\geq0}$. 
\begin{definition}
The subset $\mathbb{T}_{t}=\left\{ h\in\mathbb{D}_{t}:\left\Vert h\right\Vert _{t}=1\right\} $
is called the ($t$-scaled) unit subset of the $t$-scaled hyperbolic
ring $\mathbb{D}_{t}$. 
\end{definition}

The following corollary is a direct consequence of the above unit
characterizations (7.2), (7.4) and (7.5) on $\left\{ \mathbb{D}_{t}\right\} _{t\in\mathbb{R}}$. 
\begin{corollary}
Let $\mathbb{T}_{t}$ be the unit subset of the $t$-scaled hyperbolic
ring $\mathbb{D}_{t}$.

\medskip{}

\noindent (1) If $t<0$, then $\mathbb{T}_{t}=\left\{ \left(x,\:\pm\sqrt{\frac{1-x^{2}}{\left|t\right|}}\right):x\in\mathbb{R}\right\} $.

\medskip{}

\noindent (2) If $t=0$, then $\mathbb{T}_{t}=\left\{ \left(\pm1,\:y\right):y\in\mathbb{R}\right\} $.

\medskip{}

\noindent (3) If $t>0$, then $\mathbb{T}_{t}=\left\{ \left(x,\:\pm\sqrt{\frac{x^{2}\pm1}{t}}\right):x\in\mathbb{R}\right\} $. 
\end{corollary}

\begin{proof}
The proof of the first statement is done by (7.2). The second statement
is shown by (7.4), and the third one is proven by (7.5). 
\end{proof}
By the statement (1) of the above corollary, the unit subset $\mathbb{T}_{t}$
is an elliptic-like figure in $\mathbb{R}^{2}\overset{\textrm{set}}{=}\mathbb{D}_{t}$,
whenever $t<0$. However, if $t>0$, then the unit subset $\mathbb{T}_{t}$
is like the hyperbolic-like figure in $2$-dimensional $\mathbb{R}$-space
$\mathbb{R}^{2}\overset{\textrm{set}}{=}\mathbb{D}_{t}$, by the statement
(3). Interestingly, the unit subset $\mathbb{T}_{0}$ forms two vertical
lines in $\mathbb{R}^{2}\overset{\textrm{set}}{=}\mathbb{D}_{0}$.
Note that 
\[
\mathbb{T}_{-1}=\left\{ \left(x,\pm\sqrt{1-x^{2}}\right):-1\leq x\leq1\right\} 
\]
is equivalent to the unit circle $\mathbb{T}$ in $\mathbb{R}^{2}\overset{\textrm{set}}{=}\mathbb{C}$,
geometrically, and 
\[
\mathbb{T}_{1}=\left\{ \left(x,\:\pm\sqrt{x^{2}\pm1}\right):x\in\mathbb{R}\right\} 
\]
is equivalent to the unit hyperbola $\mathcal{T}$ in $\mathbb{R}^{2}\overset{\textrm{set}}{=}\mathcal{D}\overset{\textrm{set}}{=}\mathbb{D}_{1}$,
geometrically. These considerations illustrate, or visualize the differences
between the analyses on $\left\{ \mathbb{D}_{t}\right\} _{t\in\mathbb{R}\setminus\left\{ \pm1\right\} }$
and those on the complex analysis (on $\mathbb{D}_{-1}$), and the
hyperbolic analysis (on $\mathbb{D}_{1}$).

In the rest of this section, we consider how to express our $t$-scaled
hyperbolic numbers of $\mathbb{D}_{t}$ as in the complex numbers
of $\mathbb{C}=\mathbb{D}_{-1}$, and the hyperbolic numbers of $\mathcal{D}=\mathbb{D}_{1}$.
For an arbitrary sclae $\mathit{t}$$\in\mathbb{R}$, let $\left(x,y\right)\in\mathbb{D}_{t}$,
realized to be 
\[
\left[\left(x,y\right)\right]_{t}=\left(\begin{array}{cc}
x & ty\\
y & x
\end{array}\right),\;\;\mathrm{in\;\;}\mathcal{D}_{2}^{t}.
\]
Then this realization $\left[\left(x,y\right)\right]_{t}$ is understood
to be 
\[
\left[\left(x,y\right)\right]_{t}=x\left(\begin{array}{cc}
1 & 0\\
0 & 1
\end{array}\right)+y\left(\begin{array}{cc}
0 & t\\
1 & 0
\end{array}\right),\;\mathrm{in\;}\mathcal{D}_{2}^{t},
\]
with\hfill{}(7.6) 
\[
\left(\begin{array}{cc}
1 & 0\\
0 & 1
\end{array}\right)=\left[\left(1,0\right)\right]_{t},\;\mathrm{and\;}\left(\begin{array}{cc}
0 & t\\
1 & 0
\end{array}\right)=\left[\left(0,1\right)\right]_{t},
\]
in $\mathcal{D}_{2}^{t}$. Let $j_{t}=\left(0,1\right)\in\mathbb{D}_{t}$,
satisfying

\medskip{}

\hfill{}$\left[\left(x,y\right)\right]_{t}=x\left[\left(1,0\right)\right]_{t}+y\left[j_{t}\right]_{t}=\left[x+yj_{t}\right]_{t},$\hfill{}(7.7)

\medskip{}

\noindent in $\mathcal{D}_{2}^{t}$, by (7.6). Thus, our $t$-scaled
hyperbolic ring $\mathbb{D}_{t}$ is re-defined to be a set,

\medskip{}

\hfill{}$\mathbb{D}_{t}=\left\{ x+yj_{t}:x,y\in\mathbb{R},\;\mathrm{and\;}j_{t}^{2}=t\right\} ,$\hfill{}(7.8)

\medskip{}

\noindent since 
\[
\left[j_{t}^{2}\right]_{t}=\left[j_{t}\right]_{t}^{2}=\left(\begin{array}{cc}
t & 0\\
0 & t
\end{array}\right)=\left[\left(t,0\right)\right]_{t}=\left[t+0\cdot j_{t}\right]_{t}=\left[t\right]_{t},
\]
in $\mathcal{D}_{2}^{t}$, by (7.6) and (7.7). From below, we consider
the $t$-scaled hyperbolic ring $\mathbb{D}_{t}$ as a set (7.8).
Note that 
\[
j_{t}^{2n}=\left(j_{t}^{2}\right)^{n}=t^{n},\;\mathrm{and\;}j_{t}^{2n-1}=j_{t}^{2n-2}j_{t}=t^{n-1}j_{t},
\]
for all $n\in\mathbb{N}$, by (7.7) and (7.8), i.e.,

\medskip{}

\hfill{}$j_{t}^{n}=\left\{ \begin{array}{ccc}
t^{\frac{n}{2}} &  & \mathrm{if\;}n\;\mathrm{is\;even}\\
\\
t^{\frac{n-1}{2}}j_{t} &  & \mathrm{if\;}n\;\mathrm{is\;odd,}
\end{array}\right.$\hfill{}(7.9)

\medskip{}

\noindent in the hyperbolic ring $\mathbb{D}_{t}$ of (7.8), for all
$n\in\mathbb{N}$.

Now, for any $\theta\in\mathbb{R}$, define an element $e^{j_{t}\theta}$
by

\medskip{}

\hfill{}$e^{j_{t}\theta}\overset{\textrm{def}}{=}\overset{\infty}{\underset{n=0}{\sum}}\frac{\left(j_{t}\theta\right)^{n}}{n!}\in\mathbb{D}_{t}$\hfill{}(7.10)

\medskip{}

\noindent Then one can verify that 
\[
e^{j_{t}\theta}=\left(\overset{\infty}{\underset{n=0}{\sum}}\frac{j_{t}^{2n-1}\theta^{2n-1}}{(2n-1)!}\right)+\left(\overset{\infty}{\underset{n=0}{\sum}}\frac{j_{t}^{2n}\theta^{2n}}{\left(2n\right)!}\right),
\]
implying that

\medskip{}

\hfill{}$e^{j_{t}\theta}=j_{t}\left(\overset{\infty}{\underset{n=0}{\sum}}\frac{t^{n-1}\theta^{2n-1}}{\left(2n-1\right)!}\right)+\left(\overset{\infty}{\underset{n=0}{\sum}}\frac{t^{n}\theta^{2n}}{\left(2n\right)!}\right),$\hfill{}(7.11)

\medskip{}

\noindent in $\mathbb{D}_{t}$, by (7.9). 
\begin{proposition}
Let $e^{j_{t}\theta}\in\mathbb{D}_{t}$ be a $t$-scaled hyperbolic
number (7.10), for $\theta\in\mathbb{R}$. Then 
\[
e^{j_{t}\theta}=\left(\overset{\infty}{\underset{n=0}{\sum}}\frac{t^{n}\theta^{2n}}{\left(2n\right)!}\right)+j_{t}\left(\overset{\infty}{\underset{n=0}{\sum}}\frac{t^{n-1}\theta^{2n-1}}{\left(2n-1\right)!}\right),
\]
i.e.,\hfill{}(7.12) 
\[
e^{j_{t}\theta}=\left(\overset{\infty}{\underset{n=0}{\sum}}\frac{t^{n}\theta^{2n}}{\left(2n\right)!},\;\;\overset{\infty}{\underset{n=0}{\sum}}\frac{t^{n-1}\theta^{2n-1}}{\left(2n-1\right)!}\right)\in\mathbb{D}_{t},
\]
as a $t$-scaled hypercomplex number in $\mathbb{H}_{t}$. 
\end{proposition}

\begin{proof}
The relation (7.12) is shown by (7.8), (7.9) and (7.11). 
\end{proof}
Let $e^{j_{t}\theta}\in\mathbb{D}_{t}$ be a hyperbolic number (7.10),
and assume that a given scale $t$ is positive in $\mathbb{R}$, i.e.,
$t>0$. Then, by (7.12), one has

\medskip{}

\noindent 
\[
e^{j_{t}\theta}=j_{t}\left(\overset{\infty}{\underset{n=0}{\sum}}\frac{\sqrt{t}^{2n-2}\theta^{2n-1}}{\left(2n-1\right)!}\right)+\left(\overset{\infty}{\underset{n=0}{\sum}}\frac{\sqrt{t}^{2n}\theta^{2n}}{\left(2n\right)!}\right),
\]
and hence, 
\[
e^{j_{t}\theta}=\frac{j_{t}}{\sqrt{t}}\left(\overset{\infty}{\underset{n=0}{\sum}}\frac{\left(\sqrt{t}\theta\right)^{2n-1}}{\left(2n-1\right)!}\right)+\left(\overset{\infty}{\underset{n=0}{\sum}}\frac{\left(\sqrt{t}\theta\right)^{2n}}{\left(2n\right)!}\right),
\]
therefore, 
\[
e^{j_{t}\theta}=\left(cosh\left(\sqrt{t}\theta\right)\right)+j_{t}\left(\frac{1}{\sqrt{t}}sinh\left(\sqrt{t}\theta\right)\right),
\]
as a hypercomplex number, 
\[
\left(cosh\left(\sqrt{t}\theta\right),\;\frac{sinh\left(\sqrt{t}\theta\right)}{\sqrt{t}}\right)\in\mathbb{D}_{t},\;\mathrm{in\;}\mathbb{H}_{t},
\]
where\hfill{}(7.13) 
\[
cosh\left(x\right)=\frac{e^{x}+e^{-x}}{2},\;\mathrm{and\;}sinh\left(x\right)=\frac{e^{x}-e^{-x}}{2},
\]
for all $x\in\mathbb{R}$. Under (7.13), observe that

\medskip{}

$\;\;$$\left\Vert e^{j_{t}\theta}\right\Vert ^{2}=\left|\left\langle e^{j_{t}\theta},e^{j_{t}\theta}\right\rangle _{t}\right|$

\medskip{}

$\;\;\;\;\;\;$$=\left|\left\langle \left(cosh\left(\sqrt{t}\theta\right),\;\frac{sinh\left(\sqrt{t}\theta\right)}{\sqrt{t}}\right),\:\left(cosh\left(\sqrt{t}\theta\right),\;\frac{sinh\left(\sqrt{t}\theta\right)}{\sqrt{t}}\right)\right\rangle _{t}\right|$

\medskip{}

$\;\;\;\;\;\;$$=\left|\tau\left(\left(cosh\left(\sqrt{t}\theta\right),\;\frac{sinh\left(\sqrt{t}\theta\right)}{\sqrt{t}}\right)\cdot_{t}\left(cosh\left(\sqrt{t}\theta\right),\;\frac{sinh\left(\sqrt{t}\theta\right)}{\sqrt{t}}\right)^{\dagger}\right)\right|$

\medskip{}

$\;\;\;\;\;\;$$=\left|\tau\left(\left(\begin{array}{ccc}
cosh\left(\sqrt{t}\theta\right) &  & t\left(\frac{sinh\left(\sqrt{t}\theta\right)}{\sqrt{t}}\right)\\
\\
\frac{sinh\left(\sqrt{t}\theta\right)}{\sqrt{t}} &  & cosh\left(\sqrt{t}\theta\right)
\end{array}\right)\left(\begin{array}{ccc}
cosh\left(\sqrt{t}\theta\right) &  & t\left(-\frac{sinh\left(\sqrt{t}\theta\right)}{\sqrt{t}}\right)\\
\\
-\frac{sinh\left(\sqrt{t}\theta\right)}{\sqrt{t}} &  & cosh\left(\sqrt{t}\theta\right)
\end{array}\right)\right)\right|$

\medskip{}

$\;\;\;\;\;\;$$=\left|\tau\left(\left(\begin{array}{ccc}
cosh^{2}\left(\sqrt{t}\theta\right)-sinh^{2}\left(\sqrt{t}\theta\right) &  & 0\\
\\
0 &  & cosh^{2}\left(\sqrt{t}\theta\right)-sinh^{2}\left(\sqrt{t}\theta\right)
\end{array}\right)\right)\right|$

\medskip{}

$\;\;\;\;\;\;$$=\left|cosh^{2}\left(\sqrt{t}\theta\right)-sinh^{2}\left(\sqrt{t}\theta\right)\right|=1$.\hfill{}(7.14)

\medskip{}

\noindent So, the norm (7.13) of $e^{j_{t}\theta}$ is obtained for
any $\theta\in\mathbb{R}$, for a fixed scale $t\in\mathbb{R}$. 
\begin{theorem}
Let $t>0$ in $\mathbb{R}$, and let $e^{j_{t}\theta}\in\mathbb{D}_{t}$
be a $t$-scaled hyperbolic number (7.10). Then 
\[
e^{j_{t}\theta}=cosh\left(\sqrt{t}\theta\right)+j_{t}\left(\frac{sinh\left(\sqrt{t}\theta\right)}{\sqrt{t}}\right),
\]
with\hfill{}(7.15) 
\[
\left\Vert e^{j_{t}\theta}\right\Vert =1\Longleftrightarrow e^{j_{t}\theta}\in\mathbb{T}_{t}.
\]
In particular, if $x+yj_{t}\in\mathbb{D}_{t}$, then

\medskip{}

\hfill{}$x+yj_{t}=\left\Vert \left(x,y\right)\right\Vert _{t}e^{j_{t}\theta},\;\mathrm{for\;}\theta=Arg\left(\left(x,y\right)\right)\in\left[0,2\pi\right].$\hfill{}(7.16) 
\end{theorem}

\begin{proof}
The relation (7.15) is obtained by (7.13) and (7.14). The decomposition
(7.16) is obtained by (7.5) and (7.15). 
\end{proof}
Let $e^{j_{t}\theta}\in\mathbb{D}_{t}$ be a $t$-scaled hyperbolic
number (7.10), and assume that a given scale $t$ is negative in $\mathbb{R}$,
i.e., $t<0$. Then, by (7.12), one has

\medskip{}

\noindent 
\[
e^{j_{t}\theta}=\left(\overset{\infty}{\underset{n=0}{\sum}}\frac{\left(-\left|t\right|\right)^{n}\theta^{2n}}{\left(2n\right)!}\right)+j_{t}\left(\overset{\infty}{\underset{n=0}{\sum}}\frac{\left(-\left|t\right|\right)^{n-1}\theta^{2n-1}}{\left(2n-1\right)!}\right),
\]
and hence, 
\[
e^{j_{t}\theta}=\left(\overset{\infty}{\underset{n=0}{\sum}}\frac{(-1)^{n}\left|t\right|^{n}\theta^{2n}}{\left(2n\right)!}\right)+j_{t}\left(\overset{\infty}{\underset{n=0}{\sum}}\frac{(-1)^{n-1}\left(\left|t\right|\right)^{n-1}\theta^{2n-1}}{\left(2n-1\right)!}\right),
\]
implying that 
\[
e^{j_{t}\theta}=\left(\overset{\infty}{\underset{n=0}{\sum}}\frac{(-1)^{n}\left(\sqrt{\left|t\right|}\theta\right)^{2n}}{\left(2n\right)!}\right)+\frac{j_{t}}{\sqrt{\left|t\right|}}\left(\overset{\infty}{\underset{n=0}{\sum}}\frac{\left(-1\right)^{n-1}\left(\sqrt{\left|t\right|}\theta\right)^{2n-1}}{\left(2n-1\right)!}\right),
\]
therefore, 
\[
e^{j_{t}\theta}=\left(cos\left(\sqrt{\left|t\right|}\theta\right)\right)+j_{t}\left(\frac{1}{\sqrt{\left|t\right|}}sin\left(\sqrt{\left|t\right|}\theta\right)\right),
\]
as a hypercomplex number,

\medskip{}

\hfill{}$\left(cos\left(\sqrt{\left|t\right|}\theta\right),\;\frac{sin\left(\sqrt{\left|t\right|}\theta\right)}{\sqrt{\left|t\right|}}\right)\in\mathbb{D}_{t},\;\mathrm{in\;}\mathbb{H}_{t},$\hfill{}(7.17)

\medskip{}

\noindent Under (7.17), observe that

\medskip{}

$\;\;$$\left\Vert e^{j_{t}\theta}\right\Vert ^{2}=\left|\left\langle e^{j_{t}\theta},e^{j_{t}\theta}\right\rangle _{t}\right|$

\medskip{}

$\;\;\;\;\;\;$$=\left|\left\langle \left(cos\left(\sqrt{\left|t\right|}\theta\right),\;\frac{sin\left(\sqrt{\left|t\right|}\theta\right)}{\sqrt{\left|t\right|}}\right),\:\left(cos\left(\sqrt{\left|t\right|}\theta\right),\;\frac{sin\left(\sqrt{\left|t\right|}\theta\right)}{\sqrt{\left|t\right|}}\right)\right\rangle _{t}\right|$

\medskip{}

$\;\;\;\;\;\;$$=\left|\tau\left(\left(cos\left(\sqrt{\left|t\right|}\theta\right),\;\frac{sin\left(\sqrt{\left|t\right|}\theta\right)}{\sqrt{\left|t\right|}}\right)\cdot_{t}\left(cos\left(\sqrt{\left|t\right|}\theta\right),\;\frac{sin\left(\sqrt{\left|t\right|}\theta\right)}{\sqrt{\left|t\right|}}\right)^{\dagger}\right)\right|$

\medskip{}

$\;\;\;\;\;\;$$=\left|\tau\left(\left(\begin{array}{ccc}
cos\left(\sqrt{\left|t\right|}\theta\right) &  & t\left(\frac{sin\left(\sqrt{\left|t\right|}\theta\right)}{\sqrt{\left|t\right|}}\right)\\
\\
\frac{sin\left(\sqrt{\left|t\right|}\theta\right)}{\sqrt{\left|t\right|}} &  & cos\left(\sqrt{\left|t\right|}\theta\right)
\end{array}\right)\left(\begin{array}{ccc}
cos\left(\sqrt{\left|t\right|}\theta\right) &  & t\left(-\frac{sin\left(\sqrt{\left|t\right|}\theta\right)}{\sqrt{\left|t\right|}}\right)\\
\\
-\frac{sin\left(\sqrt{\left|t\right|}\theta\right)}{\sqrt{\left|t\right|}} &  & cos\left(\sqrt{\left|t\right|}\theta\right)
\end{array}\right)\right)\right|$

\medskip{}

$\;\;\;\;\;\;$$=\left|\tau\left(\left(\begin{array}{ccc}
cos^{2}\left(\sqrt{\left|t\right|}\theta\right)+sin{}^{2}\left(\sqrt{\left|t\right|}\theta\right) &  & 0\\
\\
0 &  & cos^{2}\left(\sqrt{\left|t\right|}\theta\right)+sin{}^{2}\left(\sqrt{\left|t\right|}\theta\right)
\end{array}\right)\right)\right|$

\medskip{}

$\;\;\;\;\;\;$$=\left|cos^{2}\left(\sqrt{\left|t\right|}\theta\right)+sin{}^{2}\left(\sqrt{\left|t\right|}\theta\right)\right|=1$.\hfill{}(7.18)

\medskip{}

\begin{theorem}
Let $t<0$ in $\mathbb{R}$, and let $e^{j_{t}\theta}\in\mathbb{D}_{t}$
be a $t$-scaled hyperbolic number (7.10). Then 
\[
e^{j_{t}\theta}=cos\left(\sqrt{\left|t\right|}\theta\right)+j_{t}\left(\frac{sinh\left(\sqrt{\left|t\right|}\theta\right)}{\sqrt{\left|t\right|}}\right),
\]
with\hfill{}(7.19) 
\[
\left\Vert e^{j_{t}\theta}\right\Vert =1\Longleftrightarrow e^{j_{t}\theta}\in\mathbb{T}_{t}.
\]
In particular, if $x+yj_{t}\in\mathbb{D}_{t}$, then

\medskip{}

\hfill{}$x+yj_{t}=\left\Vert \left(x,y\right)\right\Vert _{t}e^{j_{t}\theta},\;\mathrm{for\;}\theta=Arg\left(\left(x,y\right)\right)\in\left[0,2\pi\right].$\hfill{}(7.20) 
\end{theorem}

\begin{proof}
The relation (7.19) is obtained by (7.17) and (7.18). The decomposition
(7.20) is obtained by (7.2) and (7.19). 
\end{proof}

\section{On the Scaled Hyperbolic Rings}

In this section, we act the matricial $\mathbb{R}$-algebra $M_{2}\left(\mathbb{R}\right)$
on our scaled-hyperbolic spaces, 
\[
\left\{ \mathbf{Y}_{t}=\left(\mathbb{D}_{t},\left\langle ,\right\rangle _{t}\right):t\in\mathbb{R}\right\} ,
\]
in a natural manner. Since all scaled hyperbolic rings $\left\{ \mathbb{D}_{t}\right\} _{t\in\mathbb{R}}$
are equipotent to the $2$-dimensional $\mathbb{R}$-vector space
$\mathbb{R}^{2}$ over $\mathbb{R}$, set-theoretically, indeed, one
can act the matricial algebra $M_{2}\left(\mathbb{R}\right)$ on $\left\{ \mathbb{D}_{t}\right\} _{t\in\mathbb{R}}$,
by regarding them as $\mathbb{R}^{2}$. i.e., if $\left(x,y\right)\in\mathbf{Y}_{t}$,
then

\medskip{}

\hfill{}$\left(\begin{array}{cc}
a_{11} & a_{12}\\
a_{21} & a_{22}
\end{array}\right)\left(\begin{array}{cc}
x\\
y
\end{array}\right)=\left(\begin{array}{ccc}
a_{11}x+a_{12}y\\
\\
a_{21}x+a_{22}y
\end{array}\right),$\hfill{}(8.1)

\medskip{}

\noindent in $\mathbf{Y}_{t}$, as usual, for all $\left[a_{ij}\right]_{2\times2}\in M_{2}\left(\mathbb{R}\right)$.
Note, however, that the resulted vector in (8.1) is contained ``in
$\mathbf{Y}_{t}$,'' realized to be

\medskip{}

\hfill{}$\left(\begin{array}{ccc}
a_{11}x+a_{12}y &  & t\left(a_{21}x+a_{22}y\right)\\
\\
a_{21}x+a_{22}y &  & a_{11}x+a_{12}y
\end{array}\right),$\hfill{}(8.2)

\medskip{}

\noindent in the $t$-hyperbolic realization $\mathcal{D}_{2}^{t}$
in the $t$-hypercomplex realization $\mathcal{H}_{2}^{t}$. i.e.,
for any $A\in M_{2}\left(\mathbb{R}\right)$, one obtains a bounded
operator of the operator space $B\left(\mathbf{Y}_{t}\right)$ (over
$\mathbb{R}$), 
\[
A:\mathbf{Y}_{t}\rightarrow\mathbf{Y}_{t}.
\]
Remark here that the operator space $B\left(\mathbf{Y}_{t}\right)$
is of course a closed subspace of the operator space $B\left(\mathbf{X}_{t}\right)$
in the $t$-hypercomplex space $\mathbf{X}_{t}\in\left\{ \mathbf{H}_{t},\mathbf{K}_{t}\right\} $,
for $t\in\mathbb{R}$.

In the following, fix $t\in\mathbb{R}$, and the $t$-scaled hyperbolic
ring $\mathbb{D}_{t}$, inducing the hyperbolic space $\mathbf{Y}_{t}$
over $\mathbb{R}$. Recall that a hypercomplex number $\left(a,b\right)\in\mathbb{H}_{t}$
is invertible, if and only if 
\[
det\left(\left[\left(a,b\right)\right]_{t}\right)=\left|a\right|^{2}-t\left|b\right|^{2}\neq0,
\]
by (2.2.1) and (2.2.2). So, one can verify that a hyperbolic number
$\left(x,y\right)\in\mathbb{D}_{t}$ is invertible, if and only if

\medskip{}

\hfill{}$det\left(\left[\left(x,y\right)\right]_{t}\right)=x^{2}-ty^{2}\neq0.$\hfill{}(8.3) 
\begin{lemma}
Let $\left(x,y\right)\in\mathbf{Y}_{t}$, and $A=\left[a_{ij}\right]_{2\times2}\in M_{2}\left(\mathbb{R}\right)$.
Then $A\left(\left(x,y\right)\right)\in\mathbf{Y}_{t}$ is invertible
as a hyperbolic number of $\mathbb{D}_{t}$, if and only if

\medskip{}

\hfill{}$\left(a_{11}x+a_{12}y\right)^{2}\neq t\left(a_{21}x+a_{22}y\right)^{2}.$\hfill{}(8.4) 
\end{lemma}

\begin{proof}
By (8.1), the resulted vector $v=A\left(\left(x,y\right)\right)$
in $\mathbf{Y}_{t}$ is realized to be the matrix (8.2) of the $t$-hyperbolic
realization $\mathcal{D}_{2}^{t}$, by regarding this vector $v\in\mathbf{Y}_{t}$
as a hyperbolic number of $\mathbb{D}_{t}$. So, 
\[
det\left(\left[v\right]_{t}\right)=\left(a_{11}x+a_{12}y\right)^{2}-t\left(a_{21}x+a_{22}y\right)^{2},
\]
by (8.3). Therefore, this determinant is nonzero, equivalently, the
condition (8.4) holds, if and only if $v\in\mathbf{Y}_{t}$ is invertible
in $\mathbb{D}_{t}$. 
\end{proof}
By (8.3), we have that a hyperbolic number $\left(x,y\right)$ is
not invertible in $\mathbb{D}_{t}$, if and only if 
\[
x^{2}=ty^{2},\;\;\;\mathrm{in\;\;\;}\mathbb{R},
\]
if and only if\hfill{}(8.5) 
\[
x=\pm\sqrt{t}\left|y\right|\;\;\;\textrm{"in}\;\mathbb{C}\textrm{."}
\]
Remark that if $t<0$, then the $\sqrt{t}=i\sqrt{\left|t\right|}$,
and hence, the quantities $\pm\sqrt{t}\left|y\right|$ become $\pm i\sqrt{\left|t\right|}\left|y\right|$,
for $x$, in (8.5), i.e.,

\medskip{}

\hfill{}$t<0\Longrightarrow\left[\textrm{(9.5) holde}\Longleftrightarrow x=\pm i\sqrt{\left|t\right|}\left|y\right|,\;\mathrm{in\;}\mathbb{C}\setminus\mathbb{R},\right]$\hfill{}(8.6)

\medskip{}

\noindent meanwhile, if $t\geq0$, then the above quantity $\pm\sqrt{t}\left|y\right|$
in (8.5) are $\mathbb{R}$-quantities in $\mathbb{R}$, i.e.,

\medskip{}

\hfill{}$t\geq0\Longrightarrow\left[\textrm{(9.5) holds}\Longleftrightarrow x=\pm\sqrt{t}\left|y\right|,\;\mathrm{in\;}\mathbb{R}.\right]$\hfill{}(8.7)

\medskip{}

\noindent So, the singularity characterization (8.5) is refined by
(8.6) and (8.7), case-by-case. Now, concentrate on the case (8.6).
By the very definition of hyperbolic numbers $\left(x,y\right)\in\mathbb{D}_{t}$,
the entries $x$ and $y$ are contained in $\mathbb{R}$. Therefore,
$x$ cannot be pure-imaginary in $\mathbb{C}$. Therefore, this sub-characterization
(8.6) can be re-stated by

\medskip{}

\hfill{}$t<0\Longrightarrow\left[\left(x,y\right)\textrm{ is not invertible in }\mathbb{D}_{t},\Longleftrightarrow\left(x,y\right)=\left(0,0\right).\right]$\hfill{}(8.8)

\medskip{}

\noindent In fact, the above sub-characterization (8.8) of (8.5) is
verified by (2.2.3), since $\mathbb{D}_{t}$ is a subring of $\mathbb{H}_{t}$. 
\begin{lemma}
A hyperbolic number $\left(x,y\right)$ is not invertible in $\mathbb{D}_{t}$,
if and only if either the condition (8.9), or (8.10) holds, where

\medskip{}

\hfill{}$t<0\Longrightarrow\left(x,y\right)=\left(0,0\right)\;\mathrm{in\;}\mathbb{D}_{t},$\hfill{}(8.9)

\noindent and

\hfill{}$t\geq0\Longrightarrow x=\pm\sqrt{t}\left|y\right|\;\;\mathrm{in\;\;}\mathbb{R}.$\hfill{}(8.10) 
\end{lemma}

\begin{proof}
Assume that $t<0$. Then $\left(x,y\right)\in\mathbb{D}_{t}$ is not
invertible, if and only if the condition (8.5) holds, if and only
if $\left(x,y\right)=\left(0,0\right)$, by (8.6) and (8.8). Meanwhile,
if $t\geq0$, then $\left(x,y\right)\in\mathbb{D}_{t}$ is not invertible,
if and only if the condition (8.5) holds, if and only if 
\[
x=\pm\sqrt{t}\left|y\right|\;\;\mathrm{in\;\;}\mathbb{R},
\]
by (8.7). 
\end{proof}
By the above lemma, one obtains the following equivalent result of
the invertibility (8.4). 
\begin{theorem}
Let $\left(x,y\right)\in\mathbf{Y}_{t}\overset{\textrm{set}}{=}\mathbb{D}_{t}$
be a ``nonzero'' hyperbolic vector, and $A=\left[a_{ij}\right]_{2\times2}\in M_{2}\left(\mathbb{R}\right)$,
non-zero matrix, and let $v=A\left(\left(x,y\right)\right)\in\mathbf{Y}_{t}$.
Then $v$ is not invertible in $\mathbb{D}_{t}$, if and only if either
the condition (8.11), or (8.12) holds, where

\medskip{}

\hfill{}$t<0\Longrightarrow A\textrm{ is not invertible in }M_{2}\left(\mathbb{R}\right),$\hfill{}(8.11)

\noindent and

\hfill{}$t\geq0\Longrightarrow y=\left(\frac{\sqrt{t}a_{21}-a_{11}}{a_{12}-\sqrt{t}a_{22}}\right)x,\;\mathrm{or\;}y=\left(\frac{-\sqrt{t}a_{21}-a_{11}}{a_{12}+\sqrt{t}a_{22}}\right)x.$\hfill{}(8.12) 
\end{theorem}

\begin{proof}
By (8.2), the hyperbolic number $v$ is identified with 
\[
\left(a_{11}x+a_{12}y,\:a_{21}x+a_{22}y\right)\;\mathrm{in\;}\mathbb{D}_{t}.
\]
Thus, $v$ is not invertible in $\mathbb{D}_{t}$, if and only if
either 
\[
\left(a_{11}x+a_{12}y,\:a_{21}x+a_{22}y\right)=\left(0,0\right),\;\mathrm{if\;}t<0,
\]
by (8.9), or 
\[
a_{11}x+a_{12}y=\pm\sqrt{t}\left|a_{21}x+a_{22}y\right|,\;\mathrm{if\;}t\geq0,
\]
by (8.10).

Suppose first that $t<0$. Then the above first condition is identical
to the linear system, 
\[
\left\{ \begin{array}{c}
a_{11}x+a_{12}y=0\\
\\
a_{21}x+a_{22}y=0,
\end{array}\right.
\]
implying that 
\[
a_{11}a_{22}x=a_{12}a_{21}x,\;\mathrm{or\;}a_{21}a_{12}y=a_{11}a_{22}y.
\]
Since $\left(x,y\right)\in\mathbf{Y}_{t}$ is assumed to be non-zero,
i.e., since either $x\neq0$, or $y\neq0$, the above relation implies
that 
\[
a_{11}a_{22}=a_{12}a_{21}\Longleftrightarrow a_{11}a_{22}-a_{12}a_{21}=0,
\]
if and only if $det\left(A\right)=0$. So, if $t<0$, then $v$ is
not invertible in $\mathbb{D}_{t}$, if and only if $det\left(A\right)=0$,
if and only if $A$ is not invertible in $M_{2}\left(\mathbb{R}\right)$.
Therefore, $v$ is not invertible in $\mathbb{D}_{t}$ if and only
if $A$ is not invertible in $M_{2}\left(\mathbb{R}\right)$, whenever
$t<0$. So, the singularity (8.11) holds for $t<0$.

Assume now that $t\geq0$. Then $v$ is not invertible in $\mathbb{D}_{t}$,
if and only if

\medskip{}

\hfill{}$a_{11}x+a_{12}y=\pm\sqrt{t}\left|a_{21}x+a_{22}y\right|.$\hfill{}(8.13)

\medskip{}

\noindent by (8.10). This condition (8.13) satisfies that: for $t\geq0$,
\[
a_{11}x+a_{12}y=\pm\sqrt{t}\left|a_{21}x+a_{22}y\right|,
\]
$\Longleftrightarrow$ 
\[
a_{11}x+a_{12}y=\sqrt{t}\left|a_{21}x+a_{22}y\right|,
\]
or\hfill{}(8.14) 
\[
a_{11}x+a_{12}y=-\sqrt{t}\left|a_{21}x+a_{22}y\right|.
\]
Take the first formula $a_{11}x+a_{12}y=\sqrt{t}\left|a_{21}x+a_{22}y\right|$
in (8.14). Then it is equivalent to 
\[
a_{11}x+a_{12}y=\sqrt{t}\left(a_{21}x+a_{22}y\right),
\]
if $a_{21}x+a_{22}y\geq0$, while,\hfill{}(8.15) 
\[
a_{11}x+a_{12}y=-\sqrt{t}\left(a_{21}x+a_{22}y\right),
\]
if $a_{21}x+a_{22}y<0$. Observe the first case of (8.15): 
\[
a_{11}x+a_{12}y=\sqrt{t}\left(a_{21}x+a_{22}y\right),
\]
$\Longleftrightarrow$ 
\[
a_{12}y-\sqrt{t}a_{22}y=\sqrt{t}a_{21}x-a_{11}x,
\]
$\Longleftrightarrow$\hfill{}(8.16) 
\[
y=\left(\frac{\sqrt{t}a_{21}-a_{11}}{a_{12}-\sqrt{t}a_{22}}\right)x.
\]
It shows that $\left(x,y\right)\in\mathbf{Y}_{t}$ satisfies the condition
(8.16), if and only if $v=A\left(\left(x,y\right)\right)\in\mathbf{Y}_{t}$
is not invertible in $\mathbb{D}_{t}$, whenever $t\geq0$. So, by
considering all four cases (including (8.16)) from the condition (8.14),
one can conclude that: $\left(x,y\right)\in\mathbf{Y}_{t}$ satisfies
either 
\[
y=\left(\frac{\sqrt{t}a_{21}-a_{11}}{a_{12}-\sqrt{t}a_{22}}\right)x,
\]
or 
\[
y=\left(\frac{-\sqrt{t}a_{21}-a_{11}}{a_{12}+\sqrt{t}a_{22}}\right)x,
\]
if and only if $v=A\left(\left(x,y\right)\right)\in\mathbf{Y}_{t}$
is not invertible as an element of $\mathbb{D}_{t}$. Therefore, the
sub-characterization (8.12) holds, too, where $t\geq0$. 
\end{proof}
Also, one can obtain the following spectral-theoretic deformations
on the $t$-scaled hyperbolic ring $\mathbb{D}_{t}$ under the action
of $M_{2}\left(\mathbb{R}\right)$ on the hyperbolic space $\mathbf{Y}_{t}$.
Recall that if $\left(x,y\right)\in\mathbb{D}_{t}$, then 
\[
\sigma_{t}\left(\left(x,y\right)\right)=x+i\sqrt{-ty^{2}}=\left\{ \begin{array}{ccc}
x+i\sqrt{-t}\left|y\right| &  & \mathrm{if\;}t<0\\
\\
x-\sqrt{t}\left|y\right| &  & \mathrm{if\;}t\geq0,
\end{array}\right.
\]
and\hfill{}(8.17) 
\[
spec\left(\left[\left(x,y\right)\right]_{t}\right)=\left\{ x+i\sqrt{-ty^{2}},\:x-i\sqrt{-ty^{2}}\right\} .
\]

\begin{theorem}
Let $\left(x,y\right)\in\mathbf{Y}_{t}$ be a hypercomplex vector,
and $A=\left[a_{ij}\right]_{2\times2}\in M_{2}\left(\mathbb{R}\right)$,
and let $v=A\left(\left(x,y\right)\right)\in\mathbf{Y}_{t}$. Then
$v\in\mathbf{Y}_{t}$ has its $t$-spectral value, 
\[
\sigma_{t}\left(v\right)=\left(a_{11}x+a_{12}y\right)^{2}+i\sqrt{-t}\left|a_{21}x+a_{22}y\right|
\]
as an element of $\mathbb{D}_{t}$, and its realization $\left[v\right]_{t}\in\mathcal{D}_{2}^{t}$
has its spectrum, 
\[
spec\left(v\right)=\left\{ \sigma_{t}\left(v\right),\:\overline{\sigma_{t}\left(v\right)}\right\} \;\;\mathrm{in\;\;}\mathbb{C}.
\]
\end{theorem}

\begin{proof}
The proof is done by (8.17). 
\end{proof}
The above theorem shows the spectral deformation on $\mathbb{D}_{t}$
(in the $t$-scaled hypercomplex ring $\mathbb{H}_{t}$) under the
action of $M_{2}\left(\mathbb{R}\right)$ on the hyperbolic space
$\mathbf{Y}_{t}$.

\bigskip{}

\noindent $\mathbf{Declaration.}$

\bigskip{}

\noindent $\mathbf{Ethical\;Statement.}$ This submitted paper is
the very original, not published, or submitted to elsewhere in any
form or language, partially, or fully.

\bigskip{}

\noindent $\mathbf{Competing\;Interests.}$ None.

\bigskip{}

\noindent $\mathbf{Authors'\;contribution.}$ The co-authors, Daniel
Alpay and Ilwoo Cho, contributed to this manuscript equally, in theoretical
and structural points.
%The written form of the manuscript is ``typed''
%majorly by the corresponding author, the second-named author, Ilwoo
%Cho, however, this written form is checked and confirmed by the co-authors
%together in detail.

\bigskip{}

\noindent $\mathbf{Funding.}$ No fund for co-authors.

\bigskip{}

\noindent $\mathbf{Data\;Availability\;Statement.}$ The authors confirm
that no data known is used in the manuscript.

\bigskip{}


\begin{thebibliography}{10}


\bibitem{bibref1-2}D. Alpay and I. Cho, Operators Induced by Certain
Hypercomplex Systems, Opuscula Math., 43, no. 3, (2023) 275 - 333.

\bibitem{bibref1}D. Alpay, M. E. Luna-Elizarraras, M. Shapiro, and
D. Struppa, Gleason's Problem, Rational Functions and Spaces of Left-regular
Functions: The Split-Quaternioin Settings, Israel J. Math., 226, (2018)
319 - 349.




\bibitem{bibref2}I.Cho, and P. E. T. Jorgensen, Multi-Variable Quaternionic
Spectral Analysis, Opuscula Math., 41, no.3, (2021) 335 - 379.

\bibitem{bibref3}M. L. Curtis, Matrix Groups, ISBN: 0-387-90462-X,
(1979) Published by Springer-Verlag, NY.

\bibitem{bibref4}F. O. Farid, Q. Wang, and F. Zhang, On the Eigenvalues
of Quaternion Matrices, Linear \& Multilinear Alg., Issue 4, (2011),
451 - 473.

\bibitem{bibref5}C. Flaut, Eigenvalues and Eigenvectors for the Quaternion
Matrices of Degree Two, An. St. Univ. Ovidius Constanta, vol. 10,
no. 2, (2002) 39 - 44.

\bibitem{bibref6}P. R. Girard, Einstein's Equations and Clifford
Algebra, Adv. Appl. Clifford Alg., 9, no. 2, (1999) 225 - 230.

\bibitem{bibref7}P. R. Halmos, Linear Algebra Problem Book, ISBN:
978-0-88385-322-1, (1995) Published by Math. Assoc. Amer.

\bibitem{bibref8}P. R. Halmos, Hilbert Space Problem Book, ISBN:
978-038-790685-0, (1982) Published by Springer-Verlag (NY).

\bibitem{bibref9}W, R. Hamilton, Lectures on Quaternions, Available
on http://books.google.com, (1853) Published by Cambridge Univ. Press.

\bibitem{bibref10}I. L. Kantor, and A. S. Solodnikov, Hypercomplex
Numbers, an Elementary Introuction to Algebras, ISBN: 0-387-96980-2,
(1989) Published by Springer-VErlag (NY).

\bibitem{bibref11}V. Kravchenko, Applied Quaternionic Analysis, ISBN:
3-88538-228-8, (2003) Published by Heldemann Verlag.

\bibitem{bibref12}M. Kobayashi, Hyperbolic Hopfield Neural Networks
with Directional Multistate Activation Function, Neurocomputing, 275,
(2018) 2217 - 2226.

\bibitem{bibref13}J. Ratcliffe, Foundations of Hyperbolic Manifolds,
Grad. Text in Math., vol. 149, ISBN:978-0-3873319-73, (2006) Published
by Springer.

\bibitem{bibref14}W. F. Reynolds, Hyperbolic Geometry on a Hyperboloids,
Amer. Math. Monthly, 100, (1993) 442 - 455.

\bibitem{bibref15}L. Rodman, Topics in Quaternion Linear Algebra,
ISBN:978-0-691-16185-3, (2014) Published by Prinston Univ. Press,
NJ.

\bibitem{bibref16}B. A. Rozenfeld, The History of non-Eucledean Geometry:
Evolution of the Concept of a Geometric Spaces, ISBN: 978-038-796458-4,
(1988) Published by Springer.

\bibitem{bibref16-1}R. Speicher, Combinatorial Theory of the Free
Product with Amalgamation and Operator-Valued Free Probability Theory,
Amer. Math. Soc., Memoire, ISBN:978-0-8218-0693-7, (1998) Published
by Amer. Math. Soc..

\bibitem{bibref17}A. Sudbery, Quaternionic Analysis, DOI: 10.1017/s0305004100053638,
Math. Proc. Cambridge Philosop. Soc., (1998)

\bibitem{bibref18}L. Taosheng, Eigenvalues and Eigenvectors of Quaternion
Matrices, J. Central China Normal Univ., 29, no. 4, (1995) 407 - 411.

\bibitem{bibref19} J. A. Vince, Geometric Algebra for Computer Graphics,
ISBN: 978-1-84628-996-5, (2008) Published by Springer.

\bibitem{bibref19-1}D.V. Voiculescu, K. J. Dykema, and A. Nica, Free
Random Variables, ISBN:978-0-8218-1140-5, (1992) Published by Amer.
Math. Soc..

\bibitem{bibref20}J. Voight, Quaternion Algebra, Available on http://math.dartmouth.edu/\symbol{126}jvoight/quat-book.pdf,
(2019) Dept. of Math., Dartmouth Univ.. 
\end{thebibliography}
\end{document}